\documentclass[francais]{smfart} 

\usepackage{baskervald}
\usepackage[utf8]{inputenc} 
\usepackage[francais,english]{babel}
\usepackage{textcomp} 
\usepackage{stmaryrd}
\usepackage{graphicx} 
\usepackage{verbatim} 
\usepackage[retainorgcmds]{IEEEtrantools}
\usepackage[nice]{nicefrac}  
\usepackage{url}
\usepackage{tikz}
\usetikzlibrary{matrix,arrows}
\definecolor{cqcqcq}{rgb}{0.752941176471,0.752941176471,0.752941176471}
\definecolor{ffqqqq}{rgb}{0.333333333333,0.333333333333,0.333333333333}
\definecolor{cqcqcq}{rgb}{0.752941176471,0.752941176471,0.752941176471}
\definecolor{qqqqff}{rgb}{0.,0.,1.}
\definecolor{cqcqcq}{rgb}{0.752941176471,0.752941176471,0.752941176471}
\definecolor{ffqqqq}{rgb}{1.,0.,0.}
\usepackage{mathabx,epsfig}

\usepackage[all]{xy}
 \def\dar[#1]{\ar@<2pt>[#1]\ar@<-2pt>[#1]}
 \def\tar[#1]{\ar@<4pt>[#1]\ar@<0pt>[#1]\ar@<-4pt>[#1]}

\usepackage{amsthm} 
\usepackage{amsmath} 
\usepackage{amsfonts}
\usepackage{amssymb}
\usepackage[shortlabels]{enumitem}
\usepackage{todonotes}

\renewcommand\epsilon\varepsilon 
\renewcommand\phi\varphi 
\def\quotient#1#2{
  \raise0ex\hbox{$#1$}\big/\!\lower1ex\hbox{$#2$}
}

\newcommand\NN{\mathbb{N}} 
\newcommand\ZZ{\mathbb{Z}} 
\newcommand\QQ{\mathbb{Q}} 
\newcommand\CC{\mathbb{C}} 
 
\newcommand\FF{\mathbb{F}} 
\newcommand\FP{\mathbb{F}_p}

\newcommand\fleche{\longrightarrow} 
\newcommand{\dans}{ \!\in\! }					
\newcommand{\priv}{ \!\smallsetminus\! }		
\newcommand{\eps}{ \varepsilon}				
\newcommand{\defi}{ \begin{defin} }
\newcommand{\edefi}{ \end{defin} }
\newcommand{\exo}{ \begin{exercice} }		
\newcommand{\eexo}{ \end{exercice} }
\newcommand{\thr}{ \begin{theor} }
\newcommand{\ethr}{ \end{theor} }
\newcommand{\hyp}{ \begin{hypothese} }
\newcommand{\ehyp}{ \end{hypothese} }
\newcommand{\ques}{ \begin{question} }
\newcommand{\eques}{ \end{question} }
\newcommand{\nott}{ \begin{note} }
\newcommand{\enott}{ \end{note} }
\newcommand{\pro}{ \begin{prop} }
\newcommand{\epro}{ \end{prop} }
\newcommand{\prodefi}{ \begin{propdefi} }
\newcommand{\eprodefi}{ \end{propdefi} }
\newcommand{\rem}{ \begin{rema} }
\newcommand{\erem}{ \end{rema} }
\newcommand{\eppte}{ \end{propriete} }
\newcommand{\ppte}{ \begin{propriete} }
\newcommand{\cor}{ \begin{corr} }
\newcommand{\ecor}{ \end{corr} }
\newcommand{\lem}{ \begin{lemm} }
\newcommand{\elem}{ \end{lemm} }
\newcommand{\exen}{ \begin{exemple} }
\newcommand{\eexen}{ \end{exemple} }
\newcommand{\exe}{ \begin{exemple} }
\newcommand{\eexe}{ \end{exemple} }
\newcommand{\rap}{ \begin{rappel} }
\newcommand{\erap}{ \end{rappel} }

\newcommand{\fact}{ \begin{factt} }
\newcommand{\efact}{ \end{factt} }
\newcommand{\conj}{ \begin{conjecture} }
\newcommand{\econj}{ \end{conjecture} }
\newcommand{\dem}{ \begin{proof}[\textsc{Démonstration}] }  
\newcommand{\edem}{ \end{proof} } 

\DeclareMathOperator{\Ker}{Ker}

\DeclareMathOperator{\Coker}{Coker}

\DeclareMathOperator{\rg}{rg}
\DeclareMathOperator{\Ht}{ht}
\DeclareMathOperator{\Hom}{Hom}

\DeclareMathOperator{\Id}{Id}

\DeclareMathOperator{\Spec}{Spec}
\DeclareMathOperator{\Spf}{Spf}

\DeclareMathOperator{\id}{id}

\DeclareMathOperator{\Gal}{Gal}

\DeclareMathOperator{\Fil}{Fil}
\DeclareMathOperator{\im}{Im}

\DeclareMathOperator{\Gr}{Gr}
\DeclareMathOperator{\Ha}{Ha}
\DeclareMathOperator{\Frob}{Frob}

\DeclareMathOperator{\Hdg}{Hdg}
\DeclareMathOperator{\Fitt}{Fitt}
\DeclareMathOperator{\Newt}{\mathcal N\!ewt}

\DeclareMathOperator{\Def}{Def}
\DeclareMathOperator{\Cris}{Cris}
\DeclareMathOperator{\Diag}{Diag}
\DeclareMathOperator{\End}{End}

\DeclareMathOperator{\HN}{HN}
\DeclareMathOperator{\Deg}{Deg}

\theoremstyle{definition} 
\newtheorem{defin}{Définition}[section]
\newtheorem{propriete}[defin]{Propriété}
\newtheorem{hypothese}[defin]{Hypothèse}
\newtheorem{question}[defin]{Question}
\newtheorem{exercice}[defin]{Exercice}

\theoremstyle{plain} 
\newtheorem{theor}[defin]{Théoreme}
\newtheorem{lemm}[defin]{Lemme}  
\newtheorem{prop}[defin]{Proposition}
\newtheorem{propdefi}[defin]{Proposition-Définition}
\newtheorem{corr}[defin]{Corollaire}

\newtheorem{rappel}[defin]{Rappel}
\newtheorem{conjecture}[defin]{Conjecture}

\theoremstyle{remark} 
\newtheorem{rema}[defin]{Remarque}

\newtheorem*{exemple*}{Exemple}
\newtheorem{exemple}[defin]{Exemple}
\newtheorem{note}[defin]{Note}
\newtheorem{factt}[defin]{Fait}

\title{ \textsc{La Filtration canonique des $\mathcal O$-modules} $p$-\textsc{divisibles.}  }
\date{}
\author{Valentin Hernandez}
\address{Bureau 509, Tour 15-16\\4 Place Jussieu\\75005 Paris}
\email{valentin.hernandez@imj-prg.fr}
\urladdr{}

\begin{document}

\frontmatter
\altkeywords{}
\thanks{}

\maketitle

\begin{abstract}
In this article we associate to $G$, a truncated $p$-divisible $\mathcal O$-module of given signature, where $\mathcal O$ is a finite unramified extension of $\ZZ_p$, a filtration of $G$ by sub-$\mathcal O$-modules under the conditions that his Hasse $\mu$-invariant is smaller than an explicite bound. This filtration generalise the one given when $G$ is $\mu$-ordinary.
The construction of the filtration relies on a precise study of the cristalline periods of a $p$-divisible $\mathcal O$-module. We then apply this result to families of such groups, in particular to stricts neighbourhoods of the $\mu$-ordinary locus inside some PEL Shimura varieties.

\end{abstract}

\selectlanguage{french}
\begin{abstract} Dans cet article, à $G$ un groupe $p$-divisible tronqué muni d'une action d'une extension finie non ramifiée $\mathcal O$ de $\ZZ_p$, et de signature donnée, on associe sous une condition explicite sur son $\mu$-invariant de Hasse, une filtration de $G$ par des sous-$\mathcal O$-modules qui étend la filtration canonique lorsque $G$ est $\mu$-ordinaire.
La construction se fait en étudiant les périodes cristallines des groupes $p$-divisibles avec action de $\mathcal O$. On applique ensuite cela aux familles de tels groupes, en particulier des voisinages stricts du lieu $\mu$-ordinaire dans des variétés de Shimura PEL.
\end{abstract}

\tableofcontents

\section{Introduction}

\subsection{Historique et énoncé}
En 1972, à la suite de l'introduction par Serre des formes modulaires $p$-adiques, Katz dans \cite{Katz} introduit la notion de formes modulaires surconvergentes. Ces dernières, cas particulier des formes modulaires $p$-adiques de Serre, sont des sections de fibrés automorphes sur des voisinages stricts du lieu ordinaire. Il est possible d'utiliser les formes surconvergentes pour interpoler les formes modulaires classique, et ceci a des applications considérables en théorie des nombres.

Dans ce même article, Katz, construit un opérateur compact sur l'ensemble des formes surconvergentes en utilisant un théorème de Lubin, sur l'existence d'un sous-groupe canonique 
dans la $p$-torsion d'une courbe elliptique, à condition que celle-ci ait bonne réduction, et que son invariant de Hasse ne soit pas trop grand (cf. \cite{Katz} Theorem 3.1). Cet opérateur 
compact est d'une importance capitale dans la théorie des formes modulaires, et a en particulier mené à la construction par Coleman et Mazur de la "Eigenvariety", \cite{CM}.

La théorie du sous-groupe canonique de Lubin a depuis été généralisée à la torsion des groupes $p$-divisibles (tronqués) dont l'invariant de Hasse n'est pas trop grand, en utilisant les travaux de nombreux auteurs, dont Andreatta-Gasbarri, Abbes-Mokrane, Conrad, Tian, et Fargues, \cite{Far}, qui démontre le théorème suivant,
\thr[Fargues]
Soit $p>3$ un nombre premier, $K/\QQ_p$ une extension, et $G/\Spec(\mathcal O_K)$ un groupe $p$-divisible tronqué d'échelon $r$. Supposons que
\[ \Ha(G) < \frac{1}{2p^{r-1}}.\]
Alors il existe un sous-groupe fini et plat $C \subset G$ de degré $r\dim G - \frac{p^r-1}{p-1}\Ha(G)$ tel que $C(\mathcal O_{\overline K}) = (\ZZ/p^n\ZZ)^{\dim G[p]}$, que modulo $p^{1-\Ha(G)}$, $C$ coincide avec $\Ker F^r$ le noyau de Frobénius iteré $n$-fois. De plus $C$ est un cran de la filtration de Harder-Narasihman de $G$ (\cite{FarHN}) et vérifie donc de nombreuses compatibilités.
\ethr

De plus, à l'aide des filtrations de Harder-Narasihman, le théorème précédent peut se mettre en famille, par exemple sur des voisinages stricts du lieu ordinaire de variété de Shimura.
Le théorème précédent est aussi valide pour $p=3$ avec une borne un peu moins précise, et il a été redemontré par Hattori, qui autorise $p=2$ (et donne une meilleure borne pour $p=3$). Il a été aussi redémontré par Scholze d'une manière différente à l'aide du complexe cotangeant \cite{Sch3}.

Malheureusement le théorème précédent ne s'applique pas sur de nombreuses variétés de Shimura, comme toutes celles dont le lieu $\mu$-ordinaire est vide ; dans ce cas, tous les groupes $p$-divisibles qui apparaissent ont leur invariant de Hasse égal à 1.
L'idée est alors d'utiliser un autre invariant, l'invariant $\mu$-ordinaire tel que par exemple introduit dans \cite{GN}, voir aussi \cite{KW,Box} et \cite{Her1}. 
On utilisera la construction de \cite{Her1} puisqu'elle s'applique à notre situation (voir section \ref{sectmuha}).
Soit alors $\mathcal O$ une extension finie non ramifiée de $\ZZ_p$ de degré $f$, $G$ un groupe $p$-divisible tronqué sur $\mathcal O_K$, $K$ une extension valuée de $\QQ_p$, muni d'une action de $\mathcal O$ qui permet 
d'associer à $G$ une signature $(p_\tau,q_\tau)_{\tau \in \Hom(\mathcal O,\overline \QQ_p)}$. La $\mathcal O$-hauteur de $G$ (définie comme celle de $G[p]$) est alors $h = \Ht_\mathcal O(G) = p_\tau+q_\tau$. On peut alors associer à $G$ son degré (cf. \cite{Far}), mais aussi ses degrés partiels,
\[ \deg_\tau G = \deg \omega_{G,\tau},\]
où l'on écrit $\omega_G = \bigoplus_\tau \omega_{G,\tau}$ et le degré d'un $\mathcal O_K$-module de la forme $\bigoplus_{i=1}^r \mathcal O_K/a_i\mathcal O_K$ est $\sum_i v(a_i)$ (normalisée par $v(p) = 1$). On peut aussi construire des filtrations de Harder-Narasihman pour des schémas en $\mathcal O$-modules, en remplaçant la fonction degré utilisée dans \cite{Far}, par, pour tout $\tau$,
\[ \Deg_\tau(G) = \sum_{i=1}^f p^{f-i} \deg_{\sigma^i\tau}(G),\]
où $\sigma$ est le relèvement à $\mathcal O$ du Frobénius $x \dans \longmapsto x^p$. On note alors $\HN_\tau$ la filtration de Harder-Narasihman construite à partir de la fonction $\Deg_\tau$.
Notre théorème principal est – de manière simplifiée – le suivant,

\thr
On suppose que $p>4h$ et $k = (f-1)\max_\tau q_\tau$, pour simplifier. Pour tout $\mathcal O$-module $p$-divisible tronqué $G$ d'échelon $r + k$ sur $\Spec(\mathcal O_C)$ ($C =\widehat{\overline{\QQ_p}}$) comme précédemment. 
vérifiant que,
\[ {^\mu}\Ha(G) < \frac{1}{2p^{f(r-1)}},\]
alors il existe une filtration de $G[p^r]$ par des sous-$\mathcal O$-modules finis et plats,
\[ \Fil_\tau(G[p^r]) \subset G[p^r],\]
telle que $\Fil_\tau(G[p^r]) \subset \Fil_{\tau'}(G[p^r])$ si et seulement si $p_\tau \leq p_{\tau'}$ et $\Fil_\tau(G[p^r])(\mathcal O_{C}) \simeq (\mathcal O/p^r\mathcal O)^{p_\tau}$. Les degrés des sous-$\mathcal O$-modules de la filtration sont alors controlés par,
\[ \Deg_\tau(\Fil_\tau(G[p^r])) \geq r \sum_{i = 1}^f \min(p_\tau,p_{\sigma^{i}\tau})p^{f-i} - \frac{p^{rf} - 1}{p^f-1}\Ha_\tau(G).\]
En particulier ceux-ci sont "de grand degré". Pour la $p$-torsion, on a la formule exacte,
\[ \Deg_\tau(\Fil_\tau(G[p])) = \sum_{i=1}^f \min(p_\tau,p_{\sigma^{i}\tau})p^{f-i} - \Ha_\tau(G).\]
De plus, chaque $\Fil_\tau(G[p^r])$ est un cran de la filtration de Harder-Narasihman modifiée $\HN_\tau$ de $G[p^r]$ (si $r = 1$ ce sont aussi des crans de la filtration de Harder-Narasihman de $G[p]$ introduite dans \cite{FarHN}), cette filtration vérifie donc toutes les compatibilités classiques.

De plus, si $r = nf$, la combinaison linéaire,
\[K_n = \sum_\tau \Fil_\tau (G[p^f])[p^{nr_\tau}] \subset G[p^{nf}],\] 
où $r_\tau =  |\{\tau' \dans \Hom(\mathcal O, \mathcal O_{C}) : q_{\tau'} \leq q_\tau\}|$, déforme le noyau de Frobenius $F^{nf}$, au sens où $K_n \otimes_{\mathcal O_C} \overline{\FP} = \Ker(F^{nf})$.

\ethr
En fait, dans les théorèmes du texte on a des meilleurs constantes que celles annoncées précedement. 
À cause de l'utilisation de l'invariant de Hasse $\mu$-ordinaire tel que décrit dans \cite{Her1}, on peut en fait prendre 
$k = \max \{ k_\tau : \tau \text{ tels que } q_\tau \neq h\} + 1$ où
\[ k_\tau = \sum_{\tau'} \max(q_\tau - q_{\tau'},0),\]
qui est relié à l'échelon de $G$ et n'est nécéssaire que pour définir l'invariant de Hasse. On peut probablement contourner ce problème (i.e. avoir $k = 0$) et n'utiliser qu'un groupe d'échelon $r$ en modifiant un peu la construction de l'invariant de Hasse. 
L'hypothèse sur $p$ est dans les fait moins restrictive. On doit tout d'abord supposer $p > \max\{q_\tau : q_\tau \neq h\} + 1$ pour appliquer le théorème de Faltings \cite{Fal} – central dans notre construction.
De plus, et c'est un problème technique, on doit supposer que pour tout $\tau$ tel que $q_\tau \not\in \{0,h\}$ $p> \frac{2q_\tau}{1+K_\tau}$ – où $K_\tau < 1$ est une constante 
(explicite) dépendant de la signature et de $\tau$ – pour calculer les "degrés" des sous-groupes $\Fil_\tau(G[p^r])$. 
Sous cette hypothèse sur $p$, on a alors une borne sur ${^\mu}\Ha$ explicite, mais un peu différente de la borne du théorème précédent (un peu plus compliquée). Néanmoins si $p > 4q_\tau$, pour tout $\tau$, cette borne devient celle annoncée, à savoir $\frac{1}{2p^{f(r-1)}}$.
Le théorème précis est le théorème \ref{thrfinO}.

\subsection{Construction de la filtration}

On commence par procéder localement (i.e. sur $\mathcal O_C$, $C = \widehat{\overline{\QQ}}_p$) en étudiant finement l'application de Hodge-Tate $\alpha_G$.
Dans le cas où $G/\Spec(\mathcal O_C)$ est un $\mathcal O$-module $p$-divisible tronqué de signature $(p_\tau,q_\tau)$, on peut décomposer son complexe conormal,
\[\omega_{G^D} = \bigoplus_\tau \omega_{G^D,\tau},\]
où les constituants sont de dimensions $q_\tau$.
On peut alors regarder l'application de Hodge-Tate modifiée,
\[\alpha_{G[p],\tau} : G[p](\mathcal O_C) \overset{\alpha_G[p]}{\fleche} \omega_{G[p]^D} \fleche \omega_{G[p]^D,\tau}.\]
Lorsque $G/\Spec(\mathcal O_C)$ est un $\mathcal O$-module $p$-divisible (non tronqué) $\mu$-ordinaire, on a montré dans \cite{Her1} – en utilisant essentiellement les travaux de \cite{SW}– que $G$ est explicite,
\[ G \simeq \prod_{\ell = 1}^r \mathcal{LT}_{A_\ell}^{q_{\ell +1}-q_\ell},\]
(on renvoie à \cite{Her1} ou à la sous-section \ref{sect31} pour les notations) et on peut donc explicitement calculer la filtration de Harder-Narasihman de $G$, qui est donnée par les sous-groupes,
\[ \Fil_j = \prod_{\ell = j}^r \mathcal{LT}_{A_\ell}^{q^{(\ell +1)}-q^{(\ell)}},\]
et chacun de ces groupes correspond à certains plongements $\tau$ (tous ceux tels que $q_\tau = q^{(j)}$).
On peut alors explicitement calculer l'application de Hodge-Tate modifiée dans ce cas, et remarquer que $T_p\Fil_j = \Ker \alpha_{G,\tau}$ pour tout $\tau$ tel que $q_\tau = q^{(j)}$.

Dans Fargues (i.e. dans le cas d'une signature parallèle) on pouvait alors utiliser les résultats sur l'annulation de la cohomologie de la suite de Hodge-Tate pour en déduire qu'un certain noyau de $\alpha_G$ était un bon candidat pour être le sous-groupe canonique. Néanmoins dans notre cas cette stratégie n'est plus du tout suffisante, puisque l'application $\alpha_G$, même pour un groupe $p$-divisible $\mu$-ordinaire, est loin d'être surjective ! Il faut alors modifier l'application de Hodge-Tate par certaines périodes cristallines pour espérer contrôler son image (et son noyau).
En utilisant les résultats de Fatlings sur les Frobénius-cristaux filtrés, \cite{Fal}, et c'est le coeur technique de la construction, on peut alors prouver sous l'hypothèse sur ${^\mu}\Ha$ que, quitte à réduire un petit peu, l'image de 
$\alpha_{G[p],\tau}$ est de la 
taille espérée, et donc que son noyau définit un sous-$\mathcal O$-module fini et plat de $\mathcal O$-hauteur $p_\tau$, $\Fil_\tau(G[p])$. Les Frobénius cristaux étant des modules 
sur l'anneau $A_{cris}$ de Fontaine, l'énoncé précédent se ramène à des multiplications et divisions habiles par des périodes de modules de Lubin-Tate $t_\mathcal O$ dans 
$A_{cris}$. Une fois que l'on contrôle l'application de Hodge-Tate, on peut alors en déduire, grâce à une méthode semblable à celle de \cite{Far}, une formule sur les degrés partiels de $\Fil_\tau(G[p])$.

On procède ensuite par récurrence pour construire $\Fil_\tau(G[p^r]) \subset G[p^r]$, le point technique étant de prouver que les sous-groupes ainsi construit 
sont effectivement des crans de filtrations de type Harder-Narasihman, la difficulté étant principalement d'ordre combinatoire et il faut controler les degrés de manière plus précise 
que dans \cite{Far}, ce qui rend la preuve plus technique. On utilise alors de manière centrale les propositions élémentaires mais extrêmement astucieuses introduites dans 
\cite{Bij} sur les degrés partiels, que l'on rappelle en annexe.
Notons aussi que Bijakowski dans \cite{Bij} démontre l'existence d'une filtration canonique sur les variétés de Shimura (PEL) en niveau Iwahorique, et la caractérise complètement en terme de "degrés partiels", en utilisant uniquement les propriétés du degré d'un schéma en groupes.

\subsection{Applications aux variétés de Shimura}

De la même manière que dans \cite{Far}, on peut utiliser les filtrations de Harder-Narasihman $\HN_\tau$ pour mettre en famille les filtrations précédentes sur un espace rigide 
$\mathcal X$, en particulier sur (des voisinages stricts du lieu $\mu$-ordinaire) des variétés de Shimura PEL non ramifiées en $p$. En particulier le théorème précédent permet de 
construire un opérateur compact $U_p$ sur les formes surconvergentes des variétés de Shimura (PEL non ramifiée en $p$) sans lieu ordinaire.

Soit $\mathcal D = (B,\star,V,<,>)$ une donnée de Shimura PEL, que l'on suppose non ramifiée en $p$ (voir \cite{VW} section 1.1). On peut alors associer à $\mathcal D$ et à un niveau $K$, un schéma $X_K$ sur $\Spec(O_E)$, où $E/\QQ_p$ est une extension finie, qui est un espace de module de variétés abéliennes (\cite{KotJams}).
Soit $A \fleche X$ la variété abélienne universelle, elle est munie d'une action de $O_B$, un ordre dans l'algèbre $B$. D'après l'hypothèse sur $p$, $\mathcal O_B\otimes_{\ZZ}\ZZ_p$ se scinde en un produit (fini) d'algèbres de matrices sur des entensions non ramifées de $\ZZ_p$ et nous permet d'écrire,
\[A[p^\infty] = \prod_{i=1}^r A[\pi_i^\infty],\]
où chaque $A[\pi^\infty]$ est un groupe $p$-divisible muni d'une action de $M_n(\mathcal O_{K_i})$, $K_i/\QQ_p$ non ramifiée. Par l'équivalence de Morita, on peut alors écrire,
\[A[\pi^\infty] = \mathcal O_{K_i}^n \otimes_{\mathcal O_K} G_i,\]
où $G_i$ est un $\mathcal O_{K_i}$-module $p$-divisible du type considéré précédemment, dont on note $(p^i_\tau,q_\tau^i)$ la signature.
On a alors pour tout $i$, une section d'un fibré en droite, ${^\mu}\Ha(G_i)$ qui détermine un ouvert (dense) de $X_K$, le lieu $\mu$-ordinaire de $G_i$. L'intersection de ces ouverts 
est alors le lieu $\mu$-ordinaire de $X_K$. Supposons maintenant que le groupe réductif $G$ associé à $\mathcal D$ a un modèle réductif $\mathcal G$ sur $\ZZ_p$, notons le niveau $K = K^pK_p$ où $K^p$ est un niveau hors $p$ assez petit, et $K_p \subset \mathcal G(\ZZ_p)$. On note $\mathcal P \subset \mathcal G$ un sous-groupe parabolique, determiné par les signatures $(p_\tau^i,q_\tau^i)$. On peut alors construire une suite décroissante de sous-groupes de congruence, $\mathcal P_0 = K^p\mathcal G(\ZZ_p)$ et $\mathcal P_n = K^p\mathcal P_{n,p}$, où
\[\mathcal P_{n,p} = \{ g \dans \mathcal G(\ZZ_p) : g \pmod p^n \dans \mathcal P(\ZZ/p^n\ZZ)\},\]
On peut alors mettre en famille les résultats précédents (voir théorème \ref{thrfam}), et en déduire,

\thr
Supposons que $p$ est assez grand devant la signature $(p_\tau^i,q_\tau^i)_{\tau,i}$, au sens de la section \ref{sect9}. 
Soit $X_{\mathcal P_k}^{rig}$ la variété rigide associée au schéma $X_{\mathcal P_k}$. Alors pour tout $i$, il existe une constante explicite $\eps_n^i$ tel que sur l'ouvert,
\[ X_{\mathcal P_k}(\eps_n) = \{ x \dans X_{\mathcal P_k}^{rig}: {^\mu}\Ha(G_i)(x) < \eps_n^i\},\]
il existe une filtration de $A[p^n]$ par des sous-groupes finis et plats de grands degrés qui étend la filtration canonique sur le lieu ordinaire.
En particulier, on en déduit une section,
\[ X_{\mathcal P_0}(\eps_n) \overset{s_n}{\fleche} X_{\mathcal P_n}(\eps_n),\]
qui étend la section donnée sur le lieu ordinaire par la filtration canonique.
\ethr

\subsection{$\mathcal O$-modules stricts et applications futures}

Un cas particulier où notre théorème s'applique est le cas des $\mathcal O$-modules stricts de Faltings, où la signature est donnée par $(d,h-d),(0,h),\dots,(0,h)$.
Dans ce cas il y a un seul plongement intéressant (celui pour lequel $p_\tau = d$), et la filtration précédente est réduite à un cran. Le groupe dans le cas $\mu$-ordinaire est un produit 
de modules de Lubin-Tate $\mathcal{LT}_\tau$.

Bien que dans ce cas le théorème de \cite{Far} ne s'applique pas, on pourrait probablement appliquer sa démonstration en utilisant une théorie des cristaux pour les 
$\mathcal O$-modules stricts (i.e. utiliser $V_\pi$ le Verschiebung modifié de \cite{Fal3} à la place du Verschiebung) et – probablement – en déduire une version de notre théorème 
avec une meilleure borne sur $p$ (et peut-être ne pas réellement recourir aux résultats complets de \cite{Fal}).

On espère utiliser nos résultats pour construire des familles de formes modulaires surconvergentes propres pour les opérateurs de Hecke, dans la cohomologie cohérente des variétés de Shimura PEL non ramifiée en $p$. On reviendra sur cette question dans un futur proche.

\subsection{Description de l'article}

Dans la section \ref{sect2} on rappelle les notations classiques, les polygones associés aux $\mathcal O$-modules, ainsi que le cristal associé à de tels groupes.
Dans la section \ref{sect3} on rappelle la construction des invariants de Hasse de \cite{Her1}, ainsi que la structure des $\mathcal O$-modules $p$-divisibles $\mu$-ordinaires. On 
explique ensuite pourquoi ceux-ci sont munis d'une filtration "naturelle".

La section \ref{sect4} est consacrée à l'application de Hodge-Tate. Après l'avoir définie ainsi que sa variante cristalline, on la calcule explicitement sur les groupes $\mu$-ordinaires, 
on rappelle les résultats de Faltings et on explique la stratégie de l'article.
La section \ref{sect5} est le coeur technique de l'article, qui détermine la structure de l'image de l'application de Hodge-Tate. Elle est basée sur l'article de Faltings \cite{Fal}, 
on y montre que l'image de l'application de Hodge-Tate contient suffisamment de périodes suivant la signature de $G$, de telle sorte que l'on peut construire une application de Hodge-Tate divisée (en passant à la puissance extérieure), et on relie l'image de cette application divisée (et de l'application de départ) avec les invariants de Hasse partiels de \cite{Her1}.
Dans la section \ref{sect6} on construit la filtration dans la $p$-torsion, et on calcule ses degrés partiels à la manière de \cite{Far}.
La section \ref{sect7} décrit de nouvelles filtrations de Harder-Narasihman, basées sur l'article \cite{FarHN}, à l'aide de fonctions degrés "modifiées". 
On donne aussi des conditions pour que des sous-groupes de grand degrés soient des crans de ces filtrations.
La section \ref{sect8} est consacré au théorème général et utilise les constructions des sections précédentes. Un grande partie de la démonstration est dédiée à des minorations de degrés, pour appliquer les résultats de la section \ref{sect7}. On y prouve aussi que l'on peut construire des déformation de Frobénius à l'aide de combinaisons linéaires des filtrations canoniques. Enfin, la section \ref{sect9} traite la mise en famille des filtrations précédentes. 

\subsection{Remerciements}

Je tiens à exprimer mes remerciements envers Laurent Fargues et Vincent Pilloni pour m'avoir introduit à ce sujet ainsi que pour leur aides et leurs encouragements tout au long de 
l'écriture de cet article. Je remercie aussi sincèrement Stéphane Bijakowski pour toutes les discussions intéressantes, ainsi que sa capacité à expliciter les choses les plus 
abstraites. En particulier il me semble qu'il aurait été difficile de conclure cet article sans les propositions de \cite{Bij}.

\section{Notations et mise en situation}
\label{sect2}
 
 Soit $\mathcal O$ l'anneau des entiers d'une extension finie $F$ non ramifiée de $\QQ_p$. Notons $f = [F:\QQ_p]$. 
 On notera $\mathcal I = \Hom(F,\CC_p)$ l'ensemble des plongements de $F$ dans $\CC_p$, et on utilisera communément la notation $\tau$ pour ces plongements. On notera $\sigma$ le relèvement de Frobenius à $F$, qui induit une action transitive sur $\mathcal I$.
 
 Soit $\mathfrak X/\Spf(\mathcal O_C)$ un schéma formel admissible (au sens de Raynaud). Soit $G \fleche \mathfrak X$ un groupe de Barsotti-Tate tronqué d'échelon $r$, 
 muni d'une action de $\mathcal O$, i.e. d'une injection,
 \[\mathcal O \fleche \End_{\mathfrak X}(G).\] On appellera un tel groupe $p$-divisible (tronqué ou non) un $\mathcal O$-module $p$-divisible (tronqué ou non). 
 Notons $\omega_G = e_G^*\Omega^1_{G/\mathfrak X}$ le faisceau (dit conormal) localement libre sur $\mathfrak X$ associé à $G$. 
 Cette action induit des décompositions,
 \[ \omega_G = \bigoplus_{\tau \in \mathcal I} \omega_{G,\tau} \quad \text{et} \quad  \omega_{G^D} = \bigoplus_{\tau \in \mathcal I} \omega_{G^D,\tau},\]
 et supposons que la signature de $G$ soit donnée par $(p_\tau,q_\tau)_{\tau \in \mathcal I}$, c'est à dire,
 \[ p_\tau = \dim_{\mathcal O_C} \omega_{G,\tau} \quad \text{et} \quad q_\tau = \dim_{\mathcal O_C} \omega_{G^D,\tau}.\]
 Notons $H$ la hauteur de $G$, elle est divisible par $f$, et notons $h =\frac{H}{f} =: \Ht_{\mathcal O}(G)$. On a alors pour tout $\tau, p_\tau + q_\tau = h$.
 
 \rem
 Soit $G/\mathfrak X$ un $\mathcal O$-module $p$-divisible (tronqué ou non) de signature $(p_\tau,q_\tau)_\tau$ et de $\mathcal O$-hauteur $h$ comme ci-dessus. Alors son dual de 
 Cartier $G^D$, lui aussi muni d'une action de $\mathcal O$, est aussi de $\mathcal O$-hauteur $h$ et sa signature est donc $(q_\tau,p_\tau)_\tau$.
 \erem
 
 Notons $\overline{X}$ la réduction de $\mathfrak X$ modulo $p\mathcal O_{\mathfrak X}$, et notons toujours $G \fleche \overline{X}$ le changement de base.
 Soit $x = \Spec(k) \dans \overline{X}$ un point géométrique, et $G_x$ le groupe de Barsotti-Tate tronqué associé au-dessus de $x$.
 Le cristal de Dieudonné de $G_x$, noté $\mathbb D(G_x) =: \mathbb D$ est un $W(k)/p^rW(k)$-module libre de rang $H$, muni d'une action de $\mathcal O$ qui induit,
 \[ \mathbb D = \bigoplus_{\tau \in \mathcal I} \mathbb D_\tau,\]
 telle que, \[F : \mathbb D_\tau \fleche \mathbb D_{\sigma\tau} \quad \text{et} \quad V : \mathbb D_{\sigma\tau} \fleche \mathbb D_\tau.\]
 
 On peut associer à $G_x$ deux polygones, le polygone de Newton,
 \[ \Newt_\mathcal O(G_x) = \frac{1}{f} \Newt(\mathbb D_\tau, V^f),\]
 où $\Newt(\mathbb D_\tau, V^f)$ est le polygone de Newton (convexe) associé au $V^f$-cristal $(\mathbb D_\tau, V^f)$ par le théorème de Dieudonné Manin.

 On peut aussi associer à $G_x$ son polygone de Hodge,
 \[ \Hdg_\mathcal O(G_x) = \frac{1}{f} \sum_{\tau \in \mathcal I} \Hdg_\tau(G_x),\]
 où $\Hdg_\tau(G_x)$ est le polygone de pente 0 sur $[0,q_\tau]$ et de pente 1 sur $[q_\tau,f]$.
 
On a alors la proposition classique, voir \cite{RR}, voir aussi \cite{Her1} section 3.

\pro
Le polygone de Newton $\Newt_\mathcal O(G_x)$ est au-dessus du polygone de Hodge $\Hdg_\mathcal O(G_x)$, et ils ont mêmes points terminaux.
\epro
 
 \begin{figure}[h]
 \caption{Exemple de polygones de Hodge et Newton dans le cas de $\mathcal O = \ZZ_{p^2}$}
 \label{HdgNewt}
 \begin{center}
 \begin{tikzpicture}[line cap=round,line join=round,>=triangle 45,x=1cm,y=1cm]
\draw[->,color=black] (-0.5,0.) -- (5.,0.);
\foreach \x in {,1.,2.,3.,4.}
\draw[shift={(\x,0)},color=black] (0pt,2pt) -- (0pt,-2pt);
\draw[->,color=black] (0.,-0.5) -- (0.,4.);
\foreach \y in {-0.5,0.5,1.,1.5,2.,2.5,3.,3.5}
\draw[shift={(0,\y)},color=black] (2pt,0pt) -- (-2pt,0pt);
\clip(-0.5,-0.5) rectangle (5.,4.);
\draw (0.,0.)-- (1.,0.);
\draw (1.,0.)-- (3.,1.);
\draw (3.,1.)-- (4.00163934426,2.20659971306);
\draw [color=ffqqqq] (0.,0.)-- (2.44169908507,0.720849542536);
\draw [color=ffqqqq] (2.44169908507,0.720849542536)-- (4.00163934426,2.20659971306);
\draw (0.209836065574,2.55093256815)-- (0.583606557377,2.55093256815);
\draw [color=ffqqqq] (0.209836065574,2.02725968436)-- (0.588524590164,2.02008608321);
\draw [dash pattern=on 1pt off 1pt] (3.,1.)-- (3.,0.);
\draw (2.99344262295,-0.0961262553802) node[anchor=north west] {$q_{\tau_2}$};
\draw (1.03606557377,-0.117647058824) node[anchor=north west] {$q_{\tau_1}$};
\draw (0.64262295082,2.15638450502) node[anchor=north west] {$\Newt_{\mathcal O}$};
\draw (0.637704918033,2.7230989957) node[anchor=north west] {$\Hdg_{\mathcal O}$};
\draw [dash pattern=on 1pt off 1pt] (4.00163934426,2.20659971306)-- (4.,0.);
\draw (3.96229508197,-0.0602582496413) node[anchor=north west] {$h$};
\end{tikzpicture}
\end{center}
\end{figure}
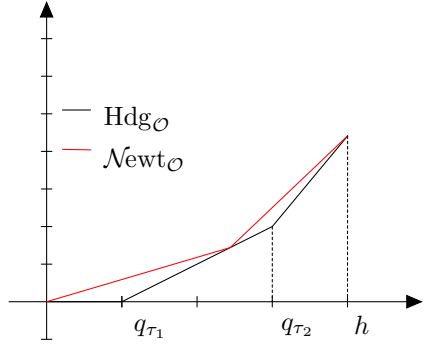

 \rem
 On peut aussi définir $\Newt_\mathcal O^\diamond$ et $\Hdg_\mathcal O^\diamond$ leurs analogues concaves évidents.
 \erem
 
 Ces polygones induisent des fonctions, si on note $\mathcal P$ l'espace des polygones convexes, à abscisses de rupture dans $\frac{1}{f}\ZZ$, 
 définis sur $[0,\Ht_{\mathcal O}(G[p])]$, muni de la topologie de la convergence uniforme,
 \[ \Newt_\mathcal O, \Hdg_\mathcal O : |\overline X| \fleche \mathcal P.\]
 Notons $\mathfrak X^{rig}$ la fibre rigide au sens de Raynaud de $\mathfrak X$, on a une application de spécialisation,
 \[ sp : \mathfrak X^{rig} \fleche \overline{X},\]
 qui nous permet de définir les applications $\Newt_\mathcal O$ et $\Hdg_\mathcal O$ sur $|\mathfrak X^{rig}|$. Ces applications sont semi-continues vers l'espace $\mathcal P$.
 
\defi
Le lieu $\mu$-ordinaire de $\overline X$, noté $\overline{X}^{\mu-ord}$, est l'ensemble des $x \dans \overline{X}$ tels que $\Newt_\mathcal O(x) = \Hdg_\mathcal O(x)$.
On définit $\mathfrak X^{rig,\mu-ord}$ par $sp^{-1}(\overline{X}^{\mu-ord})$.
\edefi

\rem
Ici et plus généralement dans tout le texte, on utilisera $\mu$ comme une simple notation, c'est à dire qu'elle ne fait référence (même si cela ne se voit pas sur la notation) 
qu'à $\mathcal O$ et pas à la signature. Bien sûr être $\mu$-ordinaire au sens précédent pour un groupe $p$-divisible dépend de l'action considérée, et de sa signature pour 
cette action, mais celles-ci seront en général fixés dans tout le texte.
Dans le cas où $\mathfrak X$ provient d'une donnée de Shimura PEL, alors il est associé à cette donnée un certain cocaractère $\mu$, qui dépend de la donnée, 
mais aussi de la signature (fixée par la condition de Kottwitz), et historiquement être $\mu$-ordinaire pour un point de cette variété de Shimura est une définition relative à 
ce cocaractère.
Néanmoins ici être $\mu$-ordinaire ne dépend que de l'action de $\mathcal O$ (qui est la même partout) et de la signature du groupe, mais on fera un abus de notation et on dira 
$\mu$-ordinaire même lorsque la signature varie.
Une meilleure notation serait de dire $\mathcal O$-ordinaire, mais il semble plus simple de garder des notations que tout le monde connaît.
\erem

\section{Invariants de Hasse $\mu$-ordinaires, rappels}
\label{sect3}

Soit $S$ un schéma sur lequel $p\mathcal O_S = 0$. Soit $G$ un groupe $p$-divisible sur $S$ tronqué d'échelon $r$, tel que 
$r > k_\tau = \sum_{\tau'} \max(0,q_\tau - q_{\tau'})$ pour un $\tau \in \mathcal I$ fixé (respectivement $r > \max_\tau k_\tau$). 
Alors on a construit dans \cite{Her1} (cf. aussi \cite{GN} dans le cas des variétés 
de Shimura) un invariant de Hasse partiel $\widetilde{\Ha_\tau}$ (respectivement, un invariant de Hasse $\mu$-ordinaire $\widetilde{^\mu\Ha}$), tels que $\widetilde{^\mu\Ha}$ est le produit
des invariants de Hasse partiels $\widetilde{\Ha_\tau}$.

\defi
Supposons donné une collection $(p_\tau,q_\tau)_{\tau \in \mathcal I}$. Pour tout $s \dans \Gal(\overline{\QQ_p}/\QQ_p)$, $s$ agit sur $(p_\tau,q_\tau)_\tau$ par,
\[s\cdot (p_\tau,q_\tau)_\tau = (p_{s\tau},q_{s\tau})_\tau.\]
On note $E$ la plus petite extension (Galoisienne) de $\QQ_p$ qui fixe $(p_\tau,q_\tau)_\tau$, que l'on appellera corps reflex (associé à la signature $(p_\tau,q_\tau)_\tau$), et
parfois corps reflex local, si l'on veut marquer la différence avec les corps reflex des variétés de Shimura.
\edefi

On note $\kappa_E$ le corps résiduel de $E$. Le champ qui classifie les groupes de Barsotti-Tate tronqué d'échelon $r$ munis d'une action de $\mathcal O$ est un champ lisse sur 
$\mathcal O_E$, cf \cite{Wed2} (il descend d'ailleurs à $\ZZ_p$). La condition qui fixe la signature en détermine une composante connexe, que l'on note 
$\mathcal{BT}_r^{\mathcal O,(p_\tau,q_\tau)}/\Spec(\mathcal O_E)$. On a alors d'après \cite{Her1} la construction suivante,

\thr
Soit $\mathcal{BT}_r^{\mathcal O,(p_\tau,q_\tau)}/\Spec(\kappa_E)$ le champ modulo $p$ des groupes $p$-divisibles tronqués d'échelon $r$, muni d'une action de 
$\mathcal O$, de signature $(p_\tau,q_\tau)$, et $G$ le $\mathcal O$-module universel.
Alors pour $\tau \in \mathcal I$ tel que $r > k_\tau$, il existe sur l'extension des scalaires, \[\mathcal{BT}_r^{\mathcal O,(p_\tau,q_\tau)}\times\Spec(\overline{\FP}),\] 
une section du fibré inversible universel $\det(\omega_{G^D,\tau})^{\otimes(p^f-1)}$,
\[ \widetilde{\Ha_\tau} \dans H^0(\mathcal{BT}_r^{\mathcal O,(p_\tau,q_\tau)},\det(\omega_{G^D,\tau})^{\otimes(p^f-1)}),\]
telle que si $x \dans \mathcal{BT}_r^{\mathcal O,(p_\tau,q_\tau)}\times \Spec(\overline{\FP})(\overline k)$, correspond à $G/\Spec(\overline k)$, 
alors $x^*\widetilde{\Ha}_\tau$ est inversible si et seulement si $\Newt_\mathcal O(G)$ et $\Hdg_\mathcal O(G)$ ont un point de contact en l'abscisse $q_\tau$.

Les sections, pour $q = q_\tau \dans \NN$ tels que $r > k_\tau$,
\[ \widetilde{\Ha}_q := \bigotimes_{\tau' : q_{\tau'} = q} \widetilde{\Ha_{\tau'}} \dans H^0(\mathcal{BT}_r^{\mathcal O,(p_\tau,q_\tau)}\times \Spec(\overline{\FP}),\bigotimes_{\tau' : q_{\tau'} = q}\det(\omega_{G^D,\tau'})^{\otimes(p^f-1)}),\]
redescendent en une section sur $\mathcal{BT}_r^{\mathcal O,(p_\tau,q_\tau)}/\Spec(\kappa_E)$ du faisceau,
\[ \det(\omega_{G^D,q})^{\otimes(p^f-1)},\]
où $\omega_{G^D,q} = \bigoplus_{\tau : q_\tau = q} \omega_{G^D,\tau}$, a priori sur $\mathcal{BT}_r^{\mathcal O,(p_\tau,q_\tau)}/\Spec(\overline{\FP})$, redescend à $\kappa_E$.

Si $r > \max_\tau k_\tau$, il existe aussi sur $\mathcal{BT}_r^{\mathcal O,(p_\tau,q_\tau)}\times \Spec(\kappa_E)$ une section de $\det(\omega_{G^D})^{\otimes(p^f-1)}$,
\[ \widetilde{^\mu\Ha} \dans H^0(\mathcal{BT}_r^{\mathcal O,(p_\tau,q_\tau)},\det(\omega_{G^D,\tau})^{\otimes(p^f-1)}),\]
telle que pour tout $x \dans \mathcal{BT}_r^{\mathcal O,(p_\tau,q_\tau)}(\overline k)$, correspondant à $G/\overline k$,
$G$ est $\mu$-ordinaire si et seulement si $x^*\widetilde{^\mu\Ha}$ est inversible.

De plus, les sections $\widetilde{\Ha_\tau}, \widetilde{\Ha_q}$ et $\widetilde{^\mu\Ha}$ sont des diviseurs de Cartier sur $\mathcal{BT}_r^{\mathcal O,(p_\tau,q_\tau)}\times \Spec(\kappa_E)$ ($\Spec(\kappa_F)$ pour $\widetilde{\Ha_\tau}$), ils sont compatibles aux changements de bases, et à la dualité, c'est à dire,
\[ \widetilde\Ha_\tau(G) = \widetilde\Ha_\tau(G^D), \forall \tau \text{ tels que }r > \max(k_\tau,k_\tau^D = \sum_{\tau'} \max(0,p_\tau-p_{\tau'}),d = \sum_{\tau'} p_{\tau'}).\]
De plus, si on a un scindage de $G$, $G = H \times H'$, qui induit un contact Hodge-Newton en un point de rupture du polygone de Newton, alors,
\[ \widetilde{\Ha_\tau}(G) = \widetilde{\Ha_\tau}(H) \otimes \widetilde{\Ha_\tau}(H'), \forall \tau \in \mathcal I.\]
\ethr

On déduit du théorème précédent, que pour tout groupe de Barsotti-Tate $G/S$ tronqué d'échelon $r$ suffisamment grand, sur une base $S$ de caractéristique $p$,
on peut par tirée en arrière lui associer des invariants de Hasse partiels, et un $\mu$-invariant de Hasse.

\rem
Même si être $\mu$-ordinaire ne dépend que de la $p$-torsion (cf \cite{Moo}), on ne sait construire $\widetilde{^\mu\Ha}$ que pour un $\mathcal O$-module $p$-divisible tronqué 
d'échelon supérieur à $\max_\tau k_\tau+1$. On ne sait pas non plus si l'invariant $\widetilde{^\mu\Ha}$ ne dépend que de la $p$-torsion, mais on espère résoudre ce 
problème dans un futur proche. Notons que c'est le cas pour les invariants construits par \cite{KW,GK}, et \cite{Box} puisqu'ils sont construit sur les champ des $G$-Zip et 
sur les strates d'Ekedahl-Oort respectivement, et que ces derniers ne dépendent "que de la $p$-torsion" des groupes de Barsotti-Tate.
\erem

On a alors le résultat classique suivant,

\thr(Wedhorn, \cite{Wed1})
Soit $\overline{X}$ la fibre spéciale d'une variété de Shimura PEL non ramifiée en $p$, alors,
\[ \overline{X}^{\mu-ord} = \{ x \dans \overline X :  {^\mu}\Ha(x) = 0\},\]
est un ouvert dense.
\ethr

\defi
Dans le cas où on a un $\mathcal O$-module $G/\mathcal O_K$ avec $K$ une extension de $\QQ_p$, comme les faisceaux $\omega_{G^D,\tau}, \omega_{G^D,q}$ sont localement libres,
on identifiera librement les invariants précédents $\widetilde{\Ha_\tau},\widetilde{\Ha_q}$, et $\widetilde{^\mu\Ha}$, avec leur valuation (dans $v(K)$) notée $\Ha_\tau,\Ha_q$, et 
$^\mu\Ha$, qui ne dépendent pas d'un choix d'une base du déterminant des différents faisceaux.
\edefi

\subsection{$\mathcal O$-modules de Lubin-Tate, groupes $p$-divisible $\mu$-ordinaire}
\label{sect31}
\defi
\label{def36}
Soit $A \subset \mathcal I$ un ensemble de plongements. Alors il existe sur $\mathcal O$ un groupe $p$-divisible $\mathcal{LT}_A$ muni d'une action de 
$\mathcal O$ tel que,
\[ \Ht_\mathcal O(\mathcal{LT}_A) = 1, \quad \text{et} \quad p_\tau = 1 \text{ si et seulement si } \tau \dans A.\]
Le display de $\mathcal{LT}_A$ est donné par produit tensoriel de celui de $\mathcal{LT}_\tau$, les $\mathcal O$-modules (stricts) de Lubin-Tate, 
pour $\tau\dans A$. De plus, un tel module sur $\widehat{\ZZ_p^{nr}} = W(\overline{\FP})$ est unique.
\edefi

On a alors le résultat de structure suivant sur $\mathcal O_C$ (\cite{Her1}, un résultat similaire sur $\overline{\FP}$ est obtenu dans \cite{Moo}),

\pro
Soit $G/\mathcal O_C$ un $\mathcal O$-module $p$-divisible (non tronqué) de signature $(p_\tau,q_\tau)$.
Notons,
\[ \{q_\tau, \tau \dans \mathcal I\} = \{q^{(1)}<\dots<q^{(r)}\}, \quad q^{(0)} = 0, q^{(r+1)} = h.\]
Alors $G$ est $\mu$-ordinaire si et seulement si,
\[ G = \prod_{l=0}^r \mathcal{LT}_{A_l}^{q^{(l+1)}-q^{(l)}},\]
où $A_l = \{ \tau \dans \mathcal I : q_\tau \leq q^{(l)}\}$. Les $(A_l)$ forment une suite strictement croissante.
\epro

\rem
Si $A \not\subset A'$ et $A' \not \subset A$, alors $\mathcal{LT}_A \times \mathcal{LT}_{A'}$ n'est pas $\mu$-ordinaire, alors que chacun des facteurs l'est (ces trois groupes ont des
 signatures différentes)...
Cet exemple montre bien l'utilité de la notation $\mu$-ordinaire, c'est-à-dire du fait qu'elle ne fait pas référence à la signature.
\erem

\exe
Si $\mathcal O = \ZZ_{p^2}$ et $\mathcal I = \{ \tau_1,\tau_2\}$, et $q_{\tau_2} \leq q_{\tau_1}$ alors "le" groupe $p$-divisible $\mu$-ordinaire sur $\mathcal O_C$ est donné par,
\[ (\QQ_p/\ZZ_p\otimes_{\ZZ_p}\mathcal O)^{q_{\tau_2}} \times \mathcal{LT}_{\tau_2}^{q_{\tau_1}-q_{\tau_2}} \times (\mu_{p^\infty}\otimes_{\ZZ_p}\mathcal O)^{p_{\tau_1}}.\]
\eexe

Rappelons les résultats de Mantovan-Viehmann \cite{MV} et X.Shen \cite{Shen} (cf. aussi \cite{Moo}),
\thr
\label{thrshen}
Soit $G/\mathcal O_K$ un groupe $p$-divisible muni d'une action de $\mathcal O$, $K$ une extension de $\QQ_p$.
Soit $\overline{G}/\kappa_K$ sa réduction et supposons que $\Newt_\mathcal O(\overline G)$ et $\Hdg_\mathcal O(\overline G)$ ont un point de rupture en l'abscisse $t$ qui
est aussi une abscisse de rupture pour $\Newt_\mathcal O(\overline G)$, alors il existe un unique scindage de $\overline{G}$ sur $\kappa_K$,
\[ \overline{G} = \overline{H_1} \times \overline{H_2},\]
où $\overline{H_1}, \overline{H_2}$ sont des $\mathcal O$-modules $p$-divisibles de hauteurs $t$ et $h-t$ respectivement, tels que leurs polygones de Newton et Hodge soient
ceux de $G$ entre $0$ et $t$ et entre $t$ et $h$ respectivement.
De plus, si $G$ est modulaire (automatique si $K$ est un corps $p$-adique ou $K = C$), alors il existe un unique $0 \subset H_2 \subset G$ $\mathcal O$-module $p$-divisible 
sur $O_K$ qui relève $\overline{H_2}$, et
$G/H_2$ relève $\overline{H_1}$.
Si de plus $G$ est polarisé, il y a une compatibilité à la polarisation.
\ethr

On en déduit en particulier,

\cor
Soit $G/O_K$ un groupe $p$-divisible modulaire $\tau$-ordinaire, i.e. $\Ha_\tau(G) = 0$, alors il existe un sous-$\mathcal O$-module $p$-divisible de $G$ de hauteur $p_\tau$.
Si de plus $G$ est $\mu$-ordinaire, il existe une filtration "canonique" de $G$,
\[ 0 \subset H_1 \subset H_2 \subset \dots \subset H_r = G,\]
où $H_i$ est un $\mathcal O$-module $p$-divisible de hauteur $h-q^{(r-i)}$.
\ecor

\rem
Grâce au travail de Scholze et Weinstein dans \cite{SW}, et en utilisant la théorie des filtrations de Harder-Narasihman de \cite{FarHN} (déjà utilisée dans \cite{Shen}) on peut 
désormais enlever l'hypothèse "modulaire" des énoncés précédents.
\erem

Donc par exemple, sur le lieu $\mu$-ordinaire d'une variété de Shimura $\mathfrak X$ (PEL et non ramifiée en $p$), le groupe $p$-divisible universel possède une filtration canonique, 
qui dans le cas où le lieu ordinaire est non-vide, redonne simplement le sous-groupe canonique, i.e. la partie multiplicative.
On voudrait étendre l'existence de cette filtration canonique à un voisinage strict du lieu $\mu$-ordinaire dans $\mathfrak X^{rig}$. Pour cela on procède comme dans \cite{Far},
en étudiant tout d'abord l'application de Hodge-Tate.

\section{Application de Hodge-Tate}
\label{sect4}
\subsection{Noyau de l'application de Hodge-Tate}

\defi
Soit $G/\mathcal O_K$ un schéma en groupes, $K/\QQ_p$. Son application de Hodge-Tate, entre faisceaux fppf sur $\mathcal O_K$, est 
\[ \alpha_G :
\begin{array}{ccc}
 \underline G &\fleche& \underline \omega_{G^D}   \\
  (\phi : G^D \fleche \mathbb G_m)  & \longmapsto  &  \phi^*\frac{\mathrm dT}{T},
\end{array}
\]
où l'on a identifié les faisceaux fppf $\underline G$ et $\mathcal Hom(\underline G^D,\underline{\mathbb G}_m)$.
\edefi

\rem
Si $G$ est muni d'une action de $\mathcal O$, alors $\omega_{G^D}$ aussi, et $\alpha_G$ est $\mathcal O$-équivariante.
\erem

On utilisera maintenant la notation $\alpha_G$ pour ses $\mathcal O_{\overline K}$-points.

\exe
\label{exeRay}
Soit $G/\mathcal O_K$ un schéma en groupe de Raynaud, \cite{Ray}, associé à la donnée $(\gamma_i,\delta_i)_{i\in \{1,\dots,f\}}$ où, pour tout $i$, $\gamma_i\delta_i = \omega$, une constante universelle de valuation $p$-adique 1, c'est à dire,
\[ G = G^{(\gamma_i,\delta_i)} =\Spec(\quotient{\mathcal O_K[X_i, i = 1 \dots f]}{(X_i^p - \gamma_{i+1}T_{i+1}})),\]
avec la (co-)multiplication est explicite, donnée par le corollaire 1.5.1 de \cite{Ray}.
On a alors,
\[G^D = G^{(\delta_i,\gamma_i)}.\]
On peut donc calculer explicitement,
\[ \omega_{G^D} = \bigoplus_{i=1}^f \quotient{\mathcal O_K}{\delta_i}\mathrm dT_i,\]
et un $\mathcal O_{\overline K}$-point de $G$ est donné par une collection $(x_i)_i \dans (\mathcal O_{\overline K})^{\{1,\dots,f\}}$  telle que, \[x_i^p = \gamma_{i+1}x_{i+1}.\]

L'application de Hodge-Tate est alors donnée par,
\[
\alpha_G :
\begin{array}{ccc}
  G(\mathcal O_{\overline K}) & \fleche & \omega_{G^D} \otimes_{\mathcal O_K} \mathcal O_{\overline K}  \\
  (x_i)_i  & \longmapsto & \sum_i x_i \pmod {\delta_i}\mathrm dT_i
\end{array}
\]
\eexe

Retournons au cas $\mu$-ordinaire. Soit $G/\mathcal O_C$ un $\mathcal O$-module $p$-divisible $\mu$-ordinaire, donc de la forme,
\[ G = \prod_{l=0}^r \mathcal{LT}_{A_l}^{q^{(l+1)}-q^{(l)}},\]
où $A_l = \{ \tau \dans \mathcal I : q_\tau \leq q^{(l)}\}$.
Soit $\tau$ un plongement de $\mathcal O$ dans $\mathcal O_C$, alors comme $G$ est $\mu$-ordinaire, les théorèmes de Shen et Mantovan-Viehmann (\ref{thrshen}) prédisent l'existence d'un sous-groupe $p$-divisible,
qui correspond à,
\[ H_\tau = \prod_{l \geq l_\tau} \mathcal{LT}_{A_l}^{q^{(l+1)}-q^{(l)}},\]
où $q_\tau = q^{(l_\tau)}$.
Or si on regarde l'application de Hodge-Tate modifiée comme suit,
\[ \alpha_{G,\tau} : T_pG \overset{\alpha_G}{\fleche} \omega_{G^D} \twoheadrightarrow \omega_{G^D,\tau},\]
alors son noyau correspond exactement au module de Tate de,
\[ \prod_{l \geq l_\tau} \mathcal{LT}_{A_l}^{q^{(l+1)}-q^{(l)}}.\]
De plus, un tel produit est rationnel, c'est-à-dire que si $G/\mathcal O_K$, alors $\Ker \alpha_{G,\tau}$ est aussi défini sur $\mathcal O_K$.
On en déduit directement la proposition suivante,

\pro
Soit $G \fleche \mathfrak X$ un $\mathcal O$-module $p$-divisible.
Sur le lieu $\mu$-ordinaire, la filtration canonique est donnée par les noyaux $\Ker \alpha_{G,\tau}$, lorsque $\tau$ varie parmi les plongements de
$\mathcal O$ dans $\mathcal O_C$.
\epro

On va essayer de déformer ces applications de Hodge-Tate $\alpha_{G,\tau}$ sur des voisinages stricts du lieu $\mu$-ordinaire afin que leurs noyaux déforment la filtration canonique.

\defi
Soit $M$ un $\mathcal O_K$-module, $K/\QQ_p$. On note, pour tout $\eps \dans v(K)$, \[M_\eps = M \otimes \mathcal O_K/(p^\eps).\]
Soit $G/\mathcal O_K$ un $\mathcal O$-module $p$-divisible tronqué. Supposons $K/\mathcal O[1/p]$. On définit l'application $\alpha_{G,\tau,\eps}$ par,
\[ \alpha_{G,\tau,\eps} : G(\mathcal O_{\overline K}) \overset{\alpha_{G,\tau}}{\fleche} \omega_{G^D,\tau} \twoheadrightarrow \omega_{G^D,\tau,\eps}.\]
\edefi

\pro
\label{proker1}
Soit $G/\mathcal O_C$ un $\mathcal O$-module de Barsotti-Tate tronqué d'échelon $1$, de signature $(p_\tau,q_\tau)_\tau$.
Soit $\eps \dans v(C)$ tel que $\frac{1}{p-1} < \eps \leq 1$, alors,
\[ \dim_{\FF_{p^f}} \Ker \alpha_{G,\tau,\eps} \leq p_\tau.\]
\epro

\dem
La démonstration est identique à \cite{Far}, proposition 9. On a que $\omega_{G^D,\tau} \simeq (\mathcal O_C/p\mathcal O_C)^{q_\tau}.$ Donc si $\delta \dans \mathcal O_C$ est de valuation $\frac{1}{p-1}$,
$\delta\omega_{G^D,\tau,\eps}$ est un $\mathcal O_{C,\eps-\frac{1}{p-1}}$ module libre de rang $q_\tau$. Or d'après le théorème 3 de \cite{Far} sur l'annulation de la cohomologie de la suite de Hodge-Tate, 
\[ \delta \omega_{G^D,\tau,\eps} \subset \mathcal O_{C,1} \otimes_{\FF_{p^f}}\im(\alpha_{G,\tau,\eps}).\]
Or $\mathcal O_{C,1}\im(\alpha_{G,\tau,\eps})$ est de présentation finie (comme sous-module de type fini de $\omega_{G^D,\eps}$ qui est libre), et donc 
$\delta \omega_{G^D,\tau,\eps}$ est engendré par moins d'éléments que $\mathcal O_{C,1}\im(\alpha_{G,\tau,\eps})$, c'est-à-dire que $\dim_{\FF_{p^f}} \im(\alpha_{G,\tau,\eps})$,
et donc,
\[ q_\tau \leq \dim_{\FF_{p^f}} \im(\alpha_{G,\tau,\eps}).\qedhere\]
\edem

On voudrait trouver un $\eps > \frac{1}{p-1}$ qui nous assure que $\Ker(\alpha_{G,\tau,\eps})$ est exactement de hauteur $p_\tau$, malheureusement la méthode utilisée dans
\cite{Far} ne peut pas s'adapter directement, puisqu'on ne veut plus relier ce noyau à $\Ha(G)$, qui est, sur les des variétés (e.g. de Shimura) sans lieu ordinaire, identiquement 1 
(c'est-à-dire que la section $\widetilde{\Ha}(G)$ est identiquement nulle).
Introduisons le relèvement cristallin de l'application de Hodge-Tate, introduit par exemple dans \cite{Far}, section 5.1.2.

\subsection{Application de Hodge-Tate cristalline}

Soit $S$ un schéma, et $G$ un $S$-schéma en groupes, fini localement libre. Supposons $p$ localement nilpotent sur $S$, et que $S$ est un $\Sigma$-schéma, où $\Sigma$ le 
spectre d'un anneau $p$-adiquement complet sans $p$-torsion. $p\mathcal O_\Sigma$ est muni de puissances divisées et on s'intéresse à $\Cris(S/\Sigma)$, le gros site cristallin fppf de $S/\Sigma$. Posons,
\[ \mathbb D(G) = \mathcal Ext^1(\underline{G}^D,\mathcal O_{S/\Sigma}),\]
le cristal de Dieudonné covariant de $G$.
Notons $\mathcal J_{S/\Sigma}$ l'idéal de $\mathcal O_{S/\Sigma}$, et plaçons nous dans la catégorie dérivée des faisceaux en groupes abéliens de $\Cris(S/\Sigma)$, $D(S/\Sigma)$. Alors d'après \cite{Far}, section 5.1.1, on a une application dans $D(S/\Sigma)$,\[\alpha_G^{cris} : G \fleche \mathbb D(G).\]
Celle-ci s'inscrit après évaluation sur l'épaississement tautologique, dans le diagramme,
\begin{center}
\begin{tikzpicture}[description/.style={fill=white,inner sep=2pt}] 
\matrix (m) [matrix of math nodes, row sep=3em, column sep=2.5em, text height=1.5ex, text depth=0.25ex] at (0,0)
{ 
 G & & \mathbb D(G)_S \\
& \omega_{G^D} = \mathcal Ext^1(G^D,\mathcal J_{S/\Sigma})_S & \\
 };

\path[->,font=\scriptsize] 
(m-1-1) edge node[auto] {$(\alpha_{G}^{cris})_S$} (m-1-3)
(m-1-1) edge node[auto] {$\alpha_G$} (m-2-2)
(m-2-2) edge node[auto] {$$} (m-1-3);
\end{tikzpicture}
\end{center}
et relève donc l'application de Hodge-Tate $\alpha_G$.

Soit de plus $(U,T,\delta)$ un ouvert de $\Cris(S/\Sigma)$. Soit $I$ l'idéal noyau de $\kappa : \mathcal O_T \fleche \mathcal O_U$. Notons $\mathcal E = \mathbb D(G)_{(U,T,\delta)}$.
On a un morphisme, $\omega_{G^D\times U} = \omega_{G^D} \otimes_{\mathcal O_S}\mathcal O_U \fleche \mathcal E/I\mathcal E,$
et notons $\Fil \mathcal E$, la préimage de $\im(\omega_{G^D\times U}  \fleche \mathcal E/I\mathcal E)$ par $\mathcal E \fleche \mathcal E/I\mathcal E$, auquel on pense comme un relèvement de la filtration de Hodge.
Alors on a la factorisation,
\[ \alpha_{G,(U,T,\delta)} : G(U) \fleche \Fil E \subset \mathcal E.\]

Si de plus sur $(U,T,\delta)$ il existe $\phi : T \fleche T$ qui relève le Frobénius de $U_0 = U \times_S S_0$, et qui commute aux puissances divisées de $U_0$ dans $T$, alors
on a un isomorphisme,
\[(\mathcal F^{(p)})_{(U,T,\delta)} = \phi^*(\mathcal F_{(U,T,\delta})),\]
pour tout cristal en $\mathcal O_{S/\Sigma}$-modules $\mathcal F$ (on utilise implicitement l'équivalence entre cristaux de $\mathcal O_{S_0/\Sigma}$-modules et $\mathcal O_{S/\Sigma}$-modules).
Dans ce cas, le morphisme de Verschiebung,
\[ V : \mathbb D(G) \fleche \mathbb D(G)^{(p)},\]
permet d'écrire plus précisément,
\[ \im((\alpha_{G}^{cris})_{(U,T,\delta)}) \subset \{ x \dans \Fil \mathcal E : Vx = x \otimes 1\}.\]

Supposons maintenant que $S = \Spec(\mathcal O_K)$, et $G$ est un groupe plat fini sur $\mathcal O_K$. Notons $S_0 = \Spec(\mathcal O_K/p\mathcal O_K)$,
$G_0 = G \times S_0$, $\overline S = \Spec(\mathcal O_{\widehat{\overline K}})$, $\overline{S_0} = \Spec(\mathcal O_{\overline K}/p\mathcal O_{\overline K})$, et enfin $\Sigma = \Spec(\ZZ_p)$.
Le morphisme $\phi : \mathcal O_{\overline{S_0}} \fleche \mathcal O_{\overline{S_0}}$ est surjectif, donc $\Cris(\overline{S_0}/\Sigma)$ possède un objet initial, l'épaississement,
\[ A_{cris} \overset{\theta}{\fleche} \mathcal O_{\widehat{\overline K}} \fleche \mathcal O_{\overline K}/p\mathcal O_{\overline K}.\]
$A_{cris}$ est muni de puissances divisées, relativement à $\ker \theta$, d'un Frobenius cristallin, et d'une action de $\Gal(\overline K/K)$.
On peut alors évaluer $\mathbb D(G)$ sur $(S_0,\Spec A_{cris})$, 
\[ E = H^0(\overline{S_0}/\Sigma, \mathbb D(G_0)) = \varprojlim_n \mathbb D(G_0)_{(A_{cris}/p^nA_{cris} \twoheadrightarrow \mathcal O_{\overline K}/p\mathcal O_{\overline K})}.\]
C'est un $A_{cris}$-module (et même un $A_{cris}/p^rA_{cris}$-module libre si $G$ est un Barsotti-Tate tronqué de rang $r$), muni d'un Verschiebung $A_{cris}$-linéaire,
\[ V : E \fleche E\otimes_{A_{cris},\phi} A_{cris} = E^{(\phi)},\]
et d'un Frobenius $A_{cris}$-linéaire,
\[ F : E \otimes_{A_{cris},\phi} A_{cris} = E^{(\phi)} \fleche E.\]
On peut définir l'application que l'on notera encore $\alpha_G^{cris}$, comme la composée,
\[ G(\mathcal O_{\overline K}) \fleche G(\mathcal O_{\overline K}/p\mathcal O_{\overline K}) \overset{(\alpha_{G}^{cris})_{(S_0,\Spec A_{cris})}}{\fleche} E,\]
qui a son image dans $\{ x \dans \Fil E : Vx = x \otimes 1\}$, et $\Fil E$ est la préimage de la filtration de Hodge par $\theta$.

Si $G$ est un groupe de Barsotti-Tate, les constructions précédentes, pour $G[p^n]$, induisent,
\[ \alpha_G^{cris} : T_pG \fleche \Fil E.\]
On peut définir $\Phi : E \fleche E$ qui est $(A_{cris},\phi)$-linéaire, par $\Phi(x) = F(x\otimes 1)$. L'application $\alpha_G^{cris}$ se factorise alors par
\[ (\Fil E)^{\Phi = p}.\]
On a alors le théorème suivant dû à Faltings,
\thr[Faltings \cite{Fal}, cf. \cite{Far} Théorème 1, \cite{Chen} Théorème 4.1]
Soit $p \neq 2$ et $G/\mathcal O_K$ un groupe $p$-divisible.
\begin{enumerate}
\item L'application de Hodge-Tate induit un isomorphisme,
\[\alpha_G^{cris} : T_pG \overset{\sim}{\fleche} (\Fil E)^{\Phi = p}.\]
\item Si on filtre $E$ par $\Fil^{-1}E = E, \Fil^0 E = \Fil E$ et \[\Fil^i E = \Fil^iA_{cris} \Fil E + \Fil^{i+1}A_{cris}E,\] 
alors on a des inclusions,
\[ tE \subset T_pG \otimes_{\ZZ_p} A_{cris} \subset E,\]
strictement compatibles aux filtrations (où $\Fil^i T_pG \otimes_{\ZZ_p} A_{cris} = T_pG\otimes \Fil^iA_{cris}$).
\end{enumerate}
\ethr

\rem
La preuve (dont on rappelle une esquisse à la section suivante dans un cadre legerement plus général) se fait par étude de la $p$-torsion puis relèvement au $\mathcal{BT}$ complet $G$.
En particulier, si $G$ est un $\mathcal{BT}_1$,
\[ \alpha_G^{cris} : G(\mathcal O_{\overline K}) \overset{\sim}{\fleche} \Fil(E/pE)^{\frac{1}{p}\Phi = \id}.\]
\erem

\subsection{Stratégie}

L'application $\theta$ réduit $\Fil E$ sur $\omega_{G^D} \otimes \mathcal O_{\overline K}$, et donc pour calculer la dimension de l'image de $\alpha_{G,\eps}$, il suffit de calculer celle de l'image de la réduction modulo $\theta$ de $(\Fil E/p^\eps E)^{\frac{1}{p} \Phi = \id}$.

Rappelons que l'on voudrait, pour $G$ un $\mathcal O$-module de Barsotti-Tate sur $\mathcal O_K$, sous une certaine hypothèse sur $\Ha_\tau(G)$, trouver un 
$\eps$ tel que le noyau de $\alpha_{G,\tau,\eps}$ soit de dimension $p_\tau$.

Dans le cas d'un lieu ordinaire non vide, cf. \cite{Far}, on peut remarquer que l'image de $\alpha_{G[p]}^{cris}$ se réduit par $\theta$ dans 
$(\omega_{G[p]^D}\otimes O_{\overline K})^{V = \id \otimes 1}$, 
et une méthode de Newton classique (déjà dans \cite{AG}, et sous la forme d'un lemme de Elkik dans \cite{Far}), permet quitte à réduire modulo $p^{1-\Ha(G)}$, de 
majorer la dimension de l'image de $\alpha_G$.
Dans notre cas, on décompose le cristal $E$ de $G$,
\[ E = \bigoplus_{\tau} E_\tau,\]
et $V$ est $\sigma^{-1}$-linéaire : $V : E_\tau \fleche E_{\sigma{-1}\tau}^{(\phi)}$. On considère alors l'application,
\[ \alpha_{G[p],\tau}^{cris} : G[p](\mathcal O_{\overline K}) \fleche E/pE \twoheadrightarrow E_\tau/pE_\tau,\]
et sa réduction modulo $\theta$ arrive dans $\omega_{G[p]^D,\tau}^{V^f = \id\otimes 1}$, seulement le déterminant de $V^f$ est nul (s'il n'y a pas de lieu ordinaire), on ne peut donc pas appliquer de méthode de Newton.

Pour contourner cela il faut ruser et essayer d'appliquer la méthode de Newton à $\widetilde{\Ha}_\tau$, à la place de $V$ ou $V^f$, qui est construit comme 
$\frac{1}{p^{k_\tau}}V^f$ sur un "relèvement du cristal", et dont le déterminant donne $\Ha_\tau(G)$. Malheureusement pour relier $\alpha_{G,\tau}$ à $\widetilde{\Ha}_\tau$, 
il faudrait que l'image de $\alpha_{G,\tau}^{cris}$ 
arrive non pas dans \[\im\left(\{ x \dans E : \Phi x = px\} \fleche \{ x \dans E_\tau : \Phi^f x = p^fx\}\right),\]
mais dans
\[\{ x \dans E_\tau : V^f x = p^{k_\tau}x \otimes 1\} = "\{ x \dans E_\tau : \widetilde{\Ha}_\tau x = x \otimes 1\}".\]
Seulement voilà, $\widetilde{\Ha}_\tau$ n'existe pas sur $E_\tau$, mais sur $\bigwedge^{q_\tau} E_\tau$, voir section \ref{sectmuha}, et les solutions de 
$V^f x = p^{k_\tau}x\otimes 1$ sont (probablement) beaucoup trop petites si on ne passe pas à la puissance extérieure.

L'idée est alors de remarquer que, pour $\eps \leq 1$,
\[ \left(\dim_{\FF_{p^f}} \Ker \alpha_{G[p],\tau,\eps} \geq p_\tau\right) \Leftrightarrow \left( \dim_{\FF_{p^f}} \im \alpha_{G[p],\tau,\eps} \leq q_\tau\right) \Leftarrow  \left(\rg_{W(\FF_{p^f})/p^{q_\tau\eps} W(\FF_{p^f})} \im \bigwedge^{q_\tau}
\alpha_{G,\tau,q_\tau\eps} \leq 1\right),\]
où le rang est le nombre minimal de générateurs du module $\im \bigwedge^{q_\tau}
\alpha_{G,\tau,q_\tau\eps}$ (qui n'est pas libre).
Et d'essayer de majorer le rang de 
$\im \bigwedge^{q_\tau}
\alpha_{G,\tau,q_\tau\eps}$, en fonction de $\Ha_\tau(G)$.

\section{Cristaux filtrés et suppression de périodes}
\label{sect5}

Dans toute cette section, on considère $G$ un $\mathcal O$-module $p$-divisible (non tronqué) sur $\mathcal O_C$, de signature $(p_\tau,q_\tau)_\tau$ et de $\mathcal O$-hauteur $h = p_\tau + q_\tau$.
Son cristal $E := H^0(\overline{S_0}/\Sigma,\mathbb D(G_0))$ est un $A_{cris}\otimes_{\ZZ_p} \mathcal O = \prod_\tau A_{cris}$-module, et se scinde donc,
\[ E = \bigoplus_{\tau} E_\tau.\]

\pro
Soit $M$ un $R_1\times R_2$-module.
Alors $M \simeq M_1 \oplus M_2$, avec $M_1$ un $R_1$-module, $M_2$ un $R_2$-module.
On a alors l'égalité,
\[\bigwedge_{R_1\times R_2}^n M =  \bigwedge^{n}_{R_1} M_1 \oplus \bigwedge^{n}_{R_2} M_2.\]
\epro

\dem
Soit $e_1 = (1,0) \dans R_1\times R_2$ et $e_2 = (0,1)$, on pose $M_1 = e_1M$ et $M_2 = e_2M$, comme $e_1 + e_2 = 1$ et $e_1e_2 = 0, M = M_1 \oplus M_2$.
Il suffit alors de remarquer que,
\begin{eqnarray*} M\otimes_{R_1 \times R_2} M = (M_1\oplus M_2) \otimes_{R_1\times R_2}(M_1\oplus M_2) \\= M_1\otimes_{R_1\times R_2} M_1 \oplus M_1 \otimes_{R_1\times R_2} M_2 \oplus M_2 \otimes_{R_1\times R_2} M_1 \oplus M_2 \otimes_{R_1\times R_2} M_2,\end{eqnarray*}  
et que $M_1 \otimes_{R_1\times R_2} M_2 = e_1M \otimes_{R_1\times R_2} e_2M = e_1e_2 M\otimes M = 0$.
De même pour $M_2 \otimes M_1$ et $M_1 \otimes_{R_1\times R_2} M_1 = M_1 \otimes_{R_1} M_1$. Le résultat s'en déduit.
\edem

 \subsection{Rappels sur $A_{cris}$.}
 \label{sect42}
 
 La construction suivante est due à Fontaine \cite{Fon}, voir aussi \cite{Chen}, section 2. Soit $C = \widehat{\overline{\QQ_p}}$, et,
 \[ \mathcal O_C^\flat = \varprojlim_{\Frob} (\mathcal O_C/p) = \varprojlim_{x \mapsto x^p} (\mathcal O_C).\]
 Considérons l'anneau des vecteurs de Witt, $W(\mathcal O_C^\flat)$, muni d'un Frobenius et d'une flèche "première coordonnée",
 \[ \theta : W(\mathcal O_C^\flat) \fleche \mathcal O_C^\flat \fleche \mathcal O_C.\]
 Complétons alors l'enveloppe à puissances divisées de $(W(\mathcal O_C^\flat),\ker\theta)$ pour la topologie $\ker \theta$-adique, on obtient donc l'anneau $A_{cris}$, 
 muni d'un Frobénius $\phi$, 
 d'une filtration 
 $(\Fil^iA_{cris})_{i \in \ZZ}$ et d'une application,
 \[ \theta : A_{cris} \fleche \mathcal O_C.\]
 La filtration vérifie que $\Fil^iA_{cris} = A_{cris}$ pour tout $i \leq 0$, $\Fil^1A_{cris} = \ker \theta$ et \[\Fil^iA_{cris}\Fil^jA_{cris} \subset \Fil^{i+j}A_{cris}, \]
 mais l'inclusion est stricte en général. On a une inclusion $\ZZ_p^{nr} = W(\overline{\FP}) \subset A_{cris}$. Il existe un élément
$t$ dans $A_{cris}$ (bien définie à un élément de $\ZZ_p^\times$ près) tel que $t \dans \Fil^1A_{cris}$ et vérifie $\phi(t) = pt$.
 En fait, on a $(\Fil^1A_{cris})^{\phi = p} = t\ZZ_p$, ce qui explique le choix dans la définition de $t$ (c'est le module de Tate du groupe multiplicatif).
 Pour chaque choix de $\tau \dans \Hom(F,C) = \mathcal I$, il existe un élément $t_\tau \dans A_{cris}$, bien défini à un élément de $\ZZ_{p^f}^\times$ (i.e. une racine $p^f-1$-ième)
 près. Cet élément $t_\tau$ est une "période" du groupe de Lubin-Tate $\mathcal{LT}_\tau$, et on en donnera une explication dans la section suivante.
 Étant fixé l'application $\theta$ précédente, et donc la filtration de $A_{cris}$, seul un des $(t_\tau)_\tau$ appartient à $\Fil^1A_{cris}$, les autres sont dans $A_{cris}\priv\Fil^1A_{cris}$.

 Notons $\tau$ le plongement correspondant et $t_\mathcal O$ ce $t_\tau$. Cela vient du fait que si l'on compose à droite $\theta$ par $\sigma$ un autormophisme de $\mathcal O_C$, on change l'application $\theta$, et donc son noyau et donc change $t_\mathcal O$. On peut construire directement $t_\mathcal O$ à partir d'un module de Lubin-Tate, voir section \ref{sectLT}.
 
 Les autres $t_{\sigma^i\tau}$ s'en déduisent par, pour $i \dans \{1,\dots,f-1\}$,
 \[ t_{\sigma^i\tau} = \frac{1}{p}\phi^i(t_\mathcal O).\]
 À partir de maintenant on oublie la notation $t_\tau$, et on garde seulement $t_\mathcal O$ et ses composées par Frobénius. Ceci à le mérite de ne plus dépendre d'aucun choix, si ce n'est uniquement du choix de $\theta$, qui est fixé partout désormais.
 Notons que, parallèlement au fait que le caractère cyclotomique soit un produit de caractères de Lubin-Tate, on a l'égalité suivante, à un inversible dans 
 $\ZZ_p^\times$ près,
 \[ t = t_\mathcal O\prod_{i=1}^{f-1} \frac{1}{p}\phi^i(t_\mathcal O).\]

 D'après \cite{FGL} proposition C.2.8 p 126, pour tout $i \dans \{1,\dots,f-1\}$, l'image de $\frac{1}{p}\phi^i(t_\mathcal O)$ dans $Gr^0A_{cris} = A_{cris}/\Fil^1A_{cris} \overset{\theta}{\simeq}\mathcal O_C$ est de valuation,
 \[ \frac{p^i}{p^f-1},\]
 et l'image de $t_\mathcal O$ dans $\Gr^1A_{cris}$, après un choix d'un isomorphisme de $\Gr^1A_{cris}$ avec $\mathcal O_C$, est de valuation (qui ne dépend pas de ce choix),
 \[ \frac{1}{p^f-1}.\]
 De ces valuation et l'équation précédente on retrouve le fait classique que l'image de $t$ dans $\Gr^1A_{cris}$ est de valuation $\frac{1}{p-1}$.
 
\rem
On peut aussi retrouver ces valuation, en considérant la $p$-torsion des $\mathcal{LT}_\tau$, pour $\tau \in \mathcal I$, celle-ci est alors un schéma en groupes de type $(p,\dots,p)$
au sens de Raynaud, et donc on peut retrouver ces valuation (qui sont les valuation de l'image de certains éléments par $\alpha_{\mathcal{LT}_\tau}$) 
grâce à l'exemple \ref{exeRay} et au théorème de Faltings.
\erem

\subsection{Le cristal $\Lambda$}

Dans cette section on fixe un $\tau$ tel que $q_\tau \not\in\{0,h\}$ et on va s'intéresser à,
\[\Lambda = \bigwedge_{A_{cris}\otimes_{\ZZ_p}\mathcal O}^{q_\tau} E = \bigoplus_{\tau'} \bigwedge^{q_\tau} E_{\tau'}.\]
Rappelons que l'on a noté $f = [\mathcal O:\ZZ_p]$.
Pour plus de simplicité, notons désormais $A := A_{cris}$, $\phi : A \fleche A$ le frobénius de $A$, $\Phi : E \fleche E$ le frobénius de $E$ qui est $(A,\phi)$-linéaire, on garde la notation $\Phi$ pour la puissance extérieure $q_\tau$-ième de $\Phi$ sur $\Lambda$, et $\Lambda_{\tau'} = \bigwedge^{q_\tau} E_{\tau'}.$
Rappelons que $\Phi : E_{\tau'} \fleche E_{\sigma\tau'}$ et donc $\Phi : \Lambda_{\tau'} \fleche \Lambda_{\sigma\tau'}$, où $\sigma$ est le Frobénius de $\mathcal O$.

On a la $\ZZ$-filtration de $E$ donnée par,
\[ \Fil^i E = 
\left\{
\begin{array}{cc}
E  & \text{si } i \leq -1 \\
\Fil E  & \text{si } i = 0    \\
\Fil^iA \cdot \Fil E + \Fil^{i+1}A \cdot E  &\text{si } i \geq 1   
\end{array}
\right.
\]
Elle induit une $\ZZ^f$-filtration sur $\Lambda$ par,
\[ \Fil^{\underline a} \Lambda = \bigoplus_{\tau'} \Fil^{a_{\tau'}}\Lambda_{\tau'}, \quad \underline a \dans \ZZ^f,\]
où
\[\Fil^{a_{\tau'}}\Lambda_{\tau'} = \im\left( \sum_{
\begin{array}{c}
i_1,\dots,i_{q_\tau} \\
i_1 + \dots + i_{q_\tau} = a_{\tau'}
\end{array}
} \Fil^{i_1}E_{\tau'} \otimes \dots \otimes \Fil^{i_{q_\tau}}E_{\tau'} \fleche \bigwedge^{q_\tau} E_{\tau'}\right).\]
En particulier, pour tout $\underline a \leq -\underline{q_\tau} = (-q_\tau,\dots,-q_\tau)$,
\[ \Fil^{\underline{a}}\Lambda = \Lambda.\]

L'image de $\alpha_G^{cris}$ est incluse dans les $x \dans E$ qui vérifient l'équation $\Phi = p$ et donc l'image de l'application,
\[ \bigwedge^{q_\tau} \alpha_{G}^{cris} : \bigwedge^{q_\tau} T_pG \fleche \Lambda,\]
est incluse dans $\Fil^{\underline 0} \Lambda$ et même dans,
\[ \{x \dans \Fil^{\underline 0} \Lambda : \Phi x= p^{q_\tau}x\}.\]
Et donc l'image de l'application $ \bigwedge^{q_\tau} \alpha_{G,\tau}^{cris}$ est incluse dans,
 \[\{ x \dans \Fil^0 \Lambda_\tau : (\Phi^f)_{|\Lambda_\tau}x = p^{fq_\tau}x\}.\]

 \subsection{Changement d'équation.}
 
Il faut penser à l'équation \[(\Phi^f)_{|\Lambda_\tau} = p^{fq_\tau},\] comme se réduisant dans $\bigwedge^{q_\tau}\omega_{G^D,\tau} =\det\omega_{G^D,\tau}$ sur l'équation $V^f = \id\otimes1$.
Or pour utiliser une méthode de Newton sur $\det\omega_{G^D,\tau}$ que l'on pourrait relier à $\Ha_\tau(G)$, il faudrait une équation qui se réduise sur $\widetilde{\Ha}_\tau = \id \otimes 1$ c'est à dire $V^f = p^{k_\tau}\id \otimes 1$ sur $\Lambda_\tau$, par exemple $\Phi^f = p^{fq_\tau - k_\tau}$.

Comment passer de $(\Phi^f)_{|\Lambda_\tau} = p^{fq_\tau}$ à $(\Phi^f)_{|\Lambda_\tau} = p^{fq_\tau - k_\tau}$ ?
Dans $A$, on a la période $t$ qui vérifie $\phi(t) = pt$ donc, si $s \leq r$ et $x \dans A$ vérifient,
\[ \phi(t^sx) = p^r(t^sx),\]
alors comme $A$ est sans $t$-torsion, on en déduit que $x \dans A$ vérifie,
\[ \phi(x) = p^{r-s}x.\]
Autrement dit, quitte à diviser par la période $t$ (et pour nous plus exactement on divisera par $t_\mathcal O$ et les autres périodes de Lubin-Tate), on peut modifier la puissance de $p$ des équations qui apparaissent naturellement sur le cristal.

\subsection{Cas des $\mathcal O$-modules de Lubin-Tate}
\label{sectLT}
Donnons l'exemple de la stratégie de modification de l'équation dans le cas des $\mathcal O$-modules de Lubin-Tate.

Soit $\tau : F \hookrightarrow C$, et $G = \mathcal{LT}_\tau$ le $\mathcal O$-module de Lubin-Tate (sur $\mathcal O_C$ disons) associé à $\tau$. Il est de hauteur $f$ et dimension 1.
On a $q_\tau = 0$ et pour tout $\tau' \neq \tau, q_{\tau'} =1$. Indexons les plongements par $\{0,\dots,f\}$ de telle manière que $\tau$ corresponde à 0 et $\sigma^i\tau$ corresponde à 
$i$.
On peut calculer explicitement son cristal :
\[ E = Ae_0 \oplus \bigoplus_{i=1}^{f-1} Ae_i,\]
\[ \Fil E = \Fil^1Ae_0 \oplus \bigoplus_{i=1}^{f-1}Ae_i,\]
et le Frobénius $\Phi$ est donné par
\[ \Phi e_{i-1} = 
\left\{
\begin{array}{cc}
 pe_i & \text{si } i \neq 1   \\
e_1  & \text{sinon.}  
\end{array}
\right.
\]
D'après le théorème de Faltings, on sait que $(\Fil E)^{\Phi = p}$ est de hauteur $f$ sur $\ZZ_p$, c'est même un $\mathcal O$-module de dimension 1.
De plus, par définition de $t_\mathcal O$, on a un isomorphisme par projection sur $e_0$,
\[ (\Fil E)^{\Phi = p}\overset{\sim}{\fleche} (\Fil^1A)^{\phi^f = p} = \ZZ_{p^f}t_{\mathcal O}.\]
Et on peut explicitement calculer,
\[ (\Fil E)^{\Phi = p} = \{ x_0t_\mathcal Oe_0 + \sum_{i=1}^{f-1} \frac{1}{p}\phi^i(x_0)\phi^i(t_\mathcal O)e_i : x_0 \dans \ZZ_{p^f}\}.\]
On peut alors explicitement construire une flèche,
\[ m : 
\begin{array}{ccc}
 \sum_i Ae_i &\fleche & \sum_i Ae_i  \\
x_ie_i , i \neq 0 & \longmapsto  & \frac{1}{p}x_i\phi^i(t_\mathcal O)e_i \\
x_0e_0 & \longmapsto &x_0t_\mathcal O e_0
\end{array}
\]
et on peut facilement vérifier qu'elle induit un isomorphisme entre,
\[ E^{\Phi = D} \overset{\overset{m}{\sim}}{\fleche} (\Fil E)^{\Phi = p},\]
où $D$ dans la base $(e_0,\dots,e_{f-1})$ est la matrice 
\[
\left(
\begin{array}{ccccc}
 p &   &   & &\\
  & 1  & &   & \\
    &  &p &   & \\
  &   & & \ddots &  \\
  & & & & p
\end{array}
\right)
\]
Sur $i \neq 0$, l'équation $\Phi = D$ devient sur $(\Phi^f)_{|E_i} = p^{f-1}$ et donc se projette sur $\omega_{G^D,i}$ sur $\frac{1}{p}V^f = \id \otimes 1$, dont les solutions par une méthode de Newton se relient à $\Ha_i(G)$ (qui est nul ici d'ailleurs).
 
\subsection{Cas des Lubin-Tate généralisés.}

Comme on a noté $A$ pour $A_{cris}$, on note $S \subset \mathcal I=Hom(\mathcal O,\mathcal O_C)$. On note $G = \mathcal{LT_S}$ le groupe $p$-divisible sur $\mathcal O_C$ avec action de $\mathcal O$ donné par la définition \ref{def36}. Sa $p$-torsion $\mathcal{LT}_S[p]$ est un schéma en groupes de type $(p,\dots,p)$ au sens de Raynaud (\cite{Ray}).
Son cristal sur $A$ est alors donné par,
\[E = \bigoplus_{\tau} Ae_\tau, \quad \Fil E = \bigoplus_{\tau} \Fil E_\tau := \bigoplus_\tau (\Fil^{\delta^{\tau \in S}}A)e_\tau,\]
c'est à dire $\Fil E_\tau = \Fil^1 A e_\tau$ si $\tau \in S$ et $\Fil E_\tau = Ae_\tau = E_\tau$ sinon.
Le Frobénius $(A,\phi)$ linéaire est donné par,
\[\Phi(e_{\sigma^{-1}\tau}) = 
\left\{
\begin{array}{cc}
e_\tau   & \text{si } \sigma^{-1}\tau \in S   \\
pe_\tau   &     \text{si } \sigma^{-1}\tau \not\in S 
\end{array}
\right.
\]
Lorsque $S = \{\tau\}$, on a bien le cristal d'un Lubin-Tate donné précédement.

On s'interresse alors au module de Tate, $T_pG = (\Fil E)^{\Phi = p}$, qui est identifié à,
\[\{ (x_\tau)_\tau \dans \prod_\tau \Fil^{\delta^{\tau \in S}}A : \phi(x_{\sigma^{-1}\tau}) = x_\tau \text{ si } \sigma^{-1}\tau \not\in S, \phi(x_{\sigma^{-1}\tau}) = px_\tau 
\text{ sinon}\}.\]
Supposons $S \neq \mathcal I$, sinon $\mathcal{LT}_\mathcal I = \mu_{p^\infty}\otimes_{\ZZ_p} \mathcal O$ et on connait bien son cristal, et soit $\tau_0 \not\in S$.
On peut alors identifier, par projection sur $e_{\tau_0}$; l'ensemble précédent à,
\[\{ x\dans A : \phi^f(x) = p^{|S|}x \text{ et }\forall j; \sigma^j\tau_0 \dans S, \phi^j(x) \dans \Fil^1A\}.\]
Cet ensemble contient,
\[ \ZZ_{p^f} \prod_{j=1}^{f-1} \left(\frac{1}{p}\phi^j(t_\mathcal O)\right)^{\delta^{\sigma^{-j}\tau_0} \in S}.\]
En effet, comme $\phi^f(t_\mathcal O) = pt_\mathcal O$ et $\tau_0 \not\in S$, on vérifie facilement que $\phi^f = p^{|S|}$ sur l'espace précédent. 
De plus, on vérifie que pour tout $x$ dans cet espace, $\phi^k(x) \dans \Fil^1A$ si et seulement si il existe $j \dans \{1,\dots,f-1\}$ tel que $k = f-j$ et 
$\sigma^{-j}\tau_0 \dans S$, i.e. si $\sigma^k\tau_0 \dans S$.

On peut donc y identifier le sous-module de $(\Fil E)^{\Phi = p}$ donné par,
\[ T_0 = \left\{ xe_{\tau_0} + \sum_{j=1}^{f-1}\phi^j(x)(\frac{1}{p})^{|S \cap \{\tau_0,\sigma\tau_0,\dots,\sigma^{j-1}\tau_0\}|}e_{\sigma^j\tau_0} : x \dans \ZZ_{p^f} \prod_{j=1}^{f-1} \left(\frac{1}{p}\phi^j(t_\mathcal O)\right)^{\delta^{\sigma^{-j}\tau_0} \in S}\right\}\]
C'est un sous-$\ZZ_{p^f}$-module libre de rang 1 de $(\Fil E)^{\Phi = p}$, qui est lui aussi libre de rang 1 sur $\ZZ_{p^f}$ d'après \cite{Fal}, or ils ont le même nombre de points modulo $p$, $(\Fil E)^{\frac{1}{p}\Phi = \id}$ en a $p^f$ d'après Faltings, et $T_0$ en a $p^f$, car la valuation de,
\[\prod_{j=1}^{f-1} \left(\frac{1}{p}\phi^j(t_\mathcal O)\right)^{\delta^{\sigma^{-j}\tau_0} \in S}\]
dans $\Gr^0A_{cris}$ est,
\[ \sum_{j=0}^{f-1} \frac{p^j \delta^{\sigma^{-j}\tau_0 \in S}}{p^f-1} < \frac{1}{p^f-1} \sum_{j=0}^{f-1} p^j = \frac{1}{p-1},\]
d'après les rappels de la sous-section \ref{sect42}.

En particulier, on peut construire une application $m$,
\[ m : 
\begin{array}{ccc}
 \sum_j Ae_{\sigma^j\tau_0} &\fleche & \sum_j Ae_{\sigma^j\tau_0}   \\
x_je_{\sigma^j\tau_0} ,  j \neq 0 & \longmapsto  & \frac{1}{p^{|S \cap \{\sigma\tau_0,\dots,\sigma^{j}\tau_0\}|}}x_i\prod_{k=1}^{f-1} \left(\frac{1}{p}\phi^{k+j}(t_\mathcal O)\right)^{\delta^{\sigma^{-k}\tau_0 \in S}}e_{\sigma^j\tau_0} \\
x_0e_{\tau_0} & \longmapsto & x_0\prod_{k=1}^{f-1} \left(\frac{1}{p}\phi^k(t_\mathcal O)\right)^{\delta^{\sigma^{-k}\tau_0} \in S} e_0
\end{array}
\]
et vérifier qu'elle induit un isomorphisme entre,
\[ E^{\Phi = D_S} \overset{\sim}{\fleche} (\Fil E)^{\Phi = p},\]
où $D_S$ a pour matrice dans la base $e_{\tau_0},e_{\sigma\tau_0},\dots,e_{\sigma^{-1}\tau_0}$, 
\[
\left(
\begin{array}{ccccc}
 p^{\delta_{\sigma^{-1}\tau_0 \not\in S}} &   &   & &\\
  & p  & &   & \\
    &  &p^{\delta_{\sigma\tau_0 \not\in S}}  &   & \\
  &   & & \ddots &  \\
  & & & & p^{\delta_{\sigma^{-2}\tau_0 \not\in S}} 
\end{array}
\right), \quad \text{i.e. } D_S(e_{\sigma^j\tau_0}) = p^{\delta_{\sigma^{j-1}\tau_0 \not\in S}}e_{\sigma^{j}\tau_0}.
\]

\rem
On peut relier le résultat précédent à l'application de Hodge-Tate cristalline, qui identifie $T_p\mathcal{LT}_S$ avec $(\Fil E)^{\Phi = p}$ et qui modulo $\ker \theta$ redonne 
l'application de Hodge-Tate. En particulier, grâce à la description explicite des groupes $p$-divisibles $\mu$-ordinaires, on peut voir que l'application de multiplication 
$m$ corrige le défaut de surjectivité de l'application de Hodge-Tate cristalline, au sens où la projection sur $\omega_{G^D,\tau}$ de la réduction modulo $\ker \theta$ 
de $m^{-1} \circ \alpha_G^{cris}$ est surjective.
On n'arrivera pas à montrer directement cet énoncé dans le cas général, c'est pourquoi on doit passer à une certaine puissance extérieure.
\erem

\subsection{Modification de périodes}

On applique les méthodes de la sous-section précédente pour un $\mathcal O$-module de Barsotti-Tate général.
Soit \[\underline{\Lambda} = (\Lambda, \Fil^{\underline{\cdot}}\Lambda, \Phi),\] le $\mathcal O$-cristal filtré introduit précédemment.
Définissons,
\[ m : 
\begin{array}{ccc}
  \Lambda = \bigoplus_{\tau'} \Lambda_{\tau'} & \fleche & \bigoplus_{\tau'} \Lambda_{\tau'}  \\
  x_{\tau'}& \longmapsto  & t_\mathcal O^{\max(0,q_\tau - q_{\tau'})}\left( \prod_{j=1}^{f-1}\left(\frac{1}{p}\phi^j(t_\mathcal O)\right)^{\max(0,q_\tau - q_{\sigma^{-j}\tau'})}\right)x_{\tau'},
\end{array}
\]
l'application de multiplication par les périodes. Elle est $\mathcal O$-équivariante, c'est-à-dire qu'elle se décompose $m = \bigoplus_{\tau'} m_{\tau'}$, où,
\[ m_{\tau'} : \Lambda_{\tau'} \fleche \Lambda_{\tau'}.\]

\rem
L'idée, à la base de la construction de $m$, est que si $\tau'$ vérifie que $q_{\tau'} < q_\tau$, on a divisé $V$ par $p^{{q_\tau} - q_{\tau'}}$ pour construire $\Ha_\tau$, et donc
pour relier notre équation $\Phi = p$ (qui devient $\Phi = p^{q_\tau}$ sur $\Lambda$ et $\Phi^f = p^{fq_\tau}$ sur $\Lambda_\tau$), on veut à chaque fois que 
$q_{\tau'} < q_{\tau}$ diminuer la puissance de $p$ dans l'équation précédente, de $q_\tau - q_{\tau'}$. Mais on veut faire cela de telle manière que l'on garde une équation (cyclique) avec $\Phi$, 
et donc à chaque fois que l'on "divise" par $t_{\mathcal O}^n$ un élément $x$ de $\bigwedge T_pG \cap \Lambda_{\tau'}$, on doit aussi diviser 
$\Phi^j(x) \dans\bigwedge T_pG \cap \Lambda_{\sigma^j\tau'}$ par $\frac{1}{p^n}\phi^j(t_\mathcal O)^n$.

Une autre façon de voir la construction de $m$ est que sur le lieu $\mu$-ordinaire, les $\mathcal O$-modules $p$-divisibles sont explicites et on veut diviser au maximum l'image de $\alpha_G^{cris}$ qui correspond à,
\[\Fil^{\underline 0} \Lambda^{\Phi = p^{q_\tau}},\]
afin de rendre $\alpha_G^{cris}\otimes 1$ surjective. Les périodes qui apparaissent dans $m$ sont alors exactement les périodes du déterminant du cristal du 
sous-$\mathcal O$-module $p$-divisible canonique de hauteur $p_\tau$ (i.e. le cran de hauteur $p_\tau$ de la filtration canonique).
\erem

Introduisons la notation suivante,
\[ \forall \tau', f_{\tau'} = \min(0,q_{\tau'}-q_\tau), \quad \text{et} \quad \underline f = (f_{\tau'})_{\tau'} \dans \ZZ^f.\]
On a que $\underline f \leq \underline 0$.
Notons aussi,
\begin{IEEEeqnarray*}{ccccc}
 D_{\underline f} & = & 
 \left(
\begin{array}{ccccc}
p^{q_\tau + f_{\sigma^{f-1}\tau}}\Id_r  & &  &  &   \\
 &p^{q_\tau + f_{\tau}}\Id_r  & & &    \\
 & & p^{q_\tau + f_{\sigma\tau}}\Id_r & & \\
 & & & \ddots  &    \\
  & & &   & p^{q_\tau+ f_{\sigma^{f-2}\tau}}\Id_r
\end{array}
\right)
\\
&=& \Diag(p^{q_\tau+f_{\sigma^{i-1}\tau}}\Id_r, i = 0 \dots f-1)  =   \Diag(p^{\min(q_\tau,q_{\sigma^{i-1}\tau})}\Id_r, i = 0 \dots f-1)
\end{IEEEeqnarray*}
dans une base adaptée à la décomposition $E_\tau \oplus E_{\sigma\tau} \oplus \dots \oplus E_{\sigma^{f-1}\tau}.$

\pro
L'application $m$ définit précédemment envoie $\Fil^{\underline f}\Lambda$ dans $\Fil^{\underline 0} \Lambda$.
De plus, $m$ envoie $\left(\Fil^{\underline f}\Lambda\right)^{\Phi = D_{\underline f}}$ dans $\left(\Fil^{\underline 0}\Lambda\right)^{\Phi = p^{q_\tau}}$.
\epro

\dem
Soit $j \dans \{0,\dots, f-1\}$. Comme $\phi^j(t_\mathcal O) \dans \Fil^1A$ si et seulement si $j = 0$, on a bien que chaque $x_{\tau'} \dans \Fil^{f_{\tau'}}\Lambda_{\tau'}$ est envoyé dans $\Fil^{f_{\tau'} + \max(0,q_\tau - q_{\tau'})} \Lambda_{\tau'} = \Fil^0\Lambda_{\tau'}$, d'où la première assertion.
Supposons maintenant que $x = \sum_{\tau'} x_{\tau'}$ vérifie $\Phi = D_{\underline f}$, on a donc que, pour tout $\tau'$,
\[ \Phi(x_{\tau'}) = p^{q_\tau + f_{\tau'}}x_{\sigma\tau'} = p^{\min(q_\tau,q_{\tau'})}x_{\sigma\tau'}.\]
Donc $t_\mathcal O^{\max(0,q_\tau - q_{\tau'})}\left( \prod_{j=1}^{f-1}\left(\frac{1}{p}\phi^j(t_\mathcal O)\right)^{\max(0,q_\tau - q_{\sigma^{-j}\tau'})}\right)x_{\tau'}$ vérifie,
\begin{eqnarray*}
\Phi\left(t_\mathcal O^{\max(0,q_\tau - q_{\tau'})}\left( \prod_{j=1}^{f-1}\left(\frac{1}{p}\phi^j(t_\mathcal O)\right)^{\max(0,q_\tau - q_{\sigma^{-j}\tau'})}\right)x_{\tau'}\right) 
\\= p^{\min(q_\tau,q_{\tau'})}
\left( \phi(t_\mathcal O)^{\max(0,q_\tau - q_{\tau'})}\right)\left( \prod_{j=1}^{f-1}\left(\frac{1}{p}\phi^{j+1}(t_\mathcal O)\right)^{\max(0,q_\tau - q_{\sigma^{-j}\tau'})}\right)x_{\sigma\tau'}\\
= p^{\min(q_\tau,q_{\tau'})} \phi(t_\mathcal O)^{\max(0,q_\tau - q_{\tau'})}\left( \prod_{j=2}^{f-1}\left(\frac{1}{p}\phi^j(t_\mathcal O)\right)^{\max(0,q_\tau - q_{\sigma^{-j}\sigma\tau'})}\right)\left(\frac{\phi^f(t_\mathcal O)}{p}\right)^{\max(0,q_\tau-q_{\sigma^{-f}\sigma\tau})}x_{\sigma\tau'} \\
= p^{\min(q_\tau,q_{\tau'})} t_\mathcal O^{\max(0,q_\tau-q_{\sigma\tau})}\phi(t_\mathcal O)^{\max(0,q_\tau - q_{\sigma^{-1}\sigma\tau'})}\left( \prod_{j=2}^{f-1}\left(\frac{1}{p}\phi^j(t_\mathcal O)\right)^{\max(0,q_\tau - q_{\sigma^{-j}\sigma\tau'})}\right)x_{\sigma\tau'} \\
= p^{q_\tau} t_\mathcal O^{\max(0,q_\tau-q_{\sigma\tau})}
\frac{\phi(t_\mathcal O)^{\max(0,q_\tau - q_{\sigma^{-1}\sigma\tau'})}}{p^{\max(0,q_\tau-q_{\tau'})}}\left( \prod_{j=2}^{f-1}\left(\frac{1}{p}\phi^j(t_\mathcal O)\right)^{\max(0,q_\tau - q_{\sigma^{-j}\sigma\tau'})}\right)x_{\sigma\tau'} \\
= p^{q_\tau}t_\mathcal O^{\max(0,q_\tau-q_{\sigma\tau})}\left(\prod_{j=1}^{f-1}\left(\frac{\phi^j(t_\mathcal O)}{p}\right)^{\max(0,q_\tau-q_{\sigma^{-j}\sigma\tau'})}\right)x_{\sigma\tau'}
\end{eqnarray*}
Donc $\sum_{\tau'} t_\mathcal O^{\max(0,q_\tau - q_{\tau'})}\left( \prod_{j=1}^{f-1}\left(\frac{1}{p}\phi^j(t_\mathcal O)\right)^{\max(0,q_\tau - q_{\sigma^{-j}\tau'})}\right)x_{\tau'}$ vérifie l'équation \[\Phi = p^{q_\tau}.\qedhere\]
\edem

On voudrait pouvoir relier $\left(\Fil^{\underline f}\Lambda\right)^{\Phi = D_{\underline f}}$ et $\left(\Fil^{\underline 0}\Lambda\right)^{\Phi = p^{q_\tau}}$ pour pouvoir caractériser 
l'image de $\alpha_{G,\tau,\eps}$, qui s'inclut dans la réduction modulo $\Fil^1A$ du second espace, en fonction de $\Ha_\tau(G)$, ce dernier apparaîssant naturellement quand on fait une 
méthode de Newton sur $\det(\omega_{G^D,\tau})^{\zeta_\tau = \id \otimes 1}$, qui est un quotient du premier espace.
L'idée est donc de montrer que $m$ est bijective, pour montrer que l'on peut relier l'image de l'application de Hodge-Tate à $\Ha_\tau(G)$.

\pro
L'application $m$ est injective.
\epro

\dem
C'est trivial puisque $A_{cris}$ est sans $t$-torsion et qu'à une unité (de $\ZZ_p^\times$) près,
\[ t = t_\mathcal O\prod_{j=1}^{f-1} \frac{1}{p}\phi^j(t_\mathcal O).\]
En effet, le terme de droite est dans $(\Fil^1A_{cris})^{\phi = p} = \ZZ_pt$, et sa valuation dans $\Gr^1A_{cris}$ est strictement plus petite que 1 (cf. \cite{FGL} page 130 ou les rappels de la section \ref{sect42}). 
\edem

On veut montrer que $m$ (plus précisément $m_\tau$) est surjective, et pour ça on va utiliser le théorème 5 de \cite{Fal}.

\subsection{Cristaux filtrés et théorème de Faltings}

On rappelle ici le coeur de l'article \cite{Fal}, qui va nous permettre de montrer que $m$ est surjective.

\defi
On définit la catégorie des cristaux filtrés admissibles d'amplitude $a$, notée $\mathcal{MF}_{[-a,0]}(\mathcal O_C)$ comme la catégorie des triplets,
\[ (M, \Fil^\cdot M, \Phi),\]
où $M$ est un $A_{cris}$-module filtré libre, les éléments $e_i$ d'une base ayant pour degrés $q_i$ avec $-a \leq q_i \leq 0$, $\Phi : M \fleche M$ est un endomorphisme $\phi$-linéaire
tel que la restriction de $\Phi$ à $\Fil^qM$ est divisible par $p^{q+a}$, vérifiant de plus la condition suivante,
\begin{equation}
\label{Ad}
\tag{Ad}
\sum_{i=-a}^0 \frac{\Phi}{p^{i+a}}(\Fil^i M) \text{ engendre sur } A_{cris} \text{ le module } M.
\end{equation}
\edefi

\rem
Strictement parlant, Faltings n'utilise pas des cristaux de la forme $\mathcal{MF}_{[-a,0]}(\mathcal O_C)$ mais plutôt, $\mathcal{MF}_{[0,a]}(\mathcal O_C)$, qui est la même catégorie, quitte à décaler la filtration (la version de Faltings est la "duale" de la notre : ses foncteurs sont contravariants).
\erem

\exe
Si $G$ est un groupe $p$-divisible, $E$ son cristal,
\[ (E, \Fil E, \Phi) \dans \mathcal{MF}_{[-1,0]}(\mathcal O_C)\]
\eexe

On a alors le théorème suivant (voir aussi \cite{Chen} dans le cas $a=1$ des groupes $p$-divisibles) :

\thr[Faltings]
\label{thrfal}
Supposons que $a \leq p-2$. Soit $M \dans \mathcal{MF}_{[-a,0]}(\mathcal O_C)$. On pose,
\[ \mathbb D(M) = (\Fil^0 M)^{\Phi = p^a}.\]
Alors,
\begin{enumerate}
\item $\mathbb D(M)$ est un $\ZZ_p$-module libre de rang égal au rang de $M$ sur $A_{cris}$.
\item On a les inclusions, strictement compatibles aux filtrations, 
\[t^a M \subset \mathbb D(M) \otimes_{\ZZ_p} A_{cris}  \subset M.\]
\end{enumerate}
\ethr

\dem[esquisse]
L'idée est de fixer une base de $M$ sur $A_{cris}$, et de se ramener à des équations. 
Comme sur $\Fil^0M$, $\Phi$ est divisible par $p^a$, on regarde la matrice de $\Phi_a = \frac{\Phi}{p^a}$, et on obtient des équations pour,
\[(\Fil^0M)^{\Phi = p^a} = (\Fil^0M)^{\Phi_a = \id},\]
et on réduit ces dernières équations modulo $p$, c'est-à-dire que l'on en cherche les équations dans $A_{cris}/p$.
Comme $a < p-1$, Faltings utilise que (cf. \cite{Fal} p127, \cite{Fal2} page 37, ou \cite{Chen} lemme 4.7 (qui s'adapte au cas général)),
\lem
\label{lemfilp}
L'application de réduction modulo $(p,\Fil^p)$ induit un isomorphisme,
\[ \left(\Fil^0 (M\otimes A_{cris}/p)\right)^{\Phi_a = \id} \overset{\sim}{\fleche} \left(\Fil^0 (M\otimes \quotient{A_{cris}}{(p,\Fil^pA_{cris})})\right)^{\Phi_a = \id}.\]
\elem

Grâce à la proposition ci-après, on est donc amené à résoudre des équations dans $\mathcal O_C/p$, qui se relèvent uniquement dans $\mathcal O_C$ par une méthode de Newton, et on peut en calculer
le nombre, qui est exactement $p^{\rg_A M}$.
Faltings montre ensuite que toute solution de $\left(\Fil^0 (M\otimes A_{cris}/p)\right)^{\Phi_a = \id}$ se relève uniquement à $A_{cris}$, par un lemme de Nakayama topologique 
.
\edem

\rem
On pourra trouver dans le cas des groupes $p$-divisibles une démonstration plus détaillée (et sur une base plus générale) dans \cite{Chen}.
\erem

\pro
On peut filtrer $\mathcal O_C/p\mathcal O_C$ par $\Fil^i \mathcal O_C/p\mathcal O_C = p^{\frac{i}{p}}\mathcal O_C/p\mathcal O_C$ et on définie,
\[ \theta : 
\begin{array}{ccc}
  \mathcal O_C/p\mathcal O_C&\fleche  & \mathcal O_C/p\mathcal O_C  \\
 x & \longmapsto  & x^p  .
\end{array}
\]
On a alors un isomorphisme, compatible aux applications $\theta$ et qui identifie les filtrations,
\[ A_{cris}/(p,\Fil^pA_{cris}) \overset{\sim}{\fleche} \mathcal O_C/p\mathcal O_C.\]
\epro

Pour la preuve, on renvoie par exemple à \cite{Chen} Lemme 4.8. Donnons peut-être les images des périodes définies dans la section \ref{sect42} sous cet isomorphisme :
Comme $t,t_\mathcal O$ sont dans le $\Fil^1A_{cris}$, grâce aux valuations rappelées dans la section \ref{sect42}, leurs images modulo $(p,\Fil^pA)$ sous l'isomorphisme précédent
sont des éléments de valuation $\frac{1}{p} + \frac{1}{p-1}$ et $\frac{1}{p} + \frac{1}{p^f-1}$ respectivement. La valuation, pour $i \dans\{1,\dots,f-1\}$, de l'image de
$\frac{1}{p}\phi^i(t_\mathcal O)$ est $\frac{p^i}{p^f-1}$ : il n'est pas dans $\Fil^1A_{cris}$.

\rem
Supposons que $M$ soit un cristal filtré dans $\mathcal{MF}_{[-a,0]}(\mathcal O_C)$ muni d'une action de $\mathcal O$.
Donc $\mathbb D(M)$ est naturellement un $\mathcal O$-module, libre, car sans torsion.
En effet, l'action de $\mathcal O$ est supposée commuter à $\Phi$,
en particulier $(\Fil^r M)^{\Phi = p^a}$ est muni d'une action de $\mathcal O$.
\erem

\subsection{Modification du cristal $\Lambda$}

Dans toute cette section on fait l'hypothèse que $q_\tau < p-1$.
On a construit un cristal filtré admissible,
\[ \underline{\Lambda} = (\Lambda, \Fil^{\cdot}\Lambda,\Phi) \dans \mathcal{MF}_{[-q_\tau,0]}(\mathcal O_C),\]
(en modifiant la filtration par une $\ZZ$-filtration en posant $\Fil^k \Lambda = \Fil^{\underline k} \Lambda$, pour tout $k \dans \ZZ$), l'admissibilité découlant moralement du fait que $E$ étant admissible, ses puissances extérieures le sont aussi (on peut le vérifier directement facilement).
Le théorème de Faltings nous dit alors que si $q_\tau < p-1$, alors,
\[ \mathbb D(\underline \Lambda) = (\Fil^{\underline 0})^{\Phi = p^{q_\tau}},\]
est un $\ZZ_p$-module libre de rang $f\binom{h}{q_\tau} = \rg_{A_{cris}} \Lambda =: fr$.

\pro
\label{pro414}
On a l'égalité suivante dans $\Lambda$,
\[ \bigwedge^{q_\tau}_{A_{cris} \otimes \mathcal O} (\Fil^0 E)^{\Phi = p} = (\Fil^{\underline0}\Lambda)^{\Phi = p^{q_\tau}}.\]
\epro

\dem
On a une inclusion naturelle évidente \[\bigwedge^{q_\tau}_{\mathcal O}(\Fil^0E)^{\Phi = p} \subset (\Fil^{\underline 0} \Lambda)^{\Phi = p^{q_\tau}} \subset \Lambda,\] simplement 
par définition de $\Fil^{\underline 0} \Lambda$.

De plus le théorème de Faltings, appliqué à $E$ et à $\Lambda$, nous dit que les deux modules ont le même rang.
Il suffit donc de voir que le conoyau (qui est automatiquement de type fini) est nul. De plus les deux modules ont le même nombre d'élément modulo $p$, on peut donc réduire modulo $p$ et vérifier que le conoyau est nul, i.e. que l'inclusion reste une inclusion modulo $p$.
Or $\bigwedge^{q_\tau}_{\mathcal O} (\Fil^0 E)^{\Phi = p}\otimes_{\mathcal O}A_{cris} \supset t^{q_\tau}\Lambda$, dont l'image dans $\Lambda/p\Lambda$ est de rang le rang de $\Lambda$ sur $A_{cris}$ ($t^{q_\tau} \neq 0 \pmod p$ car $q_\tau < p-1$).
\edem

On rappelle qu'on veut montrer que,
\[ m : \left(\Fil^{\underline f}\Lambda\right)^{\Phi = D_{\underline f}} \fleche \left(\Fil^{\underline 0}\Lambda\right)^{\Phi = p^{q_\tau}},\]
est surjective. On sait déjà qu'elle est injective, on va donc d'abord montrer que les deux $\mathcal O$-modules ont le même rang.

\defi
On modifie le cristal $\underline \Lambda$ en,
\[ \underline{\Lambda'} = (\Lambda, (\Fil^{\underline{f + r}} \Lambda)_{r \in \ZZ}, p^{q_\tau}D_{\underline f}^{-1}\Phi),\]
qui est bien définie, car $D |p^{q_\tau}I_{fr}$.
\edefi

\subsection{Admissibilité}

Rappelons que $(E,\Fil E, \Phi)$ est un $\mathcal O\otimes_{\ZZ_p}A_{cris}$-module filtré, que l'on peut décomposer,
\[ E = \bigoplus_{\tau'} E_{\tau'} \quad \text{et}\quad \Fil E = \bigoplus_{\tau'} \Fil E_{\tau'}.\]
De plus, il existe une base, pour tout $\tau'$, $(e_1^{\tau'},\dots,e_h^{\tau'})$ de $E_{\tau'}$ telle que,
\[ \Fil E_{\tau'} = (e_1^{\tau'},\dots,e_{p_{\tau'}}^{\tau'})\Fil^1A + (e_{p_{\tau'}+1}^{\tau'},\dots,e_h^{\tau'})A,\]
et donc la filtration supérieure, donnée par $\Fil^iE = (\Fil^iA)\Fil E + (\Fil^{i+1}A)E$ s'écrit simplement,
\[\Fil^i E_{\tau'} = (e_1^{\tau'},\dots,e_{p_{\tau'}}^{\tau'})\Fil^{i+1}A + (e_{p_{\tau'}+1}^{\tau'},\dots,e_h^{\tau'})\Fil^iA,\]
Notons pour tout $\tau'$,
\[ W_{\tau'} = (e_{p_{\tau'}+1}^{\tau'},\dots,e_h^{\tau'})A \subset \Fil E_{\tau'}, \quad \text{et} \quad W = \bigoplus_{\tau'} W_{\tau'}.\]
$W$ est un relèvement par $\theta$ de $\omega_{G^D}$, mais il n'est a priori pas stable sous $\Phi$. On a alors,
\begin{equation}
\label{filW}
\Fil^i E_{\tau'} =\Fil^{i+1}E_{\tau'} + (\Fil^iA)W_{\tau'}.\end{equation}

Rappelons que l'on a fixé un $\tau$ tel que $q_\tau \not\in\{0,h\}$ et 
\[\Lambda = \bigwedge_{\mathcal O\otimes A}^{q_\tau} E = \bigoplus_{\tau'} \Lambda_{\tau'},\] 
dont la filtration est donnée par la "convolution" de celle de $E$.

\pro
\label{profil}
Soit $\tau'$ tel que $q_{\tau'} < q_\tau$. Soit $0 \geq s > -q_\tau + q_{\tau'} = f_{\tau'}$, alors,
\begin{IEEEeqnarray*}{rcl} 
\Fil^s \Lambda_{\tau'} &=& \sum_{j = s + q_\tau-q_{\tau'}}^{s+q_\tau} (\Fil^{j}A) \Fil^{s-j}\Lambda_{\tau'} \\
&=& (\Fil^{s+q_\tau}A)\Lambda_{\tau'} + (\Fil^{s+q_\tau-1}A)\Fil^{-q_\tau + 1} \Lambda_{\tau'} + \dots + (\Fil^{s+q_\tau-q_{\tau'}}A)\Fil^{-q_\tau+q_{\tau'}}\Lambda_{\tau'}.\end{IEEEeqnarray*}
Autrement dit, passée la dimension de $\omega_{G^D,\tau'}$, la filtration n'est modifiée que par l'action des scalaires ($\Fil^.A$).
\epro

\dem
Il suffit d'écrire \[\Fil^iE = (\Fil^i A)W + \Fil^{i+1}E,\]
et,
\begin{equation}
\label{annW}
 \bigwedge^r W_{\tau'} = 0, \forall r > q_{\tau'}.\end{equation}
Écrivons la définition de $\Fil^s \Lambda_{\tau'}$ comme l'image d'une somme de produits tensoriels, 
\[ \Fil^s \Lambda_{\tau'} = \im\left( \sum_{
\begin{array}{c}
 i_1,\dots, i_{q_\tau} \\
 i_1 +\dots + i_{q_\tau} = s
\end{array}}
\Fil^{i_1}E_{\tau'} \otimes \dots\otimes\Fil^{i_{q_\tau}}E_{\tau'} \fleche \Lambda_{\tau'}\right),\]
et concentrons nous sur l'image de $\Fil^{i_1}E_{\tau'} \otimes \dots\otimes\Fil^{i_{q_\tau}}E_{\tau'}$, que l'on peut décomposer grâce à (\ref{filW}),
\[\Fil^{i_1}E_{\tau'} \otimes \dots\otimes\Fil^{i_{q_\tau}}E_{\tau'} = \sum_{j=0}^{q_{\tau}} (\Fil^{i_1}A)W_{\tau'} \otimes\dots\otimes (\Fil^{i_j}A)W_{\tau'} \otimes (\Fil^{i_{j+1}+1}A)E_{\tau'}
\otimes \dots \otimes (\Fil^{i_{q_\tau}+1}A)E_{\tau'}.\]
Or, si $j > q_{\tau'}$, l'image de $(\Fil^{i_1}A)W_{\tau'} \otimes\dots\otimes (\Fil^{i_j}A)W_{\tau'} \otimes (\Fil^{i_{j+1}+1}A)E_{\tau'}
\otimes \dots \otimes (\Fil^{i_{q_\tau}+1}A)E_{\tau'}$ dans $\Lambda_{\tau'}$ est nulle en vertu de (\ref{annW}), on peut donc écrire,
\begin{eqnarray*}
\Fil^{i_1}E_{\tau'} \otimes \dots\otimes\Fil^{i_{q_\tau}}E_{\tau'} \\
= \sum_{j=0}^{q_{\tau'}} (\Fil^{i_1}A)W_{\tau'} \otimes\dots\otimes (\Fil^{i_j}A)W_{\tau'} \otimes (\Fil^{i_{j+1}+1}A)E_{\tau'}
\otimes \dots \otimes (\Fil^{i_{q_\tau}+1}A)E_{\tau'}\\
 \subset \sum_{j=0}^{q_{\tau'}} (\Fil^{q_\tau - j + i_1 + \dots + i_j + i_{j+1} + \dots + i_{q_\tau}}A)\underbrace{W_{\tau'} \otimes\dots\otimes W_{\tau'}}_{j}\otimes E_{\tau'}
\otimes \dots \otimes E_{\tau'},
\end{eqnarray*}
dont l'image dans $\Lambda_{\tau'}$ est par définition incluse dans,
\[ \sum_{j=0}^{q_{\tau'}} (\Fil^{q_\tau-j + s}A) \Fil^{-q_\tau + j}\Lambda_{\tau'} \qedhere.\]
\edem

\pro
\label{prophifil}
Soit $\Lambda$ le cristal précédent.
Alors pour $r+s \leq (p-1) - q_\tau$, les sous-modules engendrés sur $A_{cris}$ par
\[\Phi((\Fil^rA)(\Fil^s\Lambda_{\tau'})) \quad \text{et} \quad \Phi(\Fil^{r+s}\Lambda_{\tau'}),\]
sont les mêmes.
\epro

\dem
C'est vrai pour $A_{cris}$ : l'espace engendré par $\phi(\Fil^iA_{cris})$ est $p^iA_{cris}$ pour tout $i \leq p-1$ (cf. \cite{Chen} Lemme 2.8 par exemple). Utilisons $<S>$ pour dénoter le sous-module engendré sur $A_{cris}$ par $S$. 
Alors \[<\phi(\Fil^rA\Fil^sA)> = <\phi(\Fil^rA)><\phi(\Fil^sA)> = p^rAp^sA = p^{r+s}A = <\phi(\Fil^{r+s}A)>,\]
pour tout $r,s \geq 0$ tels que $r+s \leq p-1$.
Soit $E$ le cristal d'un $\mathcal O$-module $p$-divisible et,
\[ \Fil^iE_{\tau'} = \Fil^iW_{\tau'} + \Fil^{i+1}E_{\tau'},\]
donc, grâce au cas de $A_{cris}$,
\begin{eqnarray*} 
<\Phi(\Fil^rA\Fil^sE_{\tau'})> = <\Phi(\Fil^{r}A(\Fil^sAW_{\tau'}) + \Fil^{r}A(\Fil^{s+1}AE_{\tau'}))>\\
= <\Phi((\Fil^{r+s}A)W_{\tau'} + (\Fil^{r+s+1}A)E_{\tau'})>=<\Phi(\Fil^{r+s}E_{\tau'})>.\end{eqnarray*}
Et donc grâce à la définition sur la filtration sur les puissances extérieures, on a encore le résultat sur $\Lambda_{\tau'}$.
\edem

\rem
Bien sûr c'est faux sans appliquer $\Phi$. Par exemple, sur $A_{cris}$, $\Fil^{r+s}A_{cris} \supset \Fil^rA_{cris}\Fil^sA_{cris}$, mais l'inclusion est stricte, même pour $r+s \leq p-1$.
\erem

Ces deux dernières propositions vont nous permettre de vérifier que le cristal filtré modifié $\underline{\Lambda}'$ est bien admissible.

\pro
Le cristal filtré $\underline{\Lambda'}$ appartient à $\mathcal{MF}_{[-q_\tau,0]}(\mathcal O_C)$.
\epro

\dem
Il faut vérifier que la restriction de $\Phi' = p^{q_\tau}D_{\underline f}^{-1}\Phi$ à $\Fil^{\underline{f + r}}\Lambda$ est divisible par $p^{q_\tau + r}$ pour $r \dans [-q_\tau,0]\cap \ZZ$,
et la condition d'admissibilité.
Soit $-q_\tau \leq r \leq 0$, et $i_1,\dots,i_{q_\tau} \geq -1$ tels que $i_1 + \dots + i_{q_\tau} = f_{\tau'} + r$.
Alors restreint à \[\Fil^{i_1}E_{\tau'} \otimes\dots\otimes \Fil^{i_{q_\tau}}E_{\tau'},\] $\Phi \otimes \dots \otimes \Phi$ est divisible par 
\[p^{i_1 + 1 + i_2 + 1 + \dots + i_{q_\tau} + 1} = p^{q_\tau+f_{\tau'}+r},\]
donc $\Phi'_{|\Fil^{f_\tau' + r}\Lambda_{\tau'}}$ est divisible par,
\[ p^{q_\tau - \min(q_\tau,q_{\tau'}) + q_\tau + \min(0,q_{\tau'} -q_\tau) + r} = p^{\max(0,q_\tau-q_{\tau'}) + q_{\tau} + \min(0,q_{\tau'} - q_\tau) + r}  = p^{q_\tau + r}.\]
Il faut ensuite vérifier la condition (\ref{Ad}).
Il faut donc voir que $\Lambda$ est engendré sur $A$ par 
\[\sum_{ i = -q_\tau}^{0} \frac{\Phi'}{p^{i+q_\tau}}(\Fil^{\underline i} \Lambda').\]
Mais en fait $\Lambda$ est engendré par $\frac{\Phi}{p^{q_\tau}}(\Fil^{\underline 0}\Lambda)$ sur $A_{cris}$, puisque c'est déjà le cas pour $E$.
Mais alors, pour tout $\tau'$, on a par la proposition (\ref{profil}),
\[\Fil^0{\Lambda_{\tau'}} = (\Fil^{q_\tau}A)\Lambda_{\tau'} + (\Fil^{q_\tau-1}A)\Fil^{-q_\tau+1}\Lambda_{\tau'} + \dots + (\Fil^{q_\tau-q_{\tau'}}A)\Fil^{-q+q_{\tau'}}\Lambda_{\tau'},\]
et donc en regardant l'espace engendré par $\Phi$, d'après la proposition (\ref{prophifil}), l'espace engendré par l'image par $\Phi$ de $\Fil^0{\Lambda_{\tau'}}$ est, si $q_{\tau'}<q_\tau$,
\begin{eqnarray*} 
p^{q_\tau}\Lambda_{\sigma\tau'} = \\ <\phi((\Fil^{q_\tau}A))\Phi(\Lambda_{\tau'}) + \phi(\Fil^{q_\tau-1}A)\Phi(\Fil^{-q_\tau+1} \Lambda_{\tau'})+ \dots + \phi(\Fil^{q_\tau-q_{\tau'}}A) \Phi(\Fil^{-q_{\tau}+q_{\tau'}}\Lambda_{\tau'})>\\
 = p^{q_\tau-q_\tau'}<\Phi(\Fil^{-q_\tau+q_{\tau'}}\Lambda_{\tau'})>
\end{eqnarray*}
Et si $q_\tau \leq q_{\tau'}$, alors $\Phi(\Fil^0\Lambda_{\tau'})$ engendre $p^{q_\tau}\Lambda_{\sigma\tau'}$. 

Autrement dit, dans tous les cas, $\Phi(\Fil^{f_{\tau'}}\Lambda_{\tau'})$ engendre 
\[p^{\min(q_\tau,q_{\tau'})}\Lambda_{\sigma\tau'}.\]
On en déduit que $D_{\underline f}^{-1}\Phi(\Fil^{\underline f}(\Lambda)) = \frac{\Phi'}{p^{q_\tau}}(\Fil^{\underline 0}\Lambda')$ engendre $\Lambda$, et donc l'admissibilité.
\edem

\subsection{Division}

On a montré dans la section précédente que le cristal modifié $\underline{\Lambda'}$ est dans $\mathcal{MF}_{[-q_\tau,0]}(\mathcal O_C)$, et il en est de même de $\underline{\Lambda}$.

En appliquant le théorème (\ref{thrfal}), on en déduit le résultat suivant,
\pro
Supposons $q_\tau \leq p-2$. 
Le $\mathcal O$-module $\mathbb D(\underline{\Lambda'}) = \left(\Fil^{0}\Lambda'\right)^{\Phi' = p^{q_\tau}}=
\left(\Fil^{\underline f}\Lambda\right)^{\Phi = D_{\underline f}}$ est libre de rang $\binom{q_\tau}{h}$ sur $\mathcal O$.
On a de plus les inclusions strictement compatibles aux filtrations,
\[ t^{q_\tau}\Lambda \subset \mathbb D(\underline{\Lambda'})\otimes_{\ZZ_p}A_{cris} \subset \Lambda.\]
\epro

\pro
Soit $(M,\Fil M,\Phi)$ un cristal filtré dans $\mathcal{MF}_{[-a,0]}(\mathcal O_C)$, muni d'une $\mathcal O$-action, c'est à dire,
\[ M = \bigoplus_{\tau} M_\tau, \quad \text{et} \quad \Phi : M_\tau \fleche M_{\sigma\tau} \quad \forall \tau.\]
Alors pour tout $\tau$ la projection,
\[\pi_\tau : M \fleche M_\tau,\]
induit un isomorphisme,
\[ \pi_\tau : \mathbb D(M) =\left(\Fil^0M\right)^{\Phi = p^a} \overset{\sim}{\fleche} \pi_\tau(\left(\Fil^0M\right)^{\Phi = p^a}).\]
En particulier ces deux $\mathcal O$-modules ont le même rang, $\rg_{A_{cris}} M$.
\epro

\dem
On a deux $\mathcal O$-modules de rang fini, et une surjection,
\[ \pi_\tau : \mathbb D(M) =\left(\Fil^0M\right)^{\Phi = p^a} \overset{}{\fleche} \pi_\tau(\left(\Fil^0M\right)^{\Phi = p}).\]
Soit $x = \sum_\tau x_\tau$ tel que $\Phi(x) = p^ax$, c'est à dire,
\[\forall j \geq 0, \Phi^j(x_\tau) = p^{ja}x_{\sigma^{j}\tau}.\]
Autrement dit, comme $A_{cris}$ est sans $p$-torsion, $x_{\sigma^j\tau}$ est entièrement déterminé par $x_\tau$.
De plus, il faudrait vérifier qu'un tel $x_{\sigma^j\tau}$ est dans $\Fil^0M$, mais c'est automatique, car il est unique, et $x_\tau$ provient (par surjectivité et unicité de l'antécédent) de 
$x = \sum_{j=0}^{f-1} x_{\sigma^j\tau} \dans \Fil^0M$.
\edem

On en déduit le théorème principal suivant, dit de 'suppression de période', ou de division, qui va nous permettre de relier l'image de $\alpha_{G,\tau}$ à l'invariant de Hasse.

\thr
L'application,
\[ m : \left(\Fil^{\underline f}\Lambda\right)^{\Phi = D_{\underline f}} \fleche \left(\Fil^{\underline 0} \Lambda\right)^{\Phi = p},\]
est un isomorphisme.
\ethr

\dem
En effet, on sait que $m$ est une application $\mathcal O$-linéaire, injective, entre deux $\mathcal O$-modules libres de même rang, on en déduit que son conoyau est fini.
Par la proposition précédente, il suffit de vérifier que,
\[ m_\tau : \pi_\tau(\left(\Fil^{\underline f}\Lambda\right)^{\Phi = D_{\underline f}}) \overset{\sim}{\fleche} \pi_\tau(\left(\Fil^{\underline 0} \Lambda\right)^{\Phi = p}),\]
est un isomorphisme, c'est-à-dire que son conoyau est nul.
En particulier, il suffit de voir que modulo $p$, $m_\tau$ est un isomorphisme.
Or d'après Faltings et le lemme (\ref{lemfilp}), il suffit de regarder ce qu'il se passe modulo $(p,\Fil^pA_{cris})$, et les deux cristaux modulo $(p,\Fil^pA)$ ont le même 
nombre de solutions à leurs équations respectives. Il suffit donc seulement de voir que $m_\tau$ est injective. Et dans ce cas on a la description de $A_{cris}/(p,\Fil^pA_{cris})$, 
et donc on peut calculer que sous l'isomorphisme, $\mathcal O_C/p = A_{cris}/(p,\Fil^pA_{cris})$, $m_\tau$ est donnée par multiplication par un élément de valuation,
\[ \sum_{j=1}^{f-1} \frac{p^j\max(0,q_\tau - q_{\sigma^{-j}\tau})}{p(p^f-1)} \leq \frac{q_\tau}{p(p-1)} - \frac{q_\tau}{p(p^f-1)} < \frac{1}{p},\]
et on note,
\[K_\tau = \sum_{j=1}^{f-1} \frac{p^j\max(0,q_\tau - q_{\sigma^{-j}\tau})}{p^f-1}.\]
Comme $q_\tau \leq p-2$, on a que $K_\tau < 1$.
Or on sait que,
\[ t^{q_\tau} \Lambda \subset \left(\Fil^{\underline f}\Lambda\right)^{\Phi = D_{\underline f}}\otimes_{\ZZ_p}A_{cris},\]
donc si on choisit un des générateurs $e$ sur $\FP$ de $\pi_\tau(\left(\Fil^{\underline f}\Lambda\right)^{\Phi = D_{\underline f}}) \pmod {p,\Fil^pA}$, et on choisit un isomorphisme,
\[\Lambda_\tau \simeq A_{cris}^h,\]
alors au moins une des coordonnées de $e$ divise $t^{q_\tau} \pmod{p,\Fil^pA}$.
Or, \[m_\tau (0,\dots,0,t^{q_\tau},0,\dots,0) \pmod{p,\Fil^pA_{cris}} = (0,\dots,0,p^{\frac{K_\tau}{p} + \frac{q_\tau}{p(p-1)} + \frac{q_\tau}{p}},0,\dots,0).\]
Mais en utilisant que $q_\tau \leq p-2$, on en déduit,
\[\frac{K_\tau}{p} + \frac{q_\tau}{p(p-1)} + \frac{q_\tau}{p} < \frac{1}{p} + \frac{p-2}{p} + \frac{1}{p} = 1,\]
donc que $m_\tau(e) \neq 0 \pmod{p,\Fil^pA_{cris}}$.
On en déduit que $m_\tau$ est injective modulo $(p,\Fil^pA_{cris})$, et donc modulo $p$ par le lemme (\ref{lemfilp}), et comme les deux espaces ont le même nombre (fini) de solutions
modulo $p$, on en déduit que $m_\tau \pmod p$ est bijective.
Donc $m_\tau$ est un isomorphisme.
\edem

\subsection{Reconstruction de l'invariant de Hasse}
\label{sectmuha}

Dans \cite{Her1}, on a construit les invariants de Hasse partiels, en supposant le groupe $p$-divisible (tronqué) $G$ sur une base lisse de caractéristique $\overline S$, 
et en regardant donc le cristal dans $\Cris(\overline S/\Spec(\ZZ_p))$,
\[ \mathbb D(G) = \mathcal Ext^1(G^D,\mathcal O_{S/\Sigma}) = \bigoplus \mathbb D(G)_\tau,\]
muni des applications $V$ et $F$. Plus précisément, on considère,
\[ \bigwedge^{q_\tau} \mathbb D(G)_\tau,\]
sur lequel l'application,
\[ V^f : \bigwedge^{q_\tau} \mathbb D(G)_\tau \fleche \bigwedge^{q_\tau} \mathbb D(G)_\tau^{(p^f)},\]
est divisible par $p^{k_\tau}$, division que l'on note $\phi_\tau$. Pour montrer ce fait, on utilise que la base est lisse, pour localement relever 
$\overline{S}$ en $S/\Spec(\ZZ_p)$ lisse, et utiliser la description du cristal 
en terme de module à connection, auquel cas il est clair que l'on peut diviser. On montre de plus que cette application se factorise à l'arrivé par un certain sous-faisceau
 (\cite{Her1} Théorème 3.7) et qu'elle est unique, et quelle induit une application 
\[\widetilde{\Ha_\tau}(G) : \det(\omega_{G^D[p],\tau}) \fleche \det(\omega_{G^D[p],\tau})^{\otimes(p^f)},\] cf \cite{Her1} Proposition 3.13, Lemme 3.14. 
 Couplé au fait que le champ des $\mathcal{BT}_r^{\mathcal O}$ est lisse, on peut donc faire la construction sur toute base.

Néanmoins étant donné un $\mathcal O$-module $p$-divisible $G$ (non tronqué pour simplifier) sur $\mathcal O_C/p$, qui n'est donc pas lisse, a priori l'expression de $\phi_\tau$ sur son cristal, bien qu'elle existe, semble un peu compliquée.
Mais en fait, si $E = \bigoplus_\tau E_\tau$ est l'évaluation de $\mathbb D(G)$ sur $A_{cris} \twoheadrightarrow \mathcal O_C/p$, on peut reconstruire $\phi_\tau$. En effet,
le principal ingrédient est que, 
\[ V(E_{\tau'}) \subset \Fil (E_{\sigma^{-1}\tau'}^{(\phi)}) + pE_{\sigma^{-1}\tau'}^{(\phi)},\]
et si $q_{\tau'} < q_\tau$, alors,
\[\im(\bigotimes^{q_\tau} \Fil (E_{\tau'}^{(\phi)}) \fleche (\bigwedge^{q_\tau} E_{\tau'})^{(\phi)}) \subset 
\left(\Fil^{q_\tau-q_{\tau'}}A_{cris}\bigwedge^{q_\tau} E_{\tau'}\right)^{(\phi)}\subset 
\bigwedge^{q_\tau} E_{\tau'} \otimes_{A_{cris},\phi} p^{q_\tau-q_{\tau'}}A_{cris}.\]
La dernière inclusion étant due au fait que si $M$ est un $A_{cris}-$module, si $z \dans \Fil^iA_{cris}$, $x \dans M$, alors dans $M^{(\phi)} = M \otimes_{A_{cris},\phi} A_{cris}$,
\[ (zx) \otimes 1 = x \otimes \phi(z), \quad \text{et} \quad \phi(z) \dans p^i A_{cris}.\]
On peut donc, pour chaque $j$ tel que $q_{\sigma^j\tau} < q_\tau$ diviser $V^j$ par $p^{q_\tau - q_{\sigma^j\tau}}$, et donc il existe,
\begin{equation}
\label{zetaOC}
 \phi_\tau' : \bigwedge^q_\tau E_\tau \fleche  \bigwedge^q_\tau E_\tau\otimes_{A_{cris},\phi^f} A_{cris},\end{equation}
tel que $p^{k_\tau}\phi_\tau' = V^f$. De plus comme $A_{cris}$ est sans $p$-torsion, un tel $\phi_\tau'$ est unique.
Si on note $G^{univ}/\mathcal{BT}_r^\mathcal O$ le groupe $p$-divisible universel sur le champ des $\mathcal {BT}_r^\mathcal O$, alors il existe un unique morphisme
de cristaux modulo $p^{r-k_\tau}$, 
\[ \phi_\tau : \bigwedge^q_\tau \mathcal E_\tau \fleche  \left(\bigwedge^q_\tau \mathcal E_\tau\right)^{(\phi^f)},\]
qui vérifie que $p^{k_\tau} \phi_\tau = V^f$, et donc si $G/\mathcal O_C$ vérifie $G[p^r] = x^*G^{univ}$ pour $x \dans \mathcal{BT}_r^\mathcal O(\mathcal O_C)$, le morphisme (\ref{zetaOC}) coïncide avec 
l'invariant de Hasse associé à $\tau$ modulo $p^{r-k_\tau}$.

\pro
La réduction modulo $(\ker\theta,p)$ de $\pi_\tau\left(\left(\Fil^{\underline f}\Lambda\right)^{\Phi = D_{\underline f}}\right)$ est incluse dans,
\[\{ x \dans \det\omega_{G[p]^D,\tau} : \widetilde{\Ha}_\tau(x) = x \otimes 1\}.\]
\epro

\dem
En effet, soit $x \dans \left(\Fil^{\underline f}\Lambda\right)^{\Phi = D_{\underline f}}$ alors il vérifie que \begin{equation}\label{Phif}\Phi^f(x_\tau) = p^{fq_\tau -k_\tau}\cdot x\otimes 1.\end{equation}
Rappelons que,
\[ V : \Lambda \fleche \Lambda^{(\phi)}, \quad \text{et} \quad F: \Lambda^{(\phi)} \fleche \Lambda,\]
avec $VF = FV = p^{q_\tau}$ et $\Phi(x) = F(x\otimes 1)$. On en déduit que $V^f(\Phi^f(x)) = p^{fq_\tau}x\otimes 1$ et par (\ref{Phif}), que $V^f(x) = p^{k_\tau} x\otimes 1$.
Comme $A_{cris}$ est sans $p$-torsion (et $\bigwedge^{q_\tau} E_\tau$ est libre), on en déduit que $\phi_\tau'(x) = x \otimes 1$. Le résultat s'en déduit par réduction modulo $(\ker\theta,p)$.
\edem

\subsection{Image de l'application de Hodge-Tate}

Essayons maintenant de relier l'image de Hodge-Tate à $\left(\Fil^{\underline f}\Lambda\right)^{\Phi = D_{\underline f}}$. On peut calculer explicitement, grâce à l'expression de 
$m_\tau$ et aux valuations rappelées dans la section (\ref{sect42}), que modulo $\Fil^1A_{cris} = \ker\theta$, $m_\tau$ est donné par multiplication par un élément de valuation $K_\tau$, où 
\[K_\tau = \sum_{j=1}^{f-1} \frac{p^j\max(0,q_\tau - q_{\sigma^{-j}\tau})}{p^f-1}.\]
Comme $q_\tau \leq p-2$, on a que $K_\tau < 1$.
On a donc le diagramme commutatif suivant, où $u$ est inversible dans $\mathcal O_C$,

\begin{center}
\begin{tikzpicture}[description/.style={fill=white,inner sep=2pt}] 
\matrix (m) [matrix of math nodes, row sep=3em, column sep=2.5em, text height=1.5ex, text depth=0.25ex] at (0,0)
{ 
\left(\Fil^{\underline f}\Lambda\right)^{\Phi = D_{\underline f}} & &\left(\Fil^{\underline 0} \Lambda\right)^{\Phi = p} \\
\bigwedge^{q_\tau} \omega_{G^D,\tau} & & \bigwedge^{q_\tau} \omega_{G^D,\tau} \\
 };

\path[->,font=\scriptsize] 
(m-1-1) edge node[auto] {$m$} (m-1-3)
(m-1-1) edge node[auto] {$\theta\circ\pi_\tau$} (m-2-1)
(m-1-3) edge node[auto] {$\theta\circ\pi_\tau$} (m-2-3)
(m-2-1) edge node[auto] {$up^{K_\tau}$} (m-2-3);
\end{tikzpicture}
\end{center}

De plus ce diagramme peut se factoriser par la flèche inversible $m_\tau$, et sous l'identification 
$\left(\Fil^{\underline 0} \Lambda\right)^{\Phi = p} = \bigwedge^{q_\tau}_\mathcal OT_pG$ de la proposition (\ref{pro414}), on en déduit que l'image de 
$\bigwedge^{q_\tau} \alpha_{G,\tau}$ se déduit de celle de l'image de la réduction modulo $\Fil^1A_{cris}$ de 
$\pi_\tau(\left(\Fil^{\underline f}\Lambda\right)^{\Phi = D_{\underline f}})$ par multiplication par $up^{K_\tau}$.

Or les éléments de $\pi_\tau(\left(\Fil^{\underline f}\Lambda\right)^{\Phi = D_{\underline f}})$ vérifient l'équation,

\[ \Phi^f = p^{fq_\tau - k_\tau}, \quad \text{où} \quad k_\tau = \sum_{j=0}^{f-1} \max(0,q_\tau - q_{\sigma^{j}\tau}).\]

C'est-à-dire l'équation $V^f = p^{k_\tau}\otimes 1$ ou encore, en notant $\phi_\tau$ la division de $V^f$ sur $\Lambda_\tau$ par $p^{k_\tau}$, construite dans \cite{Her1} et la section précédente,
\[\phi_\tau = \id.\]

Autrement dit, d'après la proposition précédente, comme $(\phi_\tau)_{|\Fil^0\Lambda_\tau}$ se réduit par $(\ker(\theta),p)$ sur, 
\[\widetilde{\Ha_\tau} : \bigwedge^{q_\tau} \omega_{G^D,\tau}\pmod p \fleche \bigwedge^{q_\tau} \omega_{G^D,\tau}^{(p^f)}\pmod p,\]
dont la valuation du déterminant est exactement $\Ha_\tau(G)$.
L'image de la réduction par $(\ker\theta,p)$ de $\pi_\tau(\left(\Fil^{\underline f}\Lambda\right)^{\Phi = D_{\underline f}})$ est incluse dans,
\[ \{ x \dans \bigwedge^{q_\tau} \omega_{G^D,\tau} \pmod p : \widetilde{\Ha}_\tau(x) = x \otimes 1\}.\]

\pro
\label{pro7}
Notons $N = \{ x \dans \bigwedge^{q_\tau} \omega_{G^D,\tau} \pmod p : \widetilde\Ha_\tau(x) = x \otimes 1\}$. C'est un sous-$\FF_{p^f}$-module du $\mathcal O_{C}/p$ module libre de rang 1, 
$\bigwedge^{q_\tau} \omega_{G^D,\tau}\pmod p$. Supposons que $\Ha_\tau(G) < 1 - \frac{1}{p^f}$.
Alors,
\[\im\left(N \fleche \bigwedge^{q_\tau} \omega_{G^D,\tau}\pmod {p^{1-\Ha_\tau(G)}}\right),\]
est un $\FF_{p^f}$-module de rang 1 engendré par un élément de valuation $\frac{\Ha_\tau(G)}{p^f-1}$.

Si $\Ha_\tau(G) \geq 1 - \frac{1}{p^f}$ cette image est toujours incluse dans l'image de 
$$p^{1/p^f}\bigwedge^{q_\tau}\omega_{G^D,\tau}.$$
\epro

\dem
On peut directement appliquer \cite{Far} proposition 7, lorsque $\Ha_\tau(G) < \frac{1}{2}$. On plonge $\FF_{p^f}$ dans $\mathcal O_C$ via le morphisme multiplicatif de Teichmuller. En général, quitte à choisir une base de $\bigwedge^{q_\tau} \omega_{G^D,\tau}$, on est ramené à résoudre une équation dans
$\mathcal O_C/p$ du type,
\[ X^{p^f} \equiv aX \pmod{p},\]
où $a= a_0p^{v(a)}$ avec $a_0 \dans \mathcal O_C^\times$ et $v(a) = \Ha_\tau(G)$. On écrit $x = up^w$, avec $u \dans \mathcal O_C^\times$ et $w \dans [0,1]$, et on en déduit,
\[ \min(1,wp) = \min(v(a) + w,1).\]
En analysant les quatre possibilités, on trouve que si $v(a) < 1 - \frac{1}{p^f}$, alors les solutions sont dans l'image de 
$p^{\frac{v(a)}{p^f-1}}a_o^{\frac{1}{p^f-1}}\mathbb F_{p^f} + p^{1-v(a)}\mathcal O_C$, où $a_o^{\frac{1}{p^f-1}}$ est un choix d'une racine de $a_0$ (on retrouve en particulier l'énonce de 
la proposition 7 de Fargues lorsque le module est de dimension 1) et si $v(a) \geq 1 - \frac{1}{p^f}$, elles sont dans $p^{\frac{1}{p^f}}\mathcal O_C$. 
\edem

On en déduit donc après identification $\bigwedge^{q_\tau} \omega_{G^D,\tau} = \mathcal O_C$ que l'image par $\theta$ de $\pi_\tau(\left(\Fil^{\underline f}\Lambda\right)^{\Phi = D_{\underline f}})$ est, à un inversible près, incluse dans \[p^{\frac{\Ha_\tau(G)}{p^f-1}}[\FF_{p^f}] + p^{1-\Ha_\tau(G)}\mathcal O_C.\]
Donc grâce au diagramme précédent, que l'image de $\bigwedge^{q_\tau}_\mathcal O \alpha_{G,\tau}$ est, à un inversible près, incluse dans,
\[ p^{K_\tau + \frac{\Ha_\tau(G)}{p^f-1}}[\FF_{p^f}]  + p^{K_\tau+1-\Ha_\tau(G)}\mathcal O_C.\]

\rem
Pour plus de simplicité, on va prendre $\frac{1}{2}$ comme borne pour $\Ha_\tau(G)$ au lieu de $1- \frac{1}{p^f}$, pour rendre les bornes plus cohérente avec celles de Fargues. De plus, on n'aurait probablement pas été capable de montrer que  nos sous-groupes "canoniques" sont des crans d'une certaine filtration de Harder-Narasihman sous l'hypothèse $\Ha_\tau(G) < 1 - \frac{1}{p^f}$.
Malheureusement, le passage à la puissance extérieure sur $\omega_{G^D,\tau}$ ne nous permettra pas en general de garder comme borne $\frac{1}{2}$, ce qui fait que la condition sur $\Ha_\tau(G)$ pour avoir un sous-groupe canonique sera malgré tout compliquée (mais ce sera simplement $\frac{1}{2}$ lorsque $p$ sera assez grand).
\erem

En particulier, on en déduit la proposition suivante,

\pro
Soit $G/\mathcal O_C$ un $\mathcal O$-module $p$-divisible de hauteur $h$ et signature $(p_\tau,q_\tau)_\tau$. Soit $\tau$ tel que $q_\tau\not\in\{0,h\}$.
Supposons de plus que $q_\tau < p-1$ et $\Ha_\tau(G) < \frac{1}{2}$.
Alors,
\[ \im\left(\bigwedge^{q_\tau}_\mathcal OT_pG \fleche \bigwedge^{q_\tau} \omega_{G^D,\tau}\pmod{1+K_\tau - \Ha_\tau(G)}\right),\]
est un $\FF_{p^f}$-module libre de rang 1, engendré par un élément de valuation $K_\tau + \frac{\Ha_\tau(G)}{p^f-1}$.
\epro

\section{Filtration canonique de la $p$-torsion}
\label{sect6}
Dans cette section, on fixe un plongement $\tau$. 
Soit $G$ un $\mathcal O$-module $p$-divisible tronqué d'échelon $r$ de hauteur $h$ et signature $(p_\tau,q_\tau)$.
Supposons que $r > k_\tau +1$ (pour pouvoir définir l'invariant de Hasse) et $q_\tau < p-1$.
Rappelons le théorème suivant,

\thr[Wedhorn \cite{Wed2} Theorem 2.8]
Si $p > 2$ et si $\mathcal D$ est une donnée PEL non ramifiée, et $\underline{X_0}$ un $BT$ avec $\mathcal D$-structure. Alors pour tout $1 \leq n \leq m \leq \infty$, 
le morphisme de foncteurs,
\[ \Def(\underline{X_0[p^m]}) \overset{[p^n]}{\fleche}  \Def(\underline{X_0[p^n]}),\]
est formellement lisse.
Les déformations sont à prendre avec la $\mathcal D$-structure.
\ethr

En particulier, si $G$ est un $\mathcal O$-module de Barsotti-Tate tronqué sur $\mathcal O_C$, alors il existe $\widehat{G}$ un $\mathcal O$-module $p$-divisible sur 
$\mathcal O_C$ tel que $G = \widehat{G}[p^r]$, et on peut donc appliquer les résultats de la section précédente à $G$ ! 

\subsection{Noyau de Hodge-Tate}

\pro
\label{prodegwedge}
Supposons que $\Ha_\tau(G) < \frac{1}{2}$. Alors pour tout $\eps$ tel que $K_\tau + \frac{\Ha_\tau(G)}{p^f-1} < \eps < 1 + K_\tau - \Ha_\tau(G)$,
\[\im\left(\bigwedge_\mathcal O^{q_\tau} \alpha_{G,\tau}\right)\pmod{p^\eps}\]
 est un $\FF_{p^f}$-module de rang 1, et
\[ \deg\Coker\left( \bigwedge^{q_\tau} \alpha_{G,\tau}\otimes 1\right) = K_\tau + \frac{\Ha_\tau(G)}{p^f-1}.\]
\epro

\dem
C'est essentiellement la proposition \ref{pro7} que l'on applique à $\hat G$, un $\mathcal O$-module $p$-divisible tel que $\hat G[p^r] = G$, puisque l'image est non triviale et incluse dans un $\FF_{p^f}$ module de rang 1. 
\edem
On en déduit en particulier que $\im(\alpha_{G,\tau,\frac{1+K_\tau-\Ha_\tau(G)}{q_\tau}})$ est engendré par au plus $q_\tau$ éléments sur $\mathcal O$.
Posons $\eps_\tau = \min(1,\frac{1+K_\tau-\Ha_\tau(G)}{q_\tau}) \leq 1$. La proposition centrale est alors la suivante,

\pro
Si $\Ha_\tau(G) < 1 + K_\tau - \frac{q_\tau}{p-1}$, alors $\frac{1}{p-1} < \eps_\tau$.
Supposons que \[\Ha_\tau(G) < \min(\frac{1}{2},  1 + K_\tau - \frac{q_\tau}{p-1}),\] alors pour tout $\frac{1}{p-1} < \eps < \eps_\tau$,
on a que,
\[ \dim_{\FF_{p^f}} \Ker \alpha_{G[p],\tau,\eps} = p_\tau.\]
\epro

\dem
En effet, on a par la proposition (\ref{proker1}) que,
\[ \dim_{\FF_{p^f}} \Ker \alpha_{G[p],\tau,\eps} \leq p_\tau.\]
Or la proposition précédente assure que $\im(\alpha_{G[p],\tau,\eps})$ est engendré par au plus $q_\tau$ éléments sur $\FF_{p^f}$ donc que,
\[ \dim_{\FF_{p^f}} \Ker \alpha_{G[p],\tau,\eps} = p_\tau.\qedhere\]
\edem

\rem\begin{enumerate}
\item Si $p$ est assez grand devant $q_\tau$, la seule hypothèse dans la proposition est $\Ha_\tau(G) < \frac{1}{2}$.
\item Si $K_\tau + \frac{\Ha_\tau(G)}{p^f-1} < \frac{1}{p-1}$, alors la proposition précédente s'applique encore avec $K_\tau + \frac{\Ha_\tau(G)}{p^f-1} < \eps < \eps_\tau$, 
en effet, $\im\alpha_{G,\tau,\eps_\tau} \fleche \im\alpha_{G,\tau,\eps}$ est injective pour de tels $\eps$.
\end{enumerate}
\erem

\subsection{Degrés}

Notons $K/\QQ_p$ une extension valuée complète quelconque, telle que $v(p) =1$. Rappelons alors les définitions de \cite{FarHN}.

\defi
Soit $M$ un $\mathcal O_K/p^r$-module de présentation finie, annulé par une puissance de $p$. 
On note $\delta = \Fitt_0M$,  c'est un idéal fractionnaire de $\mathcal O_K$
Alors on définit le degré de $M$ par,
\[\deg M = v(\delta).\]
\edefi

\exe
Si \[ M \simeq \prod_{i=1}^r \quotient{\mathcal O_K}{x_i\mathcal O_K},\]
alors $\deg M = \sum_{i=1}^r v(x_i).$
\eexe

\lem
\label{lemdeg}
Si $M = \Coker(f : L \fleche P)$ avec $L,P$ deux $\mathcal O_K/p^r$ modules libres, et $v(det(f)) < r$, 
alors $\deg M = \det f$.
\elem

\dem
En effet, choisissons un relèvement de $f$, $\widetilde f : \mathcal O_K^r \fleche \mathcal O_K^n$ tel que $M = \Coker \widetilde f$, alors $\det f \equiv \det \widetilde f \pmod{p^r}$.
\edem

\rem
Ce n'est plus vrai sans l'hypothèse $v(\det f) < r$ puisque $v(\det f) \dans [0,r]$, et on peut avoir $v(\det f) = r$ mais $\deg M > r$ (prendre $pI_n$, avec $n$ assez grand). 
\erem

On peut en particulier déduire du lemme précédent et de la proposition (\ref{prodegwedge}) la proposition suivante sur le degré du conoyau de $\alpha_{G,\tau}$.

\pro
\label{prodeg}
Soit $G$ est un $BT_r^{\mathcal O}$ avec $r > k_\tau$. 
Supposons que \[\Ha_\tau(G) < min(\frac{1}{2}, 1 + K_\tau - \frac{q_\tau}{p-1}),\] 
alors, pour tout $\min(K_\tau + \frac{\Ha_\tau}{p^f-1},\frac{1}{p-1}) < \eps < K_\tau + 1 - \Ha_\tau$,
\[\deg \Coker\left(\alpha_{G[p],\tau,\eps}\otimes 1\right) = K_\tau + \frac{\Ha_\tau}{p^f-1}.\]
\epro

\dem
D'après l'hypothèse sur $G$, on sait que, pour $K_\tau + \frac{\Ha_\tau}{p-1} < \eps \leq K_\tau + 1 - \Ha_\tau$, 
\[\deg \Coker\left(\bigwedge^{q_\tau}\alpha_{G,\tau}\otimes 1\right)_{K_\tau + 1 - \Ha_\tau} = K_\tau + \frac{\Ha_\tau}{p^f-1}.\]
De plus, comme on sait que le conoyau de $\alpha_{G,\tau}\otimes 1$ est tué par $p^{\frac{1}{p-1}}$, on a que 
$p^{\frac{1}{p-1}}\omega_{G^D,\tau} \subset \im\alpha_{G^D,\tau}\otimes 1$, et donc
\[ \deg\Coker\alpha_{G,\tau}\otimes 1 = \deg\Coker\alpha_{G,\tau,\eps}\otimes 1, \quad \forall \eps > \frac{1}{p-1}.\]
Donc $\deg\Coker(\alpha_{G[p],\tau,\eps}\otimes 1) = \deg\Coker(\alpha_{G,\tau,q_\tau\eps_\tau}\otimes1).$
Et donc, comme $K_\tau + \frac{\Ha_\tau}{p^f-1} < 1 + K_\tau -\Ha_\tau(G)$, le lemme précédent s'applique dans la 2e égalité,
\begin{IEEEeqnarray*}{cl} \deg\left(\Coker(\alpha_{G[p],\tau,\eps}\otimes 1)\right) 
&= \deg\left(\Coker(\alpha_{G,\tau,1 + K_\tau - \Ha_\tau}\otimes 1)\right) \\
&=  \deg\Coker\left(\bigwedge^{q_\tau}\alpha_{G}\otimes 1\right)_{1 + K_\tau - \Ha_\tau} \\
&=  K_\tau + \frac{\Ha_\tau}{p^f-1}. \qedhere\end{IEEEeqnarray*}
\edem

\subsection{Le théorème principal}

\thr
\label{thrptors}
Supposons donné une signature $(p_\tau,q_\tau)_\tau$. Choisissons un plongement $\tau$, tel que $q_{\tau} \not\in \{0,h = p_\tau + q_\tau\}.$ 
Supposons $p -1 > q_\tau$ (donc $p \neq 2$ en particulier).
Soit $G$ un groupe de Barsotti-Tate tronqué sur $\mathcal O_C$, de rang $k_\tau+1$, avec action de $\mathcal O$, et de signature $(p_\tau,q_\tau)$. 
Supposons de plus que \[\Ha_\tau(G) := \omega < \min(\frac{1}{2},1+K_\tau -\frac{q_\tau}{p-1}).\]
Posons $\eps_\tau = \min(1,\frac{K_\tau + 1 - \omega}{q_\tau}) \leq 1$.
Alors,
\[ \dim_{\FF_{p^f}} \Ker \alpha_{G[p],\tau,\eps_\tau} = p_\tau.\]
De plus sous l'hypothèse, 
\begin{equation}
\label{hypdeg}
\tag{H1}
\frac{2q_\tau}{p-1} < 1 + K_\tau, \quad \text{et} \quad  \Ha_\tau(G) < 1 + K_\tau - \frac{2q_\tau}{p-1},
\end{equation}
on a alors,
\begin{enumerate}
\item Soit $C$ l'adhérence dans $G[p]$ de $\Ker \alpha_{G[p],\tau,\eps_\tau} $. On a alors, si on note $E = G[p]/C$,
\[ p^{f-1}\deg_{\sigma\tau}(E) + p^{f-2}\deg_{\sigma^2\tau}(E) + \dots + \deg_{\tau}(E) = K_\tau (p^f-1)+ \Ha_\tau(G).\]
\item Le conoyau de $\alpha_{E,\tau}\otimes 1$ est de degré $K_\tau + \frac{\Ha_\tau(G)}{p^f-1}$.
\end{enumerate}
Remarquons que $\deg_\tau C_\tau^D = \deg_\tau E$ ! (Mais c'est faux pour les autres plongements).
\ethr

\rem
L'hypothèse (\ref{hypdeg}) n'est nécessaire que pour calculer la formule des degrés de $C$, et elle est bien sûr inutile si $p$ est assez grand, grâce à $\Ha_\tau(G) < \frac{1}{2}$.
Il serait intéressant de voir si l'on peut s'en passer.
De plus, on notera parfois $C_\tau$ au lieu de $C$ pour bien préciser à quel plongement est associé ce sous-groupe. 
\erem

\dem
Posons $E = G[p]/C$. Par définition de $C$, on a que $\mathfrak m_{C, 1- \eps_\tau}\im(\alpha_{C,\tau}\otimes 1) = 0$.
Donc $\mathfrak m_{C,1-\eps_\tau + \frac{1}{p-1}}\omega_{C^D,\tau} = 0$. On en déduit que pour, $\eps$ vérifiant,
\[ \min(K_\tau + \frac{\Ha_\tau}{p^f-1},\frac{1}{p-1}) < \eps \leq \eps_\tau - \frac{1}{p-1},\]
(ce qui est possible par l'hypothèse (\ref{hypdeg}) sur $p$), on a que,
\[\omega_{G^D,\tau,\eps} \simeq \omega_{E^D,\tau,\eps}.\]
Choisissons $e_1,\dots,e_{q}$ une $\FF_{p^f}$ base de $E(O_C)$, et $E_i$ l'adhérence schématique dans $E$ et $\FF_{p^f}(e_1,\dots,e_i)$.
On a alors, \[ 0 = E_0 \subset E_1 \subset E_2 \subset \dots \subset E_q = E,\]
une filtration dont les gradués sont des $p$-groupes de Raynaud munis d'une action de $\mathcal O$.
Regardons la filtration $\Fil_i\omega_{E^D,\tau} = \im(\omega_{E_i^D,\tau} \fleche \omega_{E^D,\tau}).$
On a une flèche naturelle $q_i : \omega_{(E_i/E_{i-1})^D,\tau} \fleche \Gr_i\omega_{E^D,\tau}$.
Comme l'image de $e_i$ dans $\Gr_i \omega_{E^D,\tau}$ est non nulle, donc $\alpha_{E_i/E_{i-1},\tau}(e_i) \neq 0$.
De plus par calcul direct sur les groupes de Raynaud, on a,
\[ \deg\Coker\alpha_{E_i/E_{i-1},\tau} = \frac{p^{f-1}\deg_\tau(E_i/E_{i-1})+p^{f-2}\deg_{\sigma\tau}(E_i/E_{i-1}) + \dots + \deg_{\sigma^{f-1}\tau}(E_i/E_{i-1})}{p^f-1}.\]
Et donc grâce à $q_i$ qui est surjective entre modules monogènes, 
\[ \deg\Coker\alpha_{E_i/E_{i-1},\tau}  = \deg(\Gr_i \omega_{E^D,\tau}/q_i\circ\alpha_{E_i/E_{i-1},\tau}(e_i)).\]
Mais d'après la proposition 8 de \cite{Far}, on peut écrire, 
\begin{IEEEeqnarray*}{rcl}\deg \Coker(\alpha_{E,\tau,\eps} \otimes 1) &=& \sum_{i=1}^q \deg(\Gr_i(\omega_{E^D,\tau}/ q_i\circ\alpha_{E_i/E_{i-1},\tau}(e_i))) \\
&=& \frac{1}{p^f-1}\sum_{i=1}^q p^{f-1}\deg_{\sigma\tau}(E_i/E_{i-1})+p^{f-2}\deg_{\sigma^2\tau}(E_i/E_{i-1}) + \dots + \deg_{\tau}(E_i/E_{i-1})\\
&=& \frac{p^{f-1}\deg_{\sigma\tau}(E) + p^{f-2}\deg_{\sigma^2\tau}(E) + \dots + \deg_{\tau}(E)}{p^f-1}
\end{IEEEeqnarray*}
Mais d'après la proposition \ref{prodeg}, on sait que \[\deg \Coker(\alpha_{E,\tau,\eps} \otimes 1) = \deg \Coker(\alpha_{G,\tau,\eps} \otimes 1) = K_\tau + \frac{\Ha_\tau}{p^f-1},\] 
d'où le résultat.
\edem

\rem
La formule correspond bien au calcul explicite sur le lieu $\mu$-ordinaire.
\erem

\rem
\label{remdeg}
De la formule sur les degrés de $E$, on en déduit, comme $\deg_{\tau'} E = \deg_{\tau'} G[p] - \deg_{\tau'} C = p_{\tau'} - \deg_{\tau'} C$, que,
\begin{IEEEeqnarray*}{ccc}
\sum_{i=1}^{f} p^{f-i}\deg_{\sigma^{i}\tau}(C) &=&  \sum_{i=1}^{f} p^{f-i}p_{\sigma^{i}\tau} - \sum_{i=1}^{f} \max(0,q_\tau - q_{\sigma^{-i}\tau})p^{i} - \Ha_\tau \\
 & = & \sum_{i=1}^{f} \min(p_\tau,p_{\sigma^i\tau})p^{f-i} - \Ha_\tau
\end{IEEEeqnarray*}
En particulier, $C$ est de grand degré !
De plus, $C$ est de $\mathcal O$-hauteur $p_\tau$, donc pour tout $i$, $\deg_{\sigma^i\tau} C \leq p_\tau$. Ensuite, comme $C \subset G[p]$ et $\deg_{\sigma^i\tau}G[p] = p_{\sigma^i\tau}$, on a,
\[\deg_{\sigma^i\tau} C \leq \min(p_\tau,p_{\sigma^i\tau}).\]
Donc si \[\sum_{i=1}^{f} p^{f-i}\deg_{\sigma^{i}\tau}(C) =\sum_{i=1}^{f} \min(p_\tau,p_{\sigma^i\tau})p^{f-i} - \Ha_\tau,\]
\begin{IEEEeqnarray*}{ccccc}
\sum_{i=1}^{f} \min(p_\tau,p_{\sigma^i\tau})p^{f-i} - \Ha_\tau &=& \deg C + \sum_{i=1}^{f} (p^{f-i}-1)\deg_{\sigma^{i}\tau}(C)\\
&\leq& \deg C + \sum_{i=1}^{f} (p^{f-i}-1)\min(p_{\sigma^i\tau},p_\tau)
 \end{IEEEeqnarray*}
Et donc, $\deg C \geq \sum_{i=1}^{f} \min(p_{\sigma^i\tau},p_\tau) - \Ha_\tau$.
Donc si $\Ha_\tau < \frac{1}{2}$, on a que $C$ est canonique, cf. \cite{Bij}, ou proposition (\ref{probij}).

On a aussi les inégalités suivantes, comme $\deg_{\sigma^i\tau} C \leq  \min(p_{\sigma^i\tau},p_\tau)$, alors $\deg_{\sigma^i\tau} C \geq \min(p_{\sigma^i\tau},p_\tau) - \frac{\Ha_\tau}{p^{f-i}}$, et donc
\[ \deg_\tau(C_\tau) \geq p_\tau - \Ha_\tau(G) \quad \text{et}\quad \deg_\tau(C_\tau^D) \leq \Ha_\tau.\]
\erem

\pro
Sous les hypothèses du théorème précédent, à savoir $\frac{2q_\tau}{p-1} < 1 + K_\tau$ et,
\[ \Ha_\tau(G) < \min(\frac{1}{2},1 + K_\tau - \frac{2q_\tau}{p-1}),\]
On a en fait que $C$ est l'adhérence schématique de $\Ker(\alpha_{G,\tau,1-\Ha_\tau(G)})$.
\epro

\dem
En effet, l'égalité sur la somme coefficientée des degrés partiels de $C_\tau$ nous dit en particulier, comme $\deg_{\sigma^i\tau}C \leq \min(p_\tau,p_{\sigma^i\tau})$, que 
\[\deg_\tau(C^D) \leq \Ha_\tau(G).\]
On en déduit que $\im(\omega_{C^D,\tau} \fleche \omega_{G[p]^D,1-\Ha_\tau(G)}) = 0$ et donc $C(\mathcal O_C) \subset \Ker(\alpha_{G,\tau,1-\Ha_\tau(G)})$.
\edem

\pro
Supposons que $\Ha_\tau(G) < \min(\frac{1}{p^f},1 + K_\tau - \frac{2q_\tau}{p-1})$. Alors la suite suivante est exacte :
\[ 0 \fleche C_\tau(K) \fleche G[p](K) \overset{\alpha_{G,\tau}}{\fleche} \omega_{G^D,\tau}.\]
De plus, le conoyau de la flèche,
\[ \alpha_{G,\tau} \otimes 1 : G[p](K) \otimes O_C \fleche \omega_{G^D,\tau},\]
est exactement $K_\tau + \frac{\Ha_\tau(G)}{p^f-1}$.
\epro

\dem
Si $E$ est un schéma en groupes de Raynaud avec une action de $\mathcal O$, tel que $\sum_{i=0}^{f-1} p^{f-i}\deg_{\sigma^i\tau}(E) \geq 1 - \frac{1}{p^f}$,
alors on peut vérifier que l'application $\alpha_{E,\tau}$ est nulle.

Maintenant soit $C_\tau$ le noyau de $\alpha_{G,\tau,\eps_\tau}$, donné par le théorème. Alors dans ce cas, filtrons $C_\tau$ par des sous $\mathcal O$-modules, tel que les graduées $(E_k)_{k = 1,\dots,p_\tau}$ soit des $\mathcal O$-modules de Raynaud.
De l'égalité
 \[\sum_{i=1}^{f} p^{f-i}\deg_{\sigma^{i}\tau}(C_\tau) =\sum_{i=1}^{f} \min(p_\tau,p_{\sigma^i\tau})p^{f-i} - \Ha_\tau,\]
on en déduit en particulier que $\sum_{k=1}^{p_\tau} \deg_\tau(E_k) = \deg_{\tau}(C_\tau) \geq p_\tau - \Ha_\tau(G)$.
Donc en particulier,
\[\deg_\tau(E_k) \geq 1 - \Ha_\tau(G),\]
et donc,
\[\sum_{i=0}^{f-1} p^{f-i}\deg_{\sigma^i\tau}(E_k) \geq 1 - \Ha_\tau(G) \geq 1 - \frac{1}{p^f}.\]
Donc $\alpha_{E_k,\tau} = 0$ et donc comme on peut toujours filtrer $C_\tau(\mathcal O_C)$ par des sous-$\mathcal O$-modules dont le premier contient un $x \dans C_\tau(O_C)$ donné, 
on en déduit que $\alpha_{C_\tau,\tau} = 0$.
Et on a montré précédemment que $\dim_{\mathbb F_{p^f}} \Ker(\alpha_{G,\tau}) \leq p_\tau$.
D'où l'affirmation sur la suite exacte.
Maintenant comme $\alpha_{G,\tau}$ se factorise par $G[p]/C_\tau$, et qu'on a calculé son conoyau de Hodge-Tate dans le théorème précédent, on en déduit la proposition.
\edem

\subsection{Compatibilités}

\pro
\label{produal}
Soit $G \dans \mathcal{BT}_{k_\tau+1}^\mathcal O(O_C)$. Supposons que $\Ha(G) = \Ha(G^D)$ soit strictement inférieur à $\min(\frac{1}{2}, 1 + K_\tau - \frac{2q_\tau}{p-1})$ (donc $G$ et $G^D$ ont tous les deux un sous-groupe canonique).
Notons $C_\tau$ le sous groupe canonique de $G$ et $D_\tau$ celui de $G^D$.
Alors $C_\tau = D_\tau^\perp$.
\epro

\dem
En effet, \begin{IEEEeqnarray*}{ccccc}
\deg D_\tau^\perp &=& \Ht(G^D[p]/D_\tau) - \deg(G^D[p]/D_\tau) &=& fp_\tau - \deg G^D[p] + \deg D_\tau \\
&\geq&  fp_\tau - \sum_i q_{\sigma^i\tau} + \sum_i \min(q_\tau,q_{\sigma^i\tau}) - \Ha(G^D)\\
 &\geq& \left(\sum_{i=0}^{f-1} p_\tau - q_{\sigma^i\tau} + \min(q_\tau,q_{\sigma^i\tau})\right) - \Ha(G^D)\\
 & \geq &  \left(\sum_{i=0}^{f-1} p_\tau - q_{\sigma^i\tau} + h - \max(p_\tau,p_{\sigma^i\tau})\right) - \Ha(G^D)\\
 & \geq &  \left(\sum_{i=0}^{f-1} p_\tau + p_{\sigma^i\tau} - \max(p_\tau,p_{\sigma^i\tau})\right) - \Ha(G^D)\\
 & \geq &  \left(\sum_{i=0}^{f-1} \min(p_\tau,p_{\sigma^i\tau})\right) - \Ha(G^D)\\
 \end{IEEEeqnarray*}

 Donc $D_\tau^\perp$ est de grand degré, et par unicité, cf. \cite{Bij}, et annexe (\ref{probij}), $C_\tau = D_\tau^\perp$.
 On pourrait aussi utiliser la proposition \ref{corHNp} de la section suivante.
\edem

\rem
Si on ne connaissait pas la compatibilité de $\Ha_\tau$ à la dualité, en supposant que $\Ha_\tau(G)$ et $\Ha_\tau(G^D)$ sont assez petit (plus petits que
$\min(\frac{1}{2}, 1 + K_\tau - \frac{2q_\tau}{p-1})$) on retrouve qu'ils sont égaux, en comparant les degrés des sous-groupes de $G$ et $G^D$ et en utilisant le résultat d'unicité de \cite{Bij}.
\erem

\rem
Dans la construction de Fargues, \cite{Far}, il est prouvé que le sous-groupe canonique de $G$ déforme le noyau de Frobénius, modulo $p^{1-\Ha(G)}$. Un calcul sur les modules de Dieudonné dans le cas $\mu$-ordinaire laisse entendre que le sous groupe canonique associé à un plongement $\tau$ devrait correspondre à,
\[(p^{r_\tau -1} \Ker F^f)[p],\]
la $p$-torsion de l'image par la multiplication par $p^{r_\tau - 1}$, $r_\tau = |\{ \tau' : q_{\tau'} < q_\tau\}|$. Mais il ne semble pas évident qu'un tel sous-groupe existe en 
général (i.e. qu'il soit representable ou fini et plat).
Néanmoins on montrera à la fin de l'article un résultat partiel, qui décrit une déformation (modulo $\mathfrak m_C$...) de $\Ker F^f$. 
En particulier, les sous-groupes $C_\tau$ seront inclus dans $\Ker F^f$.
\erem

\pro[Compatibilité entre les différents plongements]
Si $C_\tau$ désigne le sous-groupe canonique associé à l'abscisse $q_\tau$, alors $q_\tau \leq q_{\tau'} \Rightarrow C_{\tau'} \subset C_\tau$.
En particulier si $\tau, \tau'$ sont associés à la même abscisse de rupture $q_\tau = q_{\tau'}$, alors $C_\tau = C_{\tau'}$.
\epro

\dem
Cela va découler de la section suivante (Corollaire \ref{corHNp}), car on va prouver que $C_\tau$ est une abscisse de rupture de Harder-Narasimhan. On pourrait aussi bien utiliser \cite{Bij}, rappelé ici 
en annexe, proposition (\ref{probij}).
\edem

\rem
En particulier dans le cas des $\mathcal O$-modules stricts, disons associés à $\tau_0$, il n'y a qu'un plongement intéressant ($\tau_0$) et donc qu'un sous-groupe canonique, et dans ce cas $k_{\tau_0}=K_{\tau_0} = 0$, donc tout est plus simple. De plus, la dualité stricte de Faltings permet de simplifier beaucoup de choses : en particulier on peut montrer que le sous-groupe canonique ainsi construit relève le noyau de $F^f$.
\erem

\section{Calculs de polygones de Harder-Narasihman}
\label{sect7}
\subsection{Polygone de Harder-Narasihman classique}

\pro
Soit $G$ un schéma en groupes avec $\mathcal O$-action. Alors $\HN_{\mathcal O}(G)(x) = \frac{1}{f}\HN(G)(fx)$ a des abscisses de rupture entières.
\epro

\dem
Les abscisses de ruptures de $\HN(G)$ sont les hauteurs des groupes apparaissant dans la filtration HN de $G$, or celle-ci est stable par $\mathcal O$, donc 
les hauteurs de ces groupes sont des multiples de $f$.
\edem

\rem
Dans \cite{FarHNpdiv} et \cite{Shen}, il est introduit des polygones de HN "renormalisés" pour des $\mathcal{BT}_n$ par,
\[ \widetilde{\HN}(G[p^n])(x) = \frac{1}{n}\HN(G[p^n])(nx).\]
Ces polygones ont des abscisses de rupture dans $\frac{1}{n}\ZZ$, elles ne sont plus nécessairement entières !
\erem

\begin{figure}[h]
\begin{center}
\caption{$\mathcal O$-polygone de Hodge renversé associée à la signature $(q_\tau)_{\tau\in \mathcal I}$.}
\label{figHN}
\begin{tikzpicture}[line cap=round,line join=round,>=triangle 45,x=0.5cm,y=0.5cm]
\draw[->,color=black] (-0.5,0.) -- (20.,0.);
\foreach \x in {,1.,2.,3.,4.,5.,6.,7.,8.,9.,10.,11.,12.,13.,14.,15.,16.,17.,18.,19.}
\draw[shift={(\x,0)},color=black] (0pt,2pt) -- (0pt,-2pt);
\draw[->,color=black] (0.,-0.5) -- (0.,10.);
\foreach \y in {,1.,2.,3.,4.,5.,6.,7.,8.,9.}
\draw[shift={(0,\y)},color=black] (2pt,0pt) -- (-2pt,0pt);
\clip(-0.5,-1.) rectangle (20.,10.);
\draw (0.,0.)-- (3.,3.);
\draw (3.,3.)-- (4.99016393443,4.50789096126);
\draw (4.99016393443,4.50789096126)-- (6.9868852459,5.32855093257);
\draw (10.,6.)-- (14.0098360656,6.63845050215);
\draw (14.004934692,6.63767010007)-- (17.0049180328,6.65423242468);
\draw [dash pattern=on 1pt off 1pt] (3.,3.)-- (3.,0.);
\draw [dash pattern=on 1pt off 1pt] (4.99016393443,4.50789096126)-- (5.,0.);
\draw [dash pattern=on 1pt off 1pt] (6.9868852459,5.32855093257)-- (7.,0.);
\draw [dash pattern=on 1pt off 1pt] (10.,6.)-- (10.,0.);
\draw [dash pattern=on 1pt off 1pt] (14.004934692,6.63767010007)-- (14.,0.);
\draw [dash pattern=on 1pt off 1pt] (17.004656654,6.65423098166)-- (17.,0.);
\draw [dash pattern=on 1pt off 1pt] (8.10573770492,5.54949784792)-- (8.8631147541,5.7230989957);
\draw (2.8,-0.1) node[anchor=north west] {$p_r$};
\draw (4,-0.1) node[anchor=north west] {$p_{r-1}$};
\draw (6,-0.1) node[anchor=north west] {$p_{r-2}$};
\draw (13.5,-0.1) node[anchor=north west] {$p_1$};
\draw (9.5,-0.1) node[anchor=north west] {$p_2$};
\draw (0.996721311475,2.5) node[anchor=north west] {1};
\draw (15.3352459016,7.6) node[anchor=north west] {0};
\draw (11.5655737705,7.6) node[anchor=north west] {$\frac{n_1}{f}$};
\draw (1.8,5.5) node[anchor=north west] {$\frac{n_1+\dots+n_{r-1}}{f}$};
\draw (4.3,6.5) node[anchor=north west] {$\frac{n_1+\dots+n_{r-2}}{f}$};
\draw (16.6,-0.100430416069) node[anchor=north west] {$h$};
\draw (0.101639344262,-0.116212338594) node[anchor=north west] {$0$};
\end{tikzpicture}
\end{center}
\end{figure}
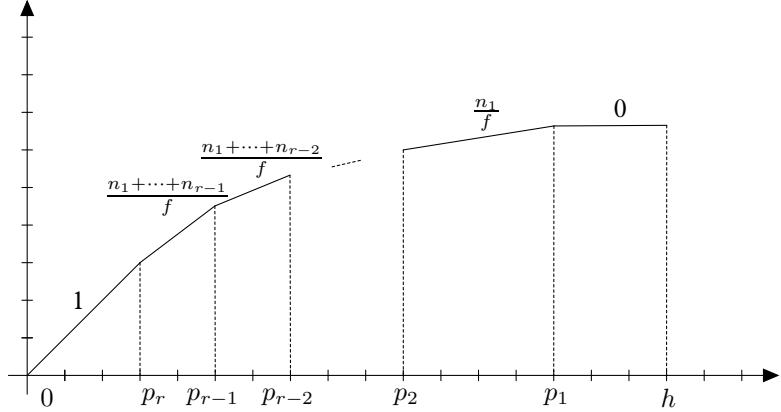

\pro
On peut tracer le $\mathcal O$-polygone de Hodge renversé (voir figure \ref{figHN}) d'un $\mathcal{BT}^\mathcal O$ de signature $(p_\tau,q_\tau)$, $\Hdg^\diamond$. 
Il a pour abscisses de rupture 
\[ 0 \leq p^r < p^{r-1} < \dots < p^1 \leq h,\]
et où les pentes sont données par $(1,\frac{ |\{\tau : p_\tau \geq p^{r-1}\}|}{f},\frac{ |\{\tau : p_\tau \geq p^{r-2}\}|}{f},\dots,\frac{ |\{\tau : p_\tau \geq p^{1}\}|}{f},0)$.
On vérifie facilement que c'est aussi $\widetilde{\HN}_{\mathcal O}(G^{\mu-ord}[p^n])$, le $\mathcal O$-polygone de Harder-Narasimhan (renormalisé) de la 
$p^n$-torsion du groupe $\mu$-ordinaire associé,
$\forall n \dans \NN$.
\epro

\rem
On retrouve un cas (très) particulier de \cite{Shen}, qui prédit que lorsque le polygone de Hodge et de Newton se touche en une abscisse de rupture du polygone de Newton (comme dans le cas $\mu$-ordinaire) alors ces polygones touchent aussi le polygone de Harder-Narasimhan, qui a lui aussi une rupture en cette abscisse.
\erem

\rem
On a bien sur l'égalité des ensembles (et donc de leurs cardinaux),
\[ \{\tau : p_\tau \geq p^{i}\} = \{\tau : q_\tau \leq q^{i}\}.\]
\erem

\pro
\label{prokern}
Si  $\eps < 1 - \frac{1}{p-1}$, le noyau de $\alpha_{G,\tau,n-\eps}$ est engendré sur $\mathcal O$ par moins de $p_\tau$ éléments, i.e.
\[ \Ker \alpha_{G,\tau,n-\eps} = \bigoplus_{i=1}^{p_\tau} \quotient{\mathcal O}{p^{a_i}\mathcal O},\quad 0 \leq a_i \leq n.\]
\epro

\dem
C'est exactement \cite{Far} proposition 13, en remplaçant $\ZZ_p$ par $\mathcal O$ (la proposition 12 se généralise trivialement).
\edem

\pro
\label{pro14}
Soit $G$ un $\mathcal{BT}_n^\mathcal O$ de signature $(p_\tau,q_\tau)_\tau$. Supposons qu'il existe un sous-$\mathcal O$-module $C$ tel que,
\[ \Ht_{\mathcal O}(C) = np_\tau \quad \text{et} \quad \deg_\tau(G/C) < 1 - \frac{1}{p-1}.\]
Alors, si $\eps = \deg_\tau(G/C), C(\mathcal O_{\overline K}) = \Ker \alpha_{G,\tau, n-\eps}$ qui est un $\mathcal O/p^n\mathcal O$-module libre.
\epro

\dem
Si $\eps = \deg_\tau(G/C) = \deg_\tau(G) - \deg_\tau(C) = np_\tau - \deg_\tau(C) = \deg(\omega_{C^D,\tau})$ donc 
$\omega_{C^D,\tau}\fleche \omega_{G^D,\tau} \fleche \omega_{G^D,\tau,\eps}$ est nulle, et par conséquent $C(\mathcal O_K) \subset \Ker(\alpha_{G,\tau,n-\eps})$,
 la proposition (\ref{prokern}) conclut.
\edem

\pro
\label{propoly0}
Soit $\tau \dans \mathcal I$.
Soit $C \subset G[p^n]$ un sous-$\mathcal O$-module de hauteur $fnp_\tau$ (i.e. $\rg(C) = p^{nfp_\tau}$). 
Supposons que 
\[\deg C > n\sum_{j = 0}^{f-1} \min(p_\tau,p_{\sigma^j\tau}) - \frac{|\{\tau' : q_\tau = q_\tau'\}|}{2},\]
Alors le polygone de Harder-Narasimhan $\widetilde{\HN}_{\mathcal O}(G[p^n])$ a un point de rupture en l'abscisse $p_\tau$.
\epro

\dem
Voir la démonstration de la proposition \ref{propoly1} qui s'applique aussi ici.
\edem

On déduit en particulier de cette dernière proposition,

\cor
\label{corHNp}
Soit $G$ un $\mathcal{BT}_{k_\tau +1}^\mathcal O$ tel que $\Ha_\tau(G) < \min(\frac{1}{2}, 1 + K_\tau - \frac{2q_\tau}{p-1}).$
Soit $C_\tau$ le sous-groupe canonique donné par le théorème (\ref{thrptors}).
Alors $C_\tau$ est un cran de la filtration de Harder-Narasihman de $G[p]$. 
En particulier on retrouve que $C_\tau$ est compatible à la dualité au sens de la proposition (\ref{produal}), mais aussi que, si
le théorème (\ref{thrptors}) s'applique pour $G$ pour deux plongements $\tau$ et $\tau'$, alors,
\[ q_{\tau'} \leq q_\tau \Rightarrow C_\tau \subset C_{\tau'}.\]
\ecor

\dem
On a $\deg C_\tau \geq \sum_{i=1}^{f} \min(p_{\sigma^i\tau},p_\tau) - \Ha_\tau(G)$ et par hypothèse, $\Ha_\tau(G) < \frac{1}{2}$, la proposition précédente s'applique. Il reste à montrer que le sous-groupe $C'$ qui induit la rupture – i.e. le cran de la filtration de HN à l'abscisse $p_\tau$ – est $C_\tau$. Or $\deg_{\tau'}(C') \leq \min(p_\tau,p_{\tau'})$ car il est de hauteur $p_\tau$ et que $C' \subset G$ qui est de $\tau'$-degré $p_{\tau'}$.
Donc si $\deg_\tau(C') < p_\tau - \frac{1}{2}$, alors $\deg(C') = \sum_{\tau'} \deg_{\tau'}(C') < \sum_{\tau'} \min(p_\tau,p_{\tau'} - \frac{1}{2} < \deg(C_\tau)$, ce qui est absurde puisque $C'$ est un cran HN à l'abscisse $p_\tau$. Donc $\deg_\tau(C') > p_\tau - \frac{1}{p-1}$ et la proposition \ref{pro14} assure que $C'$ est l'adhérence schématique de 
$\Ker \alpha_{G,\tau,1-\eps}$, or c'est aussi le cas de $C_\tau$, donc $C_\tau = C'$.
\edem

\rem
Malheureusement pour la $p^n$-torsion avec $n >1$, on n'arrivera pas à montrer la généralisation du corollaire précédent (parce que les bornes sur $\Ha_\tau$ seront moins grossières, on pourrait néanmois y arriver quitte à sacrifier les bornes), mais on va devoir changer un peu la filtration de 
Harder-Narasihman en fonction du plongement $\tau$, c'est l'objet de la sous-section suivante. Cela peut s'expliquer en partie par le fait que le théorème \ref{thrptors} ne nous 
donne pas les degrés des sous-groupes canoniques, mais seulement des combinaisons linéaires des degrés partiels.
\erem

\subsection{Fonction degré et polygones de Harder-Narasihman modifiés}

Notons $\mathfrak{Gr}_p^{\mathcal O}(\mathcal O_K)$ la catégorie (exacte) des schémas en groupes finis et plats sur $O_K$ (d'ordre une puissance de $p$) où 
$K$ est une extension valuée de $\ZZ_p$.
On suppose qu'il existe un plongement $K \supset F$.

\defi
Pour tout $\tau \in \mathcal I$, et tout $G/\mathcal O_K$ un schéma en groupes avec action de $\mathcal O$, on définie une nouvelle fonction degrée $\Deg_\tau$ par,
\[ \Deg_\tau(G) = \sum_{j=1}^f p^{f-j}\deg_{\sigma^i\tau}(G).\]
Cette fonction degré vérifie les propriétés
\begin{enumerate}
\item $\Deg_\tau$ est additive sur les suites exactes dans $\mathfrak{Gr}_p^{\mathcal O}(\mathcal O_K)$.
\item Si $u : G \fleche G'$ est un morphisme qui devient un isomorphisme en fibre générique, alors $\Deg_\tau(G') \geq \Deg_\tau(G)$, avec égalité si et seulement si
$u$ est un isomorphisme.
\end{enumerate}
\edefi

\dem
Voir \cite{BPS} Proposition 1.19. 
\edem

Ces propriétés sont analogues à celles vérifiées par la fonction degré de \cite{FarHN}, et elles permettent de développer un formalisme Harder-Narasihman.

On note la fonction de pente,
\[ \mu_\tau = \frac{\Deg_\tau}{f\Ht_{\mathcal O}},\]
elle est à valeurs dans $[0,\frac{p^f-1}{f(p-1)}]$. On aurait pu la renormaliser pour la rendre à valeurs dans $[0,1]$, mais cela aurait (inutilement) alourdi les formules qui suivront.

À partir de maintenant, dans cette sous-section, fixons un $\tau \dans \mathcal I$.

On a alors la proposition, voir \cite{FarHN}, Théorème 1, ou \cite{And},

\pro
Soit $G$ un groupe fini plat sur $O_K$ (d'ordre une puissance de $p$),  muni d'une action de $\mathcal O$, il possède une unique filtration par des sous-groupes finis et plats,
\[0 = G_0 \subsetneq G_1 \subsetneq G_2 \subsetneq \dots \subsetneq G_r = G\]
telle que,
\begin{enumerate}
\item Pour tout $i$, $G_{i+1}/G_{i}$ est semi-stable pour la fonction de pente $\mu_\tau$.
\item Pour tout $i \geq 1$, $\mu_\tau(G_i/G_{i-1}) > \mu_\tau(G_{i+1}/G_i)$.
\end{enumerate}
On fera référence à cette filtration comme $\tau$-filtration de Harder-Narasihman. On notera $\HN_\tau(G)$ le polygone concave de Harder-Narasihman, 
définie par les pentes $(\mu_\tau(G_i/G_{i-1}))_{i = 1,\dots,r}$ avec multiplicités $(\Ht_\mathcal O(G_i/G_{i-1}))_{i = 1,\dots,r}$. C'est un polygone à abscisses de ruptures entières et à 
pentes rationnelles.
\epro

\rem
On peut vérifier que si $G$ est $\mu$-ordinaire, alors sa filtration de Harder-Narasihman "classique" (i.e. définie par la fonction $\deg$, cf. \cite{FarHN}) vérifie les deux propriétés 
de la proposition ci-dessus, en particulier sur le lieu $\mu$-ordinaire, les filtrations de Harder-Narasihman données par $\mu$ ou $\mu_\tau$ sont égales, pour tout $\tau \dans \mathcal I$ (mais les polygones sont différents).
Dans le cas général, les filtrations sont différentes, mais on peut montrer que si $^\mu\Ha(G)$ est suffisamment petit (mais a priori sans borne précise, sauf pour un sous-groupe bien précis de ces filtrations), alors elles coïncident.
\erem

\exe
Supposons que $G$ soit un $\mathcal O$-module $p$-divisible $\mu$-ordinaire de signature $(p_\tau,q_\tau)$, alors on calcule explicitement,
\[ \HN_\tau(G[p])(p_{\tau'}) = \frac{1}{f}\sum_{i=1}^f p^{f-i}\min(p_{\tau'},p_{\sigma^i\tau}).\]
\eexe

\rem
La différence entre deux pentes consécutives du polygone $\HN_\tau$ $\mu$-ordinaire est donc au moins $\frac{1}{f}$, cela servira dans la démonstration du théorème \ref{thrntors}, 
pour montrer que le sous-groupe canonique est un cran de la $\tau$-filtration Harder-Narasihman.
\erem

On va utiliser ces filtrations "modifiées" pour mettre en famille la filtration canonique. En effet, les sous-groupes de la filtration canonique seront des sous-groupes apparaissant dans 
les $\tau$-filtrations de Harder-Narasihman, pour différents plongements $\tau$. A priori il n'est pas clair qu'ils apparaissent aussi dans la filtration de Harder-Narasihman classique, 
d'où la nécessité d'introduire ces nouvelles filtrations.

\defi
Si $G$ est un $\mathcal{BT}_n$ avec action de $\mathcal O$, on note $\widetilde{\HN}_\tau(G)$ la renormalisation (en fonction de $n$) de $\HN_\tau(G)$, c'est à dire,
\[\widetilde{\HN}_\tau(G)(x) = \frac{1}{n}\HN_\tau(G)(nx),\]
est donc un polygone à abscisses entre $0$ et $\Ht_\mathcal O(G[p])$.
\edefi

On va avoir besoin d'utiliser l'analogue de la proposition \ref{propoly0} pour les nouvelles filtrations :

\pro
\label{propoly1}
Soit $\tau' \dans \mathcal I$. Soit $G$ un $\mathcal{BT}_n$.
Soit $C \subset G[p^n]$ un sous-$\mathcal O$-module de hauteur $fnp_{\tau'}$ (i.e. $\rg(C) = p^{nfp_{\tau'}}$). 
Supposons que 
\[\Deg_\tau C > n\sum_{j = 0}^{f-1} p^{f-j}\min(p_{\tau'},p_{\sigma^j\tau}) - \frac{\sum_{j=1}^f p^{f-j}\delta_{q_{\sigma^j\tau} = q_{\tau'}}}{2},\]
Alors le polygone de Harder-Narasimhan $\widetilde{\HN}_{\tau}(G[p^n])$ a un point de rupture en l'abscisse $p_{\tau'}$.
\epro

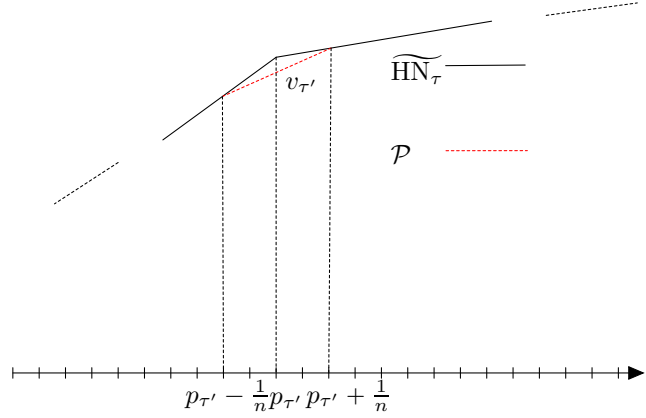
\begin{figure}[h]
\caption{Rupture autour de $\widetilde{\HN}_\tau(G^{\mu-ord}[p^n])(p_{\tau'})$.}
\label{fig2}
\begin{center}
\begin{tikzpicture}[line cap=round,line join=round,>=triangle 45,x=0.7cm,y=0.7cm]
\draw[->,color=black] (3.,0.) -- (15.,0.);
\foreach \x in {3.,3.5,4.,4.5,5.,5.5,6.,6.5,7.,7.5,8.,8.5,9.,9.5,10.,10.5,11.,11.5,12.,12.5,13.,13.5,14.,14.5}
\draw[shift={(\x,0)},color=black] (0pt,2pt) -- (0pt,-2pt);
\clip(3.,-1.) rectangle (15.,8.);
\draw (8.,6.)-- (5.84885245902,4.42898134864);
\draw (8.,6.)-- (12.0923770492,6.68579626973);
\draw [dash pattern=on 1pt off 1pt,color=ffqqqq] (6.9803922743,5.25536360309)-- (9.04412536448,6.17497343756);
\draw [dash pattern=on 1pt off 1pt] (5.,4.)-- (3.79073770492,3.2137733142);
\draw [dash pattern=on 1pt off 1pt] (13.1300819672,6.7962697274)-- (14.8595901639,7.03299856528);
\draw (11.2196218159,5.84072805214)-- (12.7305234552,5.85220581398);
\draw [dash pattern=on 1pt off 1pt,color=ffqqqq] (11.2264121951,4.22405901773)-- (12.7719039984,4.22549373796);
\draw (7.7,-0.2) node[anchor=north west] {$p_{\tau'}$};
\draw (8,5.8) node[anchor=north west] {$v_{\tau'}$};
\draw (10,6.24338106169) node[anchor=north west] {$\widetilde{\HN_\tau}$};
\draw (10,4.51383905308) node[anchor=north west] {$\mathcal P$};
\draw [dash pattern=on 1pt off 1pt] (9.04412536448,6.17497343756)-- (9.,0.);
\draw [dash pattern=on 1pt off 1pt] (8.,6.)-- (8.,0.);
\draw [dash pattern=on 1pt off 1pt] (6.9803922743,5.25536360309)-- (7.,0.);
\draw (6.1,0.04) node[anchor=north west] {$p_{\tau'} - \frac{1}{n}$};
\draw (8.4,0.04) node[anchor=north west] {$p_{\tau'} + \frac{1}{n}$};
\end{tikzpicture}
\end{center}
\end{figure}

\dem
Si $\widetilde{\HN}_{\tau}(G[p^n])$ (renormalisé donc) n'a pas de rupture en $p_{\tau'}$, alors il est en dessous du polygone $\mathcal P$, voir figure \ref{fig2}.

(polygone $\widetilde{\HN}(G^{\mu-ord})$ avec une droite reliant $p_\tau - \frac{1}{n}$ et $p_\tau +\frac{1}{n}$).

On peut explicitement calculer 
\[\mathcal P(p_{\tau'}) = v_{\tau'} =\frac{\sum_{j = 0}^{f-1} p^{f-j}\min(p_{\tau'},p_{\sigma^j\tau})}{f} - \frac{\sum_{j=1}^f p^{f-j}\delta_{q_{\sigma^j\tau} = q_{\tau'}}}{2nf}.\]
En particulier, s'il existe un sous-groupe comme dans l'énoncé, et qu'il n'y a pas de rupture, on contredit la proposition 13 de \cite{FarHN}, qui s'applique encore dans ce cadre 
(c'est-à-dire pour $\HN_\tau$).
\edem

\pro
\label{propoly}
Supposons maintenant que le polygone de Harder-Narasihman $\widetilde{\HN}_\tau(G[p^n])$ a une rupture en $p_{\tau}$, et seulement l'hypothèse,
\begin{equation}
\label{hyp2}
\tag{H2}
\Deg_{\tau} C > n\sum_{j = 0}^{f-1} p^{f-j}\min(p_{\tau},p_{\sigma^j\tau}) - \frac{p-2}{p-1}.
\end{equation}
Alors en plus $C$ est un cran de la $\tau$-filtration de Harder-Narasimhan de $G[p^n]$.
\epro

\rem
L'hypothèse (\ref{hyp2}) permettra de traiter le cas où $\tau = \tau'$ et $n_{\tau} = |\{\theta: q_{\tau} = q_\theta\}| = 1$ lorsque l'hypothèse de la proposition précédente n'est pas vérifiée.
\erem

\dem
Supposons (\ref{hyp2}),
\[ \deg_{\tau} C > np_{\tau} - \frac{p-2}{p-1} \quad \text{i.e.}\quad \deg_{\tau}(G/C) < \frac{p-2}{p-1} = 1 - \frac{1}{p-1},\]
donc d'après le point (1) de la proposition \ref{pro14}, $C =  \Ker(\alpha_{G,\tau,n-\eps})$, où $\eps = \deg_{\tau}(G/C) < \frac{p-2}{p-1}$.
Soit $C'$ le cran de la filtration de Harder-Narasimhan de $G[p^n]$ d'abscisse $np_{\tau}$, qui existe puisqu'on a supposé qu'il y avait une rupture. 
Dans ce cas, on a $\Deg_\tau(C') \geq \Deg_\tau(C) > n\sum_{j = 0}^{f-1} p^{f-j}\min(p_{\tau},p_{\sigma^j\tau}) - \frac{p-2}{p-1}$, donc en particulier,
\[  \deg_{\tau} C' > np_{\tau} - \frac{p-2}{p-1} \quad \text{i.e.}\quad \deg_{\tau}(G/C') < \frac{p-2}{p-1} = 1 - \frac{1}{p-1}.\]
D'après le point (1) de la proposition \ref{pro14}, on a donc $C' = \Ker(\alpha_{G,\tau,n-\eps'})$, 
où $\eps' = \deg_{\tau}(G/C')$. Si $\eps' \leq \eps$, alors $C \subset C'$, mais $\Ht C = \Ht C'$, donc $C = C'$. Idem si $\eps' \geq \eps$.
\edem

\section{Filtration canonique supérieure}
\label{sect8}
\subsection{Récurrence et théorème principal}

\pro
\label{prorecdeg}
Soit $r = k_\tau +1$ et $G \dans \mathcal{BT}_{r+1}^\mathcal O$ tel que $\Ha_\tau(G) < \frac{1}{p^f + 1}$.
Soit $C_\tau$ son $\tau-$sous-groupe canonique. Alors $\quotient{p^{-r}C_\tau}{C_\tau}$ est un $\mathcal{BT}_r^\mathcal O$, et,
\[\Ha_\tau(\quotient{p^{-r}C_\tau}{C_\tau}) \leq p^f \Ha_\tau(G).\]

Plus précisément,
\[\Ha_\tau(\quotient{p^{-r}C_\tau}{C_\tau}) = (p^f-1)\deg_\tau(C_\tau^D) + \Ha_\tau(G).\]
\epro

\dem
Considérons la flèche sur $\mathcal O_K$,
\[ G \overset{\pi}{\fleche} \quotient{G}{C_\tau}.\]
Elle induit une suite exacte,
\[\omega_{G^D,\tau} \overset{\pi_*}{\fleche} \omega_{(G/C_\tau)^D,\tau} \fleche \omega_{C_\tau^D,\tau} \fleche 0.\]
En effet, regardons nos groupes sur $\mathcal O_C/p$ et considérons le triangle distingué (cf. \cite{Ill}, VII.3.1.1.5),
\begin{center}
\begin{tikzpicture}[description/.style={fill=white,inner sep=2pt}] 
\matrix (m) [matrix of math nodes, row sep=3em, column sep=2.5em, text height=1.5ex, text depth=0.25ex] at (0,0)
{ 
\ell_{C_\tau^D} & & \ell_{G^D} \\
 & \ell_{(G/C_\tau)^D} & \\
 };

\path[->,font=\scriptsize] 
(m-1-1) edge node[auto] {$$} (m-1-3)
(m-1-3) edge node[auto] {$$} (m-2-2)
(m-2-2) edge node[auto,left] {$+1$} (m-1-1)
;
\end{tikzpicture}
\end{center}
$G$ est un $\mathcal{BT}_r$, donc son complexe de co-lie $\ell_{G^D}$ est $\omega_{G^D} \oplus \omega_{G^D}[1]$. Bien sur, $G/C_\tau$ n'est pas un $\mathcal{BT}_r$ à priori, mais on peut écrire la suite exacte,
\[0 \fleche \quotient{p^{-1}C_\tau}{C_\tau} \fleche \quotient{G}{C_\tau} \overset{p}{\fleche} \quotient{pG}{C_\tau}\fleche 0,\]
qui induit un triangle distingué,
\begin{center}
\begin{tikzpicture}[description/.style={fill=white,inner sep=2pt}] 
\matrix (m) [matrix of math nodes, row sep=3em, column sep=2.5em, text height=1.5ex, text depth=0.25ex] at (0,0)
{ 
\ell_{(p^{-1}C_\tau/C_\tau)^D} & & \ell_{(G/C_\tau)^D} \\
 & \ell_{(pG/C_\tau)^D} & \\
 };

\path[->,font=\scriptsize] 
(m-1-1) edge node[auto] {$$} (m-1-3)
(m-1-3) edge node[auto] {$p$} (m-2-2)
(m-2-2) edge node[auto,left] {$+1$} (m-1-1)
;
\end{tikzpicture}
\end{center}
Or, sur $\mathcal O_C/p$, la flèche $p$ est nulle, on en déduit un isomorphisme, $\ell_{(p^{-1}C_\tau/C_\tau)^D} \simeq \ell_{(G/C_\tau)^D}$.
Or, $\quotient{p^{-1}C_\tau}{C_\tau}$ est un $\mathcal{BT}_1$, on en déduit, $\ell_{(p^{-1}C_\tau/C_\tau)^D} \simeq 
\omega_{(p^{-1}C_\tau/C_\tau)^D}\oplus\omega_{(p^{-1}C_\tau/C_\tau)^D}[1] \simeq \omega_{(G/C_\tau)^D}\oplus\omega_{(G/C_\tau)^D}[1] $ sur $\mathcal O_C/p$, 
et donc, la suite exacte longue du premier triangle donne,
\begin{eqnarray*}
 0 \fleche \mathcal{H}^{-1}(\ell_{C_\tau^D}) \fleche  \mathcal{H}^{-1}(\ell_{G^D}) \fleche  \mathcal{H}^{-1}(\ell_{(G/C_\tau)^D}) \fleche  \mathcal{H}^{0}(\ell_{C_\tau^D}) \\
\fleche  \mathcal{H}^{0}(\ell_{G^D}) \fleche  \mathcal{H}^{0}(\ell_{(G/C_\tau)^D}) \fleche 0.\end{eqnarray*}
Il suffit alors de montrer que $ \mathcal{H}^{0}(\ell_{C_\tau^D}) \simeq \omega_{C_\tau^D}
\fleche  \mathcal{H}^{0}(\ell_{G^D}) \simeq \omega_{G^D}$ est nulle.
Or cette flèche est la réduction modulo $p$ de celle sur $\mathcal O_C$, qui est donnée par,
\[ \omega_{C_\tau^D} \fleche \omega_{G^D} \simeq (\mathcal O_C/p^r)^d,\]
et $C_\tau$ est tué par $p$, donc $\omega_{C_\tau^D}$ est de $p$-torsion, la flèche précédente est donc de la forme $p^{r-1}\phi$, et comme $r \geq 2$, elle est nulle une fois réduite modulo $p$. On a donc la suite exacte; qui est de plus $\mathcal O$-équivariante puisque c'est le cas de tous les morphismes entre les schémas en groupes,
\[ \omega_{G^D} \fleche \omega_{(G/C_\tau)^D} \fleche \omega_{C_\tau^D} \fleche 0.\]

On en déduit donc que, d'après l'inégalité \ref{remdeg} et le lemme \ref{lemdeg},
\[ \deg_\tau(C_\tau^D) = \deg(\omega_{C_\tau^D,\tau}) = \det(\pi_*) \leq \Ha_\tau(G).\]
Notons de plus que, \[k_\tau(G) = k_\tau(\quotient{p^{-r}C_\tau}{C_\tau}),\]
En effet, soit $G \dans \mathcal{BT}^\mathcal O(\mathcal O_C)$ tel que $\Ha_\tau(G) < \frac{1}{p^f + 1}$. Soit $C_\tau/\mathcal O_C$ le sous-groupe du théorème.
L'isogenie,
\[ 0 \fleche C_\tau \fleche G \fleche G/{C_\tau}\fleche 0,\]
induit,
\[\omega_{G^D,\tau'} \fleche \omega_{G/C_\tau^D,\tau'} \fleche \omega_{C_\tau^D,\tau'} \fleche 0,\]
Or $\omega_{G^D,\tau'}$ et $\omega_{G/C_\tau^D,\tau'}$ sont sans $p$-torsion (car $BT$) et $\omega_{C_\tau^D,\tau'}$ est de $p$-torsion, donc après inversion de $p$ on voit que 
$G$ et $G/C_\tau$ ont même signature.
Donc les signatures de $G$, $G/C_\tau$ sont les mêmes ! En particulier ils ont les mêmes $k_\tau$ pour tout $\tau$.

Mais le carré suivant au niveau des cristaux est commutatif,
\begin{center}
\begin{tikzpicture}[description/.style={fill=white,inner sep=2pt}] 
\matrix (m) [matrix of math nodes, row sep=3em, column sep=2.5em, text height=1.5ex, text depth=0.25ex] at (0,0)
{ 
\mathbb D(G)_\tau& &\mathbb D(G)_\tau^{(p^f)} \\
\mathbb D(G/C_\tau)_\tau & &\mathbb D(G/C_\tau)^{(p^f)}_\tau \\
 };

\path[->,font=\scriptsize] 
(m-1-1) edge node[auto] {$V^f$} (m-1-3)
(m-1-1) edge node[auto] {$\pi_*$} (m-2-1)
(m-1-3) edge node[auto] {$\pi_*^{(p^f)}$} (m-2-3)
(m-2-1) edge node[auto] {$V^{'f}$} (m-2-3);
\end{tikzpicture}
\end{center}
Et donc comme $k_\tau(G) = k_\tau(G/C_\tau)$, on en déduit le même diagramme sur les puissances extérieures $q_\tau$ et avec les divisions de $V^f$, et donc la commutativité du diagramme modulo $p$ suivant,
\begin{center}
\begin{tikzpicture}[description/.style={fill=white,inner sep=2pt}] 
\matrix (m) [matrix of math nodes, row sep=3em, column sep=2.5em, text height=1.5ex, text depth=0.25ex] at (0,0)
{ 
\det\omega_{G^D,\tau} & &\det\omega_{G^D,\tau}^{(p^f)} \\
\det\omega_{(G/C_\tau)^D,\tau} & &\det\omega_{(G/C_\tau)^D,\tau}^{(p^f)} \\
 };

\path[->,font=\scriptsize] 
(m-1-1) edge node[auto] {$\widetilde\Ha_\tau(G)$} (m-1-3)
(m-1-1) edge node[auto] {$\pi_*$} (m-2-1)
(m-1-3) edge node[auto] {$\pi_*^{(p^f)}$} (m-2-3)
(m-2-1) edge node[auto] {$\widetilde\Ha_\tau(G/C_\tau)$} (m-2-3);
\end{tikzpicture}
\end{center}
On a donc l'égalité modulo $p$,
\[ \pi_*^{(p^f)}\circ \widetilde\Ha_\tau(G) \equiv \widetilde\Ha_\tau(G/C_\tau) \circ \pi_{*},\]
et donc en passant au déterminant,
\[ p^{p^f\deg_\tau(C_\tau^D) + \Ha_\tau(G)} \equiv p^{\Ha_\tau(G/C_\tau) + \deg_\tau(C_\tau^D)} \pmod p.\]
Or par hypothèse, 
\[ p^f\deg_\tau(C_\tau^D) + \Ha_\tau(G) \leq (p^f+1)\Ha_\tau(G) < 1,\]
et donc,
\[(p^f-1)\deg_\tau(C_\tau^D) + \Ha_\tau(G) = \Ha_\tau(G/C_\tau).\qedhere\]
\edem

\rem
On a aussi que si $p^f \deg_\tau(C_\tau^D) + \Ha_\tau(G) \geq 1$, alors \[\Ha_\tau(G/C_\tau) \geq 1 - \deg_\tau(C_\tau^D) \geq 1-\Ha_\tau(G).\]
\erem

\thr
\label{thrntors}
Soit $p > q_\tau + 1$. Soit $K/\mathcal O[1/p]$ une extension valuée. Soit $G$ un $\mathcal{BT}_{n+k_\tau}^\mathcal O(\mathcal O_K)$ de signature $(p_\tau,q_\tau)_\tau$.
Supposons que,
\[\Ha_\tau(G) < \frac{1}{p^{(n-1)f}}\min(\frac{1}{2},1+K_\tau - \frac{q_\tau}{p-1}).\]
Alors il existe $C_\tau^n \subset G[p^n]$ un sous $\mathcal O/p^n$-module de $G$.
Supposons de plus que \begin{equation}
\label{hypdeg2}
\tag{H3}
\frac{2q_\tau}{p-1} < 1 + K_\tau \quad \text{et} \quad \Ha_\tau(G) < \frac{1+K_\tau}{p^{(n-1)f}} - \frac{2q_\tau}{p^{nf} - p^{(n-1)f}}.
\end{equation}
Alors
\begin{enumerate}
\item $C_\tau^n(\mathcal O_{\overline K})$ est un $\mathcal O/p^n\mathcal O$-module libre.
\item $C_\tau^n(\mathcal O_{\overline K})$ coïncide avec le noyau de l'application $\alpha_{G,\tau, n-\frac{p^{nf} - 1}{p^f-1}\Ha_\tau(G)}$.
\item On a que,
\begin{eqnarray*} \sum_{i=1}^f \deg_{\sigma^i\tau}(G[p^n]/C_\tau^n)p^{f-i} &=& nK_\tau(p^f-1) + n\Ha_\tau(G) + (p^f-1)\left( \deg_\tau(C_\tau^{1,D}) + \dots + \deg_\tau(C_\tau^{n-1,D})\right) \\
& \leq& nK_\tau(p^f-1) + \frac{p^{nf} - 1}{p^f-1}\Ha_\tau(G).\end{eqnarray*}
Ou encore,
\[\Deg_\tau(C^n_\tau) = \sum_{i=1}^f \deg_{\sigma^i\tau}(C_\tau^n)p^{f-i} \geq n\sum_{i=1}^{f} \min(p_\tau,p_{\sigma^i\tau})p^{f-i} - \frac{p^{nf} - 1}{p^f-1}\Ha_\tau(G).\]
Et donc en particulier,
\[ \deg C_\tau^n = \sum_{i=1}^f \deg_{\sigma^i\tau}(C_\tau^n) \geq n\sum_{\tau'} \min(p_\tau,p_{\tau'}) - \frac{p^{nf} - 1}{p^f-1}\Ha_\tau(G).\]
\end{enumerate}
Notons que,
\[ \Ha_\tau(G/C_\tau^n) = \Ha_\tau(G) + (p^f-1)\deg_\tau(C_\tau^{n,D}) \leq p^{nf}\Ha_\tau(G).\]
On a de plus les propriétés suivantes,
\begin{enumerate}[(a)]
\item $C_\tau^n$ est un cran de la $\tau$-filtration de Harder-Narasimhan de $G[p^n]$.
\item $C_\tau^n$ est compatible à la dualité ; si $D_\tau^n$ est le sous-groupe canonique de $G^D[p^n]$, alors $D_\tau^n = C_\tau^{n,\perp}$.
\item $C_\tau^n[p^k] = C_\tau^k, \forall k \leq n$.
\item $C_\tau^n/C_\tau^k$ est le sous-groupe canonique de rang $n-k$ de $p^{n-k}C_k/C_k$.
\end{enumerate}
\ethr

\rem
On aimerait que l'hypothèse (\ref{hypdeg2}) ne soit pas nécessaire, remarquons que c'est le cas si le nombre premier $p$ est assez grand (avec même borne que pour l'hypothèse (\ref{hypdeg})). Si $p$ est assez grand, les hypothèses du théorème deviennent simplement $\Ha_\tau(G) < \frac{1}{2p^{(n-1)f}}.$

Si de plus on suppose que $\Ha_\tau(G) < \frac{1}{p^{(n-1)f}}\min(\frac{3}{8},1+K_\tau - \frac{2q_\tau}{p-1})$, la démonstration montre (facilement) que $C_\tau^n$ est un cran de la filtration de
Harder-Narasihman "classique" de $G[p^n]$.

$G^D$ vérifie aussi les hypothèses du théorème puisque $\Ha_\tau(G) = \Ha_\tau(G^D)$.
\erem

\dem
On construit $C_\tau^n$ par récurrence comme dans \cite{Far}. Avec la formule sur le degré pour $C_\tau^1$, et la proposition (\ref{prorecdeg}) toujours par récurrence, on trouve la formule sur le degré de $C_\tau^n$, et donc l'assertion de liberté et le point (2) grâce à la proposition \ref{pro14}.

Pour le cran de la $\tau$-filtration de Harder-Narasimhan, on procède comme dans \cite{Far}, et on sait d'après la proposition (\ref{propoly}) que si,
\[ \frac{p^{nf} - 1}{p^f-1}\Ha_\tau(G) <\frac{\sum_{j=1}^f p^{f-j}\delta_{q_{\sigma^j\tau} = q_{\tau'}}}{2},\]
on a une rupture. En particulier, si $n_\tau \geq 2$, on a la rupture, puisque $\frac{p^{nf} - 1}{p^f-1}\Ha_\tau(G) < 1 - \frac{1}{p-1}$. De plus par cette dernière inégalité, s'il y a rupture, le groupe correspondant est exactement $C_\tau^n$, car la deuxième partie de la proposition (\ref{propoly}) s'applique.

Il suffit seulement de montrer que si $\mathcal P = \HN_\tau(G[p^n])$ (non renormalisé, donc à abscisses de ruptures entières, et associé à la fonction $\Deg_\tau$), 
et si $i < np_\tau < j$ désignent des abscisses de ruptures de $\mathcal P$, alors,
\begin{equation}
\label{rupture}
 \frac{np_\tau - i}{j-i}\mathcal P(j) + \frac{j - np_\tau}{j-i}\mathcal P(i) < \frac{n}{f}\sum_{j=1}^f p^{f-j}\min(p_\tau,p_{\sigma^j\tau}) - \frac{1}{f}\frac{p^{nf} - 1}{p^f - 1} \Ha_\tau(G).\end{equation}
Comme dans \cite{Far}, en raisonnant sur un dessin, on voit que si $i < np_\tau - 1$,
\[\frac{np_\tau - i}{j-i}\mathcal P(j) + \frac{j - np_\tau}{j-i}\mathcal P(i) < \frac{n}{f}\sum_{\tau'} \min(p_\tau,p_{\tau'}) - \frac{2n_\tau}{3f}\]
Et comme,
\[ \frac{p^{nf}-1}{f(p^f-1)}\frac{1}{2p^{(n-1)f}} < \frac{2}{3f},\]
on en déduit que si $i < np_\tau - 1$ (ou si $j > np_\tau + 1$) il y a une rupture en l'abscisse $np_\tau$.

\rem
Ce raisonnement, pour ces $i,j$, marchait en fait aussi avec le polygone de Harder-Narasihman $\HN_{\mathcal O}$ "classique", c'est-à-dire associé avec la fonction $\deg$ et non 
$\Deg_\tau$.
\erem

Il reste donc le cas $n_\tau = 1$, $i = np_\tau - 1$ et $j = np_\tau +1$.
On suppose donc qu'il existe des ruptures de $\HN_\tau(G[p^n])$ aux abscisses $np_\tau-1$ et $np_\tau +1$, qui correspondent donc à des groupes $D$ et $D'$.

Essayons maintenant de montrer (\ref{rupture}) avec $i = np_\tau - 1$ et $j = np_\tau +1$, c'est-à-dire que $\Ht_\mathcal O(D) = np_\tau -1$ et $\Ht_\mathcal O(D') = np_\tau +1$.
Si le polygone de Harder-Narasihman $\HN_\tau(G[p^n])$ a une rupture en $np_\tau$, c'est gagné (proposition \ref{propoly}), supposons donc par l'absurde que ce n'est pas le cas.

\begin{figure}[h]
\begin{center}
\caption{$\tau$-degré $\Deg_\tau$ minimal de $D$.}
\label{HNbreak}
\begin{tikzpicture}[line cap=round,line join=round,>=triangle 45,x=1.5cm,y=0.7cm]
\draw[->,color=black] (1.93321476589,0.) -- (8.,0.);
\foreach \x in {1.5,2.,2.5,3.,3.5,4.,4.5,5.,5.5,6.,6.5,7.,7.5}
\draw[shift={(\x,0)},color=black] (0pt,2pt) -- (0pt,-2pt);
\clip(1.93321476589,-2.48344923044) rectangle (8.,9.44281939932);
\draw (5.5,7.)-- (7.,8.);
\draw (5.5,7.)-- (4.5,4.);
\draw [dash pattern=on 4pt off 4pt,domain=4:5.996845605649142] plot(\x,{(--8.80592903973-2.0658071613*\x)/-0.488648884338});
\draw [dash pattern=on 2pt off 2pt](5.5,7.)-- (5.5,0.);
\draw [dash pattern=on 2pt off 2pt](5.,0.)-- (4.99160688666,5.47482065997);
\draw [dash pattern=on 2pt off 2pt](6.,0.)-- (5.99684560565,7.33123040377);
\draw (2.95081967213,8.48493543759)-- (3.46721311475,8.48493543759);
\draw [dash pattern=on 4pt off 4pt] (2.96721311475,7.72740315638)-- (3.46721311475,7.74461979914);
\draw (3.58692551168,8.91238277446) node[anchor=north west] {$\HN_\tau^{\mu-ord}$};
\draw [dash pattern=on 2pt off 2pt] (5.50819672131,5.26542324247)-- (4.00819672131,5.26542324247);
\draw (3.7,7.5) node[anchor=north west] {$\HN_\tau^{\mu-ord}(np_\tau)$};
\draw (2.,5.8) node[anchor=north west] {$\HN_\tau^{\mu-ord}(np_\tau) - \delta$};
\draw [dash pattern=on 2pt off 2pt] (4.99337463583,3.16)-- (4.6,3.16);
\draw (1.87,3.7) node[anchor=north west] {$\HN_\tau^{\mu-ord}(np_\tau\!\!-\!\!1)\!\!-\!\!2\delta\!\!+\!\!\mu_1\!\!-\!\!\mu_2$};
\draw (5.39541706619,-0.14) node[anchor=north west] {$np_\tau$};
\draw (5.79458862552,-0.0879289891655) node[anchor=north west] {$np_\tau\!\!+\!\!1$};
\draw (4.7,-0.0879289891655) node[anchor=north west] {$np_\tau\!\!-\!\!1$};
\begin{scriptsize}
\draw [fill=qqqqff] (5.50819672131,5.26542324247) circle (1.5pt);
\draw[color=qqqqff] (5.8,5.1) node {$C_\tau^{min}$};
\end{scriptsize}
\end{tikzpicture}
\end{center}
\end{figure}
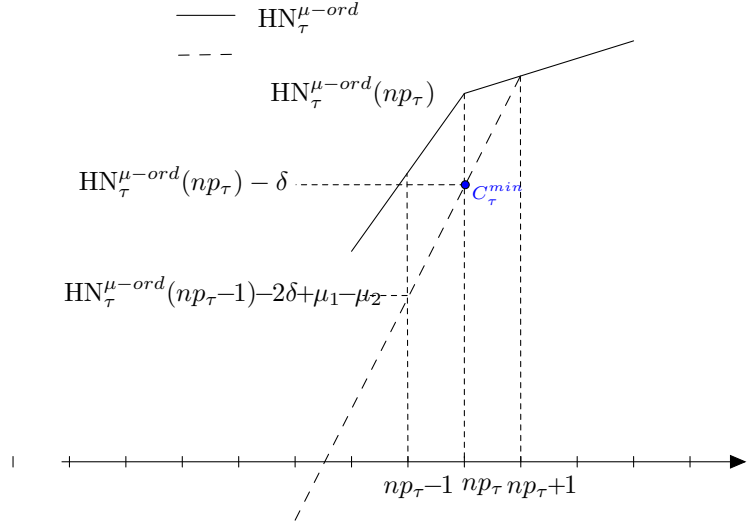

On a que,
\[ D(\mathcal O_C) \subset (\mathcal O/p^n\mathcal O)^{h}.\]
On écrit la suite, exacte en fibre générique,
\[ 0 \fleche D[p^{n-1}] \fleche D \fleche p^{n-1}D\fleche 0,\]
où $p^{n-1}D$ désigne l'adhérence schématique de $p^{n-1}D(\mathcal O_C)$, c'est un $\mathcal O$-module de hauteur $x \leq p_\tau -1$, puisque 
$\Ht_\mathcal O(D) = np_\tau -1$ et $D[p^{n-1}]$ est de hauteur supérieure ou égale à $(n-1)p_\tau$, on a donc,
\[ D(\mathcal O_C) = (\mathcal O/p^n\mathcal O)^x \oplus N,\]
où $N$ est un $\mathcal O/p^{n-1}\mathcal O$-module (de type fini). Moralement, plus $x$ est petit, plus le degré de $D$ aussi 
(puisqu'il est de plus en plus inclus dans la $p^{n-1}$-torsion de $G$). 
On va montrer qu'il est maximal, i.e. $x = p_\tau -1$.
Tout d'abord, essayons de minorer $\deg D$. Comme on a supposé que le polygone $\HN_\tau(G[p^n])$ n'avait pas de rupture en $np_\tau$, on peut donc en déduire par ce 
qui précède sur le degré du sous-groupe $C^n_\tau$ et la figure \ref{HNbreak} suivante, où la ligne pointillée représente la polygone minimal qui passe par le degré minimum
possible de $\Deg_\tau(C^n_\tau)$ autorisé par le théorème, et de telle manière qu'il n'y ait pas de rupture, et où $\mu_1,\mu_2$ sont les pentes du polygone $\mu$-ordinaire 
autour de $np_\tau$.
On en déduit donc, en ayant noté la borne donnée par le degré dans le théorème,
\[\delta = \frac{p^{nf}-1}{p^f-1}\Ha_\tau(G),\]
que,
\[\Deg_\tau D > \sum_{i=1}^f p^{f-i}\min(np_\tau-1,np_{\sigma^i\tau}) - 2\frac{p^{nf}-1}{p^f-1}\Ha_\tau(G) + \frac{1}{f},\]
et donc
\begin{equation}
\label{degD}
\deg D > \sum_{i=1}^f \min(np_\tau-1,np_{\sigma^i\tau}) - 2\frac{p^{nf}-1}{p^f-1}\Ha_\tau(G) + \frac{1}{f}.\end{equation}

Maintenant, on sait aussi que,
\[ \deg D \leq \deg D[p^{n-1}] + \deg p^{n-1}D,\]
et supposons que $x \leq p_\tau -2$.
On peut alors majorer (la seconde inégalité étant du à : si $p_{\sigma^j\tau} \leq x < p_\tau$ alors $(n-1)p_{\sigma^j\tau} < (n-1)p_\tau < np_\tau -1 - x$),
\begin{align*} \deg D &\leq \sum_{j=1}^f \min(np_\tau -1 -x,(n-1)p_{\sigma^j\tau}) + \min(x,p_{\sigma^j\tau})  \\
& \leq  \sum_{j=1}^f \min(np_\tau -1,np_{\sigma^i\tau},(n-1)p_{\sigma^j\tau} + x) \\
& \leq  np_\tau - 2 + \sum_{j=1}^{f-1}  \min(np_\tau -1,np_{\sigma^i\tau}) \\
& \leq \sum_{j=1}^{f}  \min(np_\tau -1,np_{\sigma^i\tau}) - 1
\end{align*}
ce qui contredit l'inégalité (\ref{degD}), d'après l'annexe (\ref{calculann1}).
On a donc que $x = p_\tau -1$, et on peut donc écrire,
\[ D(\mathcal O_C) = (\mathcal O/p^n\mathcal O)^{p_\tau-1} \oplus N,\]
où $N$ est un $\mathcal O$-module tué par $p^{n-1}$, de $\mathcal O$-longueur $n-1$. 

\lem
On a l'égalité,
\[ D[p^{n-1}]= C^{n-1}_\tau.\]
\elem

\dem
On a que ,
\[\deg D \leq \deg D[p^{n-1}] + \deg(p^{n-1}D).\]
Or ce dernier module est de $p$-torsion et hauteur $p_\tau -1$, on a donc que,
\[ \deg(p^{n-1}D) \leq \sum_i \min(p_\tau-1,p_{\sigma^i\tau}),\]
et donc,
\[ \deg D[p^{n-1}] \geq \deg D - \sum_i \min(p_\tau-1,p_{\sigma^i\tau}).\]
On en déduit en particulier en reprenant la minoration de $\deg D$ (\ref{degD}), que,
\begin{align*}
\deg D[p^{n-1}] &> \sum_{i=1}^f \min(np_\tau-1,np_{\sigma^i\tau}) - 2\frac{p^{nf}-1}{p^f-1}\Ha_\tau(G) + \frac{1}{f} - \sum_i \min(p_\tau-1,p_{\sigma^i\tau}) \\
		& =  \sum_{i=1}^f \min((n-1)p_\tau,(n-1)p_{\sigma^i\tau}) - 2\frac{p^{nf}-1}{p^f-1}\Ha_\tau(G) + \frac{1}{f}.
		\end{align*}
On peut donc utiliser la proposition (\ref{probij}) avec $C^{n-1}_\tau$ et $D[p^{n-1}]$, puisque,
\[ \deg C^{n-1}_\tau >  \sum_{i=1}^f \min((n-1)p_\tau,(n-1)p_{\sigma^i\tau}) - \frac{p^{(n-1)f}-1}{p^f-1}\Ha_\tau(G),\]
et on vérifie (voir annexe (\ref{calculann1})) que,
\[ \frac{p^{(n-1)f}-1}{p^f-1}\frac{1}{2p^{n-1}f} + 2 \frac{p^{nf}-1}{p^f-1}\frac{1}{2p^{n-1}f} - \frac{1}{f} \leq 1.\qedhere\]
\edem

Comme $p^{(n-1)}D$ est un $\mathcal O$-module de hauteur $p_\tau -1$,
\[ \Deg_\tau(p^{n-1}D) \leq \sum_{j=1}^f p^{f-j}\min(p_\tau-1,p_{\sigma^j\tau}).\]

On en déduit (utiliser la croissance par déformation de $\Deg_\tau$),
\begin{eqnarray*}\Deg_\tau(D) \leq \Deg_\tau(D[p^{n-1}]) + \Deg_\tau(p^{n-1}D) =\\ 
-(n-1)\Ha_\tau(G) -  (p^f-1)\left( \deg_\tau(C_\tau^{1,D}) + \dots + \deg_\tau(C_\tau^{n-2,D})\right) \\ + \sum_{j=1}^f p^{f-i}[(n-1)\min(p_\tau,p_{\sigma^i\tau}) + \min(p_\tau - 1,p_{\sigma^i\tau})] \\
= -(n-1)\Ha_\tau(G) -  (p^f-1)\left( \deg_\tau(C_\tau^{1,D}) + \dots + \deg_\tau(C_\tau^{n-2,D})\right) \\ + \sum_{j=1}^f p^{f-i}[n\min(p_\tau,p_{\sigma^i\tau}) - \delta_{p_{\sigma^i\tau} \geq p_\tau}] \\
\end{eqnarray*}

En appliquant cela à $G^D$ et $D^{'\perp}$, on trouve que,
\begin{eqnarray*}\Deg_\tau(D') \leq \Deg_\tau(G[p^n]) - \sum_{j} p^{f-j}(nq_\tau -1) + \Deg_\tau(D^{'\perp}) \\ 
\leq -(n-1)\Ha_\tau(G^D) -  (p^f-1)\left( \deg_\tau(C_\tau^{1,\perp,D}) + \dots + \deg_\tau(C_\tau^{n-2,\perp,D})\right) \\ + \sum_{j=1}^f p^{f-i}[np_{\sigma^i\tau} - nq_\tau + 1 + (n-1)\min(q_\tau,q_{\sigma^i\tau}) + \min(q_\tau - 1,q_{\sigma^i\tau})] \\
=  -(n-1)\Ha_\tau(G^D) -  (p^f-1)\left( \deg_\tau(C_\tau^{1,\perp,D}) + \dots + \deg_\tau(C_\tau^{n-2,\perp,D})\right) \\ + \sum_{j=1}^f p^{f-i}[np_{\sigma^i\tau} + np_\tau - (n-1)\max(p_\tau,p_{\sigma^i\tau}) + \max(p_\tau,p_{\sigma^i\tau}) + 1 - \delta_{q_{\sigma^j\tau} \geq q_\tau} ] \\
=  -(n-1)\Ha_\tau(G^D) -  (p^f-1)\left( \deg_\tau(C_\tau^{1,\perp,D}) + \dots + \deg_\tau(C_\tau^{n-2,\perp,D})\right) \\ + \sum_{j=1}^f p^{f-i}[n\min(p_\tau,p_{\sigma^i\tau}) + 1 
- \delta_{q_{\sigma^j\tau} \geq q_\tau}]
\end{eqnarray*}
Et on en déduit que,
\begin{eqnarray*} \Deg_\tau(D) + \Deg_\tau(D') \leq \sum_{j=1}^f p^{f-j}\left(2n\min(p_\tau,p_{\sigma^i\tau}) + 1 - \delta_{q_{\sigma^j\tau} \geq q_\tau} - \delta_{p_{\sigma^j\tau} \geq p_\tau}\right) \\
- 2(n-1)\Ha_\tau(G) -(p^f-1)\left(\deg_\tau(C_\tau^{1,D}) +\deg_\tau(C_\tau^{1,\perp,D})  + \dots + \deg_\tau(C_\tau^{n-2,D}) +\deg_\tau(C_\tau^{n-2,\perp,D}) \right).\end{eqnarray*}
Remarquons que $\deg_\tau(C_\tau^{j,\perp,D}) = \deg_\tau(C_\tau^{j,D})$.

Comme \[\Deg_\tau(C^n_\tau) =  n\sum_{i=1}^{f} \min(p_\tau,p_{\sigma^i\tau})p^{f-i} - n\Ha_\tau(G) - (p^f-1)\left( \deg_\tau(C_\tau^{1,D}) + \dots + \deg_\tau(C_\tau^{n-1,D})\right),\]
 pour montrer (\ref{rupture}) il suffit de voir que,
\[ \frac{1}{2}\sum_{j=1}^f p^{f-j}(\delta_{q_{\sigma^j\tau}\geq q_\tau} +\delta_{q_{\sigma^j\tau}\leq q_\tau}-1) > \Ha_\tau + (p^f-1)\deg_\tau C^{n-1,D}_\tau,\]
C'est-à-dire comme $n_\tau =1$ (le pire cas dans la formule précédente),
\[ \frac{1}{2} > \Ha_\tau + (p^f-1)\deg_\tau C^{n-1,D}_\tau = \Ha_\tau(G/C^{n-1}_\tau).\]
Mais comme on a supposé que $\Ha_\tau(G/C^{n-1}_\tau) \leq p^{(n-1)f}\Ha_\tau(G) < \frac{1}{2}$ pour pouvoir faire la récurrence, c'est gagné, il y a donc une rupture de $\HN_\tau(G[p^n])$ en l'abscisse $np_\tau$.
\edem

\rem
Les sous-groupes de la filtration canonique n'existent à priori qu'au dessus de $\mathcal O[1/p]$, mais le fait qu'ils correspondent à une rupture Harder-Narasihman 
va permettre de les redescendre.
\erem

On en déduit le théorème final,

\thr
\label{thrfinO}
Soit $G$ un $\mathcal O$-module $p$-divisible tronqué d'échelon $n+k$ sur $\Spec(O_K)$ de signature $(p_\tau,q_\tau)_\tau$, où $k = \max_\tau k_\tau$. 
Supposons que \[p > \max\{ \frac{2q_\tau}{1+K_\tau} : \tau \dans \mathcal I, q_\tau \neq h\} + 1.\]

Supposons de plus que,
\begin{equation}
\label{hypfinale}
\tag{$H_n$}
^\mu\Ha(G) < \frac{1}{p^{(n-1)f}}\min(\frac{1}{2},1+K_\tau - \frac{2q_\tau}{p-1}), \forall \tau \in \mathcal I \text{ tels que } q_\tau \notin \{0,h\}.\end{equation}
Alors il existe une (unique) filtration, appelée filtration canonique, $(\Fil_\tau G[p^n])_{\tau \in \mathcal I}$ de $G[p^n]$ par des sous-$\mathcal O/p^n$-modules finis et plats de 
$G[p^n]$, dont les inclusions sont données par,
\[ \Fil_\tau G[p^n] \subset \Fil_{\tau'} G[p^n] \quad \text{si et seulement si} \quad p_\tau \leq p_{\tau'},\]
telle que pour tout $\tau$,
\[\Ht_{\mathcal O} \Fil_\tau G[p^n] = np_\tau = nh-nq_\tau,\]
et pour tout $\tau$,
\[ \Deg_\tau(\Fil_\tau(G[p^n])) = \sum_{i=1}^f p^{f-i}\deg_{\sigma^i\tau}(\Fil_\tau G[p^n]) \geq n\sum_{i=1}^{f} \min(p_\tau,p_{\sigma^i\tau})p^{f-i} - \frac{p^{nf} - 1}{p^f-1}\Ha_\tau(G).\]
Et donc en particulier,
\[ \deg \Fil_\tau G[p^n] = \sum_{i=1}^f \deg_{\sigma^i\tau}(\Fil_\tau G[p^n]) \geq n\sum_{\tau'} \min(p_\tau,p_{\tau'}) - \frac{p^{nf} - 1}{p^f-1}\Ha_\tau(G).\]
De plus $\Fil_\tau G[p^n]$ coïncide avec le noyau de $\alpha_{G[p^n],\tau, n - \frac{p^{nf}-1}{p^f-1}^\mu\Ha(G)}$.
Le cran $\Fil_\tau(G[p^n])$ est de plus un cran de la $\tau$-filtration de Harder-Narasihman de $G[p^n]$, et donc la filtration canonique est compatible à la dualité, 
à la $p^k$-torsion ($k<n$) et aux quotients.
\ethr

\rem
L'hypothèse (\ref{hypfinale}) est simplement $^\mu\Ha(G) < \frac{1}{2p^{(n-1)f}}$ lorsque $p$ est assez grand.
\erem

\dem
L'existence des sous-groupes $\Fil_\tau(G[p^n])$ est assurée par le théorème (\ref{thrntors}), ainsi que les propositions sur les degrés, les hauteurs, l'application de Hodge-Tate et la 
$\tau$-filtration de Harder-Narsihman. Il ne reste donc plus qu'à montrer les inclusions, mais cela découle de la Proposition \ref{probij}, cf. \cite{Bij}, puisque l'on peut majorer,
\[ \frac{p^{nf} - 1}{p^f-1}\Ha_\tau(G) < \frac{2}{3},\]
et donc si $p_\tau < p_{\tau'}$ on vérifie que
\[ \deg\Fil_\tau(G[p^n]) + \deg\Fil_{\tau'}(G[p^n]) > \HN^{\mu-ord}(np_\tau) + \HN^{\mu-ord}(np_{\tau'}) - \frac{4}{3},\]
et donc la proposition assure que $\Fil_\tau(G[p^n]) \subset \Fil_{\tau'}(G[p^n])$. Dans le cas où $p_\tau = p_{\tau'}$ avec $\tau \neq \tau'$ il suffit de vérifier que,
\[ \deg\Fil_\tau(G[p^n]) + \deg\Fil_{\tau'}(G[p^n]) > 2\HN^{\mu-ord}(np_\tau) -2,\]
mais la même minoration que précédemment s'applique, et donc $\Fil_\tau(G[p^n]) =\Fil_{\tau'}(G[p^n])$.
\edem

\section{Application aux familles}
\label{sect9}
On fixe $p$ un nombre premier, $\mathcal O$ les entiers d'une extension non ramifiée de $\QQ_p$ et $(p_\tau,q_\tau)_\tau$ une signature telle que,
\[p > \max\{ \frac{2q_\tau}{1+K_\tau} : \tau \dans \mathcal I, q_\tau \neq h\} + 1.\]

\thr
\label{thrfam}
Soit $K$ une extension valuée complète de valuation discrete de $\QQ_p$ et $\mathfrak X$ un $\Spf(O_K)$-schéma formel topologiquement de type fini, sans $p$-torsion, et réduit.
Soit $G \fleche \mathfrak X$ un $\mathcal O$-module de Barsotti-Tate de signature $(p_\tau,q_\tau)_\tau$, tronqué d'échelon $r > \max_\tau k_\tau +n$.
Posons \[\eps_n = \frac{1}{p^{(n-1)f}}\min(\frac{1}{2},1+K_\tau - \frac{2q_\tau}{p-1}), \forall \tau \in \mathcal I \text{ tels que } q_\tau \notin \{0,h\}.\]
Soit $U = \mathfrak X^{rig}_{ord}(\overset{\circ}{\eps_n})$ le voisinage strict du lieu ordinaire de $\mathfrak X^{rig}$ où le $\mu$-invariant de Hasse est strictement plus petit 
que $\eps_n$. Alors il existe sur $U$ une filtration par des sous $\mathcal O$-modules $(\Fil_\tau G^{rig}[p^n])_\tau$ de $G^{rig}[p^n]_{|U}$ tel que,
\begin{enumerate}
\item $\Fil_\tau G^{rig}[p^n]$ est localement –pour la topologie étale– isomorphe à $(\mathcal O/p^n\mathcal O)^{p_\tau}$.
\item Si $p_\tau \leq p_{\tau'}$ on a une inclusion $\Fil_\tau G^{rig}[p^n] \subset \Fil_{\tau'}G^{rig}[p^n]$.
\item En tout point de $U$, la fibre de $\Fil_\tau G^{rig}[p^n]$ coïncide avec le cran de hauteur $np_\tau$ de la $\tau$-filtration de Harder-Narasihman de la fibre de $G^{rig}[p^n]$.
\end{enumerate}
 La filtration précédente est invariante sous $\End_\mathcal O(G)$, et si de plus $G$ est muni d'une polarisation compatible à $\mathcal O$, 
 $\lambda : G \overset{\simeq}{\fleche} G^D$, telle que $\lambda^D = \eps\lambda$, pour $\eps \dans \ZZ_p^\times (\mathcal O^\times ?)$, alors la filtration précédente vérifie en plus
 que chaque $\Fil_\tau G^{rig}[p^n]$ est totalement isotrope sous l'accouplement $G^{rig}[p^n] \times G^{rig}[p^n] \fleche \mathcal O/p^n\mathcal O(1)$ (ou $\ZZ_p$?).
 \ethr
 
 \dem
 On utilise le théorème \ref{thrfinO} et le théorème 4 de \cite{FarHN} qui nous permettent de mettre en famille chacun des groupes $\Fil_\tau(G[p^n])$ sur un éclatement formel 
 admissible $\mathfrak Y$ de $\mathfrak X$. Il faut prouver néanmoins que cela reste une filtration, puisqu'a priori on a plusieurs filtrations de Harder-Narasihman en jeu, 
mais on peut refaire comme dans la démonstration du théorème 4 de\cite{FarHN} : 
 Soit $\tau,\tau'$ tels que $p_\tau \leq p_{\tau'}$, alors considérons le morphisme de $\mathfrak Y$-schémas en groupes, 
 \[ \Fil_\tau G[p^n] \fleche G[p^n] \fleche G/\Fil_{\tau'}G[p^n].\]
 Il est nul en tout point de $\mathfrak Y^{rig} = \mathfrak X^{rig}$, mais comme $\Fil_\tau G[p^n]$ et $G/\Fil_{\tau'}G[p^n]$ sont localement libres sur $\mathfrak Y$ qui est réduit, 
 le morphisme est nul.
 Donc on a bien la filtration voulue.
 \edem

\subsection{Déformations de Frobénius}

Reprenons les notations du théorème principal \ref{thrfinO}.

\pro Soit $K/\QQ_p$ une extension finie, et $G/\Spec(O_K)$ un $\mathcal O$-module $p$-divisible tronqué d'échelon $k+f$. 
Notons pour tout $\tau$,
\[ r_\tau = |\{\tau' \dans \mathcal I : q_{\tau'} \leq q_\tau\}|.\]
Supposons 
\begin{equation}
\label{hypfin}
\tag{$H_f$}
^\mu\Ha(G) < \frac{1}{p^{(f-1)f}}\min(\frac{1}{2},1+K_\tau - \frac{2q_\tau}{p-1}), \forall \tau \in \mathcal I \text{ tels que } q_\tau \notin \{0,h\}.
\end{equation}
Alors le sous groupe $K_1 = \sum_\tau \Fil_\tau (G[p^f])[p^{r_\tau}]$ déforme le noyau de $F^f$ de $G[p^f] \otimes \overline \FP$, c'est à dire que,
\[ K_1 \otimes_{O_K} \overline \FP = \Ker F^f.\]
\epro

\dem
Soit $\tau \dans \mathcal I$. Soit $u$ une variable et notons $(M,\phi)$ le module de Kisin de $H_\tau^1 \subset G[p]$ sur $k[[u]]$.
On note $\phi^\#$ le linerarisé de $\phi$, et on décompose $M = \bigoplus_\tau M_\tau$. D'après la théorie des diviseurs élémentaires, il existe une base
$(e_1,\dots,e_{p_\tau})$ de $M_\tau$ telle que $(u^{a_1}e_1,\dots,u^{a_{p_\tau}}e_{p_\tau})$ soit une base de $\phi^\#(M_{\sigma^{-1}\tau})$ où $0 \leq a_i \leq e$.
Alors,
\[ \deg_\tau H_\tau^1 = \deg_\tau (M,\phi) = \frac{1}{e}\sum_{i=1}^{p_\tau}a_i.\]
Or on sait que \[\deg_\tau H_\tau^1 > p_\tau - \Ha_\tau.\]
On en déduit donc que $a_i \geq e(1-\Ha_\tau)$ pour tout $i$. On a donc que $u^{e(1-\Ha_\tau)}$ divise $\phi^\#_\tau$.
Or on sait que pour tout $\tau'$ tel que $q_{\tau'} \leq q_\tau$ on a $\deg_{\tau'} H_\tau \geq p_\tau - \Ha_\tau(G)$, et donc $u^{e(1-\Ha_\tau)}$ divise $\phi^\#_{\tau'}$ 
pour tout $\tau'$ tel que $q_{\tau'} \leq q_\tau$. Donc si on regarde le module de Kisin $(\widetilde M,\widetilde \phi)$ de $H_\tau^f$, la matrice de 
$\widetilde \phi^f \pmod u$ – qui correspond au module de Dieudonné de $H_\tau^f \otimes \overline \FP$ – est divisible par $p^{r_\tau}$ et donc 
$H_\tau^f[p^{r_\tau}] \otimes \overline \FP \subset \Ker F^f$.
On en déduit donc que $K_1 \otimes \overline \FP \subset \Ker F^f$ (dans $G[p^f] \otimes \overline \FP$), mais ils ont même hauteur.
\edem

\rem
Le même résultat reste vrai sous l'hypothèse ($H_{nf}$) avec \[K_n = \sum_\tau \Fil_\tau (G[p^f])[p^{nr_\tau}] \subset G[p^{nf}],\]
qui déforme alors $\Ker F^{nf}$, au sens 
précédent.

La preuve laisse entendre que le résultat est probablement vrai seulement modulo $p^{1-\max_\tau{\Ha_\tau}}$, puisque chaque $\phi^\#_\tau$ est nul modulo $u^{e(1-\Ha_\tau)}$.
Malheureusement, il ne semble pas clair à ma connaissance qu'il soit possible de relier, pour $H/\Spec(O_K)$ un schéma en groupes, $H \otimes O_K/p^w$ ($w \leq 1$) avec 
$M \otimes W(k)[[u]]/u^{ew}$. Si tel est le cas, la démonstration précédente devrait s'adapter.

Le même résultat avec $K/\QQ_p$ une extension quelconque (e.g. $K = C$) est encore vrai. On devrait probablement pouvoir faire une preuve similaire en remplaçant $W(k)[[u]]$ par $A_{cris}$, malheureusement il semble qu'une théorie des modules de Breuil-Kisin sur $A_{cris}/\mathcal O_C$ comme présentée dans \cite{FarKis,Lau} ne concerne que les groupes $p$-divisibles, éventuellement tronqués, ce que ne sont pas les crans de la filtration canonique...
Néanmoins on peut obtenir le résultat par un argument de familles.
\erem

\pro
Le résultat précédent vaut encore pour toute extension valuée $K/\QQ_p$ (en particulier $K = C$).
\epro

\dem
Soit $w' > {^\mu}\Ha(G)$ qui vérifie encore $(H_f)$, et soit $\mathcal X_{w'}$ l'ouvert de l'espace rigide $X^{rig}$ (où $X$ est la présentation du champ $\mathcal{BT}_{k+f}$) 
sur lequel ${^\mu}\Ha < w'$, et soit $G[p^f]$ la $p^f$-torsion du groupe universel.
Comme $w'$ vérifie $(H_f)$, d'après le théorème \ref{thrfam}, il existe une filtration de $G[p^f]$ par des schémas en groupes finis et plat, et on a donc sur $X_{w'}$ un $\mathcal O$-module fini et plat,
\[K_1 = \sum_\tau \Fil_\tau (G[p^f])[p^{r_\tau}] \subset G[p^f].\]
De manière équivalente, il existe $\mathfrak X_{w'}$ un ouvert d'un éclatement formel admissible de $\mathfrak X$ sur lequel $K_1$ s'étend en un schéma en groupe fini et plat.
Le groupe sur $\Spec(O_K)$ de l'énoncé, avec donc ${^\mu}\Ha < w'$, définit un $O_K$-point de $\mathfrak X_{w'}$, et donc un $\overline{\FP}$-point $\overline x$ de 
$\mathfrak X_{w'}\otimes \overline{\FP}$. Or sur $\mathfrak X_{w'}\otimes \overline{\FP}$, on a deux schémas en groupes finis et plats, $K_1 \otimes \overline{\FP}$ et $\Ker F^f$, le 
noyau du Frobenius itéré $f$-fois de $G[p^f]\otimes \overline{\FP}$, dont on sait de plus qu'ils sont égaux sur la réduction à $\overline{\FP}$ des points de 
$\mathcal X_{w'}(\overline{\QQ_p})$. Or la réduction $\mathcal X_{w'}(\overline {\QQ_p}) = \mathfrak X_{w'}(\mathcal O_{\overline{\QQ_p}}) \fleche \mathfrak X_{w'}(\overline{\FP})$ 
est surjective, donc 
$(K_1)_{\overline x} = (\Ker F^f)_{\overline x}$.
\edem

\rem
On peut vérifier (déjà dans le cas $\mu$-ordinaire) que $K_n$ n'est pas un cran de la filtration de Harder-Naraihman de $G[p^{nf}]$.
\erem
  
\appendix
\section{}

\lem
\label{calculann1} Soit $p$ un nombre premier, et $n,f \dans \NN^*$. Alors on a l'égalité,
\[\frac{p^{(n-1)f}-1}{p^f-1}\frac{1}{2p^{n-1}f} + 2 \frac{p^{nf}-1}{p^f-1}\frac{1}{2p^{n-1}f} - \frac{1}{f} \leq 1\]
\elem

\dem
Cela revient à l'équation,
\[ p^{(n-1)f} - 1 + 2(p^{nf}-1) - \frac{f+1}{f}2p^{(n-1)f}(p^f-1) \leq 0,\]
donc,
\[ 2p^{nf}(\frac{f+1}{f}-1) - p^{(n-1)f}(1 + 2\frac{f+1}{f}) + 3 \geq 0, \]
c'est à dire,
\[ p^{(n-1)f}(\frac{2p^f - 3f-1}{f}) + 3 \geq 0,\]
mais comme $2p^f \geq 3f+1$, on a bien la majoration voulue.
\edem

\pro[Bijakowski, \cite{Bij} proposition 1.25]
\label{probij}
Soit $D,C \subset G[p^n]$ deux sous-$\mathcal O$-modules de $\mathcal O$-hauteurs respectives $d \leq c$.
Supposons que \[\deg D + \deg C > \sum_{\tau'} \left(\min(np_{\tau'},d) + \min(np_{\tau'},c)\right) - |\{\tau' :  d-1\leq np_{\tau'} \leq c\}|,\]
alors $D \subset C$.
\epro

\dem
Notons $h = \Ht_\mathcal O(D\cap C)$. La $\mathcal O$-hauteur de $D+C$ est alors $d + c - h \geq h$.
On a alors,
\[ \deg(D+C) \leq \sum_{\tau'} \min(np_{\tau'},d + c - h), \quad \text{et}\quad \deg(D\cap C) \leq  \sum_{\tau'} \min(np_{\tau'},h).\]
On peut alors écrire,
\[ \deg D + \deg C \leq \deg(D+C) + \deg(D \cap C) \leq \sum_{\tau' \in A} 2np_{\tau'} + \sum_{\tau'\in B} \left(np_{\tau'} + h\right) + \sum_{\tau'\in C} \left(c + d\right),\]
où $A = \{ \tau' : np_{\tau'} < h\}, B = \{ \tau' : h \leq np_{\tau'} < d + c - h\},$ et $C = \{ \tau' : np_{\tau'} \geq d + c - h\}$.
Mais si $D \not\subset C$, alors $h \leq d - 1$, et on en déduit donc,
\[ \deg D + \deg C \leq \sum_{\tau'} \left(\min(np_{\tau'},d) + \min(np_{\tau'},c)\right) - |\{\tau' : d-1 \leq np_{\tau'} \leq c \}|.\]
Ce qui contredit l'hypothèse de l'énoncé.
\edem

 \nocite{*}
\bibliographystyle{alpha-fr} 
\bibliography{biblio} 

\backmatter

\end{document}